\documentclass[leqno]{amsart}

\usepackage{amssymb}
\usepackage{amsmath}
\usepackage{color}

\newtheorem{theorem}{Theorem}[section]
\newtheorem{lemma}[theorem]{Lemma}

\newtheorem{corollary}[theorem]{Corollary}

\newtheorem{proposition}[theorem]{Proposition}
\newtheorem{remark}[theorem]{Remark}
\newcommand{\R}{\mathbb{R}}

\newcommand{\q}{{\mathcal{C}}}
\newcommand{\N}{{\mathbb{N}}}
\newcommand{\Z}{{\mathbb{Z}}}
\renewcommand{\j}{{\mathcal{J}}}
\newcommand{\eps}{\varepsilon}

\newcommand{\K}{{\mathcal{K}}}

\newcommand{\e}{{\mathcal{E}}}

\renewcommand{\L}{{\mathcal{L}}}

\newcommand{\dem}[1]{\medskip \noindent {\bf #1}}
\newcommand{\fdem}{\hfill$\square$\medskip}

\definecolor{colorjer}{rgb}{0,0,0}
\definecolor{colorjua}{rgb}{0,0,0}

\long\def\co#1{\marginpar{\raggedright\small$\bullet$\ #1}}
\long\def\change#1{{\color{colorjua}#1}}
\long\def\changejer#1{{\color{colorjer}#1}}

\long\def\change#1{#1}
\long\def\changejer#1{#1}
\long\def\co#1{}
\long\def\changejer#1{{\color{colorjua}#1}}

\title{Pulsating \change{fronts} for nonlocal dispersion and KPP nonlinearity}
\author{J\'er\^ome Coville}
\address{\noindent J. Coville -- INRA, Equipe BIOSP, Centre de Recherche d'Avignon, Domaine Saint
Paul, Site Agroparc, 84914 Avignon cedex 9, France}
\email{jerome.coville@avignon.inra.fr}

\author{Juan D\'avila}
\address{\noindent J. D\'avila -- Departamento de
Ingenier\'{\i}a  Matem\'atica and CMM, Universidad de Chile, Casilla
170 Correo 3, Santiago, Chile.}
\email{jdavila@dim.uchile.cl}

\author{Salom\'e Mart\'\i nez}
\address{\noindent S. Mart\'\i nez -- Departamento de
Ingenier\'{\i}a  Matem\'atica and CMM, Universidad de Chile, Casilla
170 Correo 3, Santiago, Chile.}
\email{samartin@dim.uchile.cl}

\date{\today}

\numberwithin{equation}{section}

\begin{document}

\begin{abstract}
In this paper we are interested in propagation phenomena for nonlocal reaction-diffusion equations of the type:
$$\frac{\partial u}{\partial t}=J * u - u + f(x,u) \qquad t\in \R, \; x \in \R^N,$$
where $J$ is a probability density  and $f$ is a KPP nonlinearity  periodic in the $x$-variables.
Under suitable assumptions we establish the existence of pulsating fronts describing the invasion of the  0 state by an heterogeneous state. We also give a variational characterization of the minimal speed of such pulsating fronts and exponential bounds on the asymptotic behaviour of the solution.
\end{abstract}

\maketitle


\bigskip



\section{Introduction}
In this paper we are interested in propagation phenomena for nonlocal reaction-diffusion equations of the type:
\begin{equation}
\frac{\partial u}{\partial t}=J * u - u + f(x,u) \qquad t\in \R, \; x \in \R^N,\label{ptw.eq.para-intro}
\end{equation}
where $J$ is a probability density  and $f$ is a  nonlinearity which is KPP in $u$ and periodic in the $x$-variables, that is,
$$
f(x,u) = f(x+k,u)
\quad \forall x \in \R^N, \ k \in \Z^N ,  \ u \in\R.
$$
More precisely, we are interested in  the existence/non-existence and the characterization of front type solutions called pulsating fronts.
A pulsating front connecting 2 stationary periodic solutions $p_0, p_1$  of \eqref{ptw.eq.para-intro} is an entire solution
that has the form $u(x,t):=\psi(x\cdot e+ct,x)$ where $e$  is a unit vector in $\R^N$, $c\in\R$,
and $\psi(s,x)$ is periodic in the $x$ variable, and such that
 \begin{align*}
&\lim_{s\to -\infty}\psi(s ,x )=p_0(x) \quad \text{uniformly in } x\\
 &\lim_{s\to +\infty}\psi(s ,x )=p_1(x) \quad \text{uniformly in } x .
\end{align*}
\change{The real number $c$ is called the effective speed of the pulsating front.}

Using an equivalent definition, pulsating fronts were first defined and used by Shigesada, Kawasaki and Teramoto \cite{SKT1,SKT2} in their study of biological  invasions  in a heterogeneous environment  modelled by the following reaction diffusion equation
\begin{equation}
\frac{\partial u}{\partial t}=\nabla\cdot\left(A(x)\nabla u\right)+f(x,u)\quad \text{ in } \R^+\times \R^N,\label{ptw.eq.skt-equation}
 \end{equation}
where $A(x)$ and $f(x,u)$ are respectively  a periodic smooth elliptic matrix and  a smooth periodic function. Using heuristics and numerical simulations, in a one dimensional situation and for the particular nonlinearity $f(x,u):=u(\eta(x)-\mu u)$, Shigesada, Kawasaki and Teramoto were able to recover  earlier results on the minimal speed of spreading obtained by probabilistic methods by  Freidlin and Gartner \cite{Freidlin, FreidlinGartner}.


The above definition of pulsating front has been introduced by Xin \cite{Xin1,Xin2} in his study of flame propagation. This definition is a natural extension of the definition of the sheared travelling fronts  studied for example in \cite{berestycki-larroutourou-lions,berestycki-nirenberg}.  Within this framework, Xin  \cite{Xin1,Xin2} has proved  existence and uniqueness up to translation of pulsating fronts
for equation \eqref{ptw.eq.skt-equation} with a  homogeneous bistable or ignition non-linearity. Since then, much  attention  has been drawn to the study of  periodic reaction-diffusion
equations and the existence and the uniqueness  of  pulsating front have been proved  in various situations, see  for example \cite{berestycki-hamel-cpam,Berestycki-Hamel-Roques-I,Berestycki-Hamel-Roques-II,Hamel-Roques,Heinze,HPS,HZ,MNL,W2,Xin1,Xin2,Xin3}. In particular, Berestycki,  Hamel and Roques \cite{Berestycki-Hamel-Roques-I,Berestycki-Hamel-Roques-II} have showed that when $f(x,u)$ is  of KPP type, then the existence of a unique non trivial stationary solution $p(x)$ to \eqref{ptw.eq.skt-equation} is governed by the sign of the periodic principal eigenvalue of the following spectral problem
\begin{equation*}
\nabla\cdot\left(A(x)\nabla \phi\right)+f_u(x,0)\phi +\lambda_p \phi =0.
\end{equation*}
Furthermore, they have showed that there exists a critical speed $c^*$ so that a pulsating front with speed $c\ge c^*$ in the direction $e$ connecting the two equilibria 0 and $p(x)$ exists and no pulsating front with speed $c<c^*$ exists.  They also gave a precise characterisation of $c^*$ in terms of some periodic principal eigenvalue.
Versions of \eqref{ptw.eq.skt-equation} with periodicity in time, or more general media are studied in
\change{\cite{berestycki-hamel-cpam,berestycki-hamel-contempmath,Berestycki-Hamel-Nadin, Mellet,nadin,nadin-rossi,nolen-roquejoffre,Nolen-Ryzhik,shen-die, Zlatos}}.
It is worth noticing that when the matrix $A$ and $f$ are homogeneous, then the equation \eqref{ptw.eq.skt-equation}
\change{reduces} to a classical reaction diffusion equation with constant coefficients and the pulsating front $(\psi,c)$ is indeed a travelling front which have been well studied since the pioneering works of  Kolmogorov, Petrovskii and Piskunov \cite{KPP}.

 \smallskip

Here we are concerned with a nonlocal version of  \eqref{ptw.eq.skt-equation} where the classical local diffusion operator $\nabla\cdot(A(x)\nabla u)$ is replaced by the integral operator $J* u-u$.
The introduction of such type of long  range interaction finds its justification in many problems  ranging from micro-magnetism \cite{DGP,DOPT1,DOPT2}, neural network \cite{EMc} to ecology \cite{CMS,Clark,DK,KM,M,SSN}. For example, in some population dynamic models, such long range interaction is used to model the dispersal of individuals through their environment, \cite{F1,F2,HMMV}.
Regarding equation \eqref{ptw.eq.para-intro} we quote \cite{AB,BFRW,Chen,Cov2,Cov4, CD1,CDM2} for the existence and characterisation of travelling fronts for this equation with homogenous nonlinearity and \cite{bates-zhao,Cov6,cdm1,GR1,HMMV} for the study of the stationary problem.

In what follows,  we assume that $J:\R^N \to \R$ satisfies
\begin{align}
\label{hyp J}
\left\{
\begin{aligned}
&\text{$J \ge 0$, $\int_{\R^N} J =1$, $J(0)>0$,}
\\
&\text{$J$ is smooth, symmetric with support contained in the unit ball,}
\end{aligned}
\right.
\end{align}
and that $f:\R^N \times [0,\infty) \to \R$ is $[0,1]^N$-periodic in $x$ and satisfies:
\begin{align}
\label{hyp f1}
\left\{
\begin{aligned}
& \hbox{$f \in C^3(\R^N\times[0,\infty))$},
\\
& \hbox{$f(\cdot,0)\equiv0$},
\\
& \hbox{$f(x,u)/u$ is decreasing with
respect to $u$ on $(0,+\infty)$},
\\
& \hbox{there exists $M>0$ such that $f(x,u)\le0$ for all $u \ge
M$ and all $x$.}
\end{aligned}
\right.
\end{align}
The model example is
$$
f(x,u) = u ( a(x) - u )
$$
where $a(x)$ is a periodic, $C^3$  function.


Before constructing pulsating fronts, we discuss the existence
of solutions of the  stationary equation
\begin{align}
\label{stationary}
J * u - u + f(x,u) =0 \qquad x \in \R^N  .
\end{align}

Under the assumption \eqref{hyp f1},   0 is  a solution of \eqref{stationary} and, as shown in \cite{Cov6},
the existence of a positive periodic stationary  solution $p(x)$ is characterized by the sign of a {\em generalized principal eigenvalue}
of the linearisation of \eqref{stationary} around 0, defined by
\begin{align}
\label{def mu p}
\mu_0 = \sup\{\ \mu\in \R \ | \ \exists \phi \in C_{per}(\R^N), \phi>0,  \ \hbox{such that} \ J*\phi - \phi + f_u(x,0) \phi+\mu \phi\leq 0\}
\end{align}
where $C_{per}(\R^N)$ is the space of continuous periodic functions in $\R^N$.

More precisely, we have
\begin{theorem}
\label{thm 1}
The stationary equation \eqref{stationary} has a positive continuous periodic solution $p(x)$ if and only if $\mu_0 < 0$. Moreover the positive solution is Lipschitz and  unique in the class of positive bounded periodic function.
\end{theorem}
\co{ \changejer{I'm wondering what to add in the theorem. Here, We get  the Lipschitz regularity of the solution even if there exist not an eigenfunction and by the Theorem also uniqueness in the class of continuous function }}
This result is analogous to the characterization of stationary positive solutions of \change{the differential equation} \eqref{ptw.eq.skt-equation} with $f$ of type KPP in $u$.
The main difference is that $\mu_0$ is not always an eigenvalue, that is, the supremum in \eqref{def mu p} is not always achieved.
Similar results for \eqref{stationary}, but assuming that  $\mu_0$ is an eigenvalue \changejer{and for the one-dimensional case (i.e $N=1$) }, have been obtained in \cite{bates-zhao,cdm1}. \changejer{In this particular situation, the  uniqueness of the positive solution of  \eqref{stationary} in the class of \change{bounded measurable  functions} has been proved in \cite{cdm1}.  For the multidimensional case, the existence and uniqueness of a stationary solution in the class of periodic functions has been obtained by Shen and Zhang \cite{shen-zhang}  assuming that $\mu_0$ is eigenvalue and by Coville \cite{Cov6} without this assumption.}
\change{The difference of Theorem~\ref{thm 1} and \cite{Cov6} is that we obtain a Lipschitz continuous solution.}

The question whether $\mu_0$ is really a principal eigenvalue, that is, if there exists $\phi \in C_{per}(\R^N)$, $\phi>0$ such that
\begin{align}
\label{spectral-prob}
J * \phi - \phi + f_u(x,0)\phi+\mu_0\phi =0  \qquad \text{ in } \R^N
\end{align}
 has been studied in \cite{Cov6,shen-zhang} where simple criteria on $f_u(x,0)$ have been  derived to ensure the existence of a principal eigenfunction $\phi$.  For instance, \changejer{the following criterion proposed in \cite{Cov6}}
\begin{align*}
\int_{[0,1]^N} \frac{1}{A - f_u(x,0)} \, d x  = +\infty ,
\quad\text{ where }\quad
A = \max_{x\in \R^N} f_u (x,0) ,
\end{align*}
guarantees that $\mu_0$ is a principal eigenvalue. Some properties of $\mu_0$ and the   existence criteria will be discussed in Section~\ref{sec p ev}.

%

Our main result on pulsating fronts is the following:
\begin{theorem}
\label{thm main}
Assume $\mu_0<0$ and that there exists $\phi \in C_{per}(\R^N)$, $\phi>0$ satisfying \eqref{spectral-prob}.
Then, given any unit vector $e \in \R^N$
there is a number $\change{c_e^*}>0$ such that for $c\ge \change{c_e^*}$ \eqref{ptw.eq.para-intro} has a pulsating front solution $u(x,t)=\psi(x.e+ct,x)$ with effective speed $c$, and for $c<\change{c_e^*}$ there is no such solution.
\end{theorem}

The minimal speed $\change{c_e^*}$ is given by
\begin{equation}
\label{def c star}
\change{c^*_e}
:=\inf_{\lambda >0}\left(\frac{-\mu_\lambda}{\lambda}\right)
 \end{equation}
 where $\mu_\lambda$ is the periodic principal eigenvalue of the following problem
  \begin{equation}
J_\lambda * \phi -\phi +f_u(x,0)\phi +\mu \phi =0    \qquad \text{ in } \R^N
\end{equation}
with $J_\lambda(x):=J(x)e^{\lambda x.e}$.
We will see \change{in} Section~\ref{sec p ev} that this eigenvalue problem is solvable under the assumptions of Theorem~\ref{thm main}.

Shen and Zhang showed in \cite{shen-zhang}  that $\change{c_e^*}$ corresponds to the  speed of spreading for this equation in the following sense. For reasonable initial conditions, the solution of \eqref{ptw.eq.para-intro} satisfies
$$
\limsup_{t\to+\infty} \sup_{x\cdot e + c t \le 0} u(x,t) = 0
\qquad\text{if $c>\change{c_e^*}$,}
$$
while
$$
\liminf_{t\to+\infty} \inf_{x\cdot e + c t \ge 0} ( u(x,t) -p(x) ) = 0
\qquad\text{if $c<\change{c_e^*}$.}
$$
The nonexistence statement in  Theorem~\ref{thm main}
is a  consequence of the these spreading speed results.
Along our analysis, we also obtain some asymptotic behaviour of $\psi(s,x)$ as $s\to \pm\infty$ where $\psi$  is the pulsating front  constructed in Theorem~\ref{thm main}. More precisely, let $\lambda(c)$ denote the smallest positive $\lambda$ such that $c=\frac{-\mu_\lambda}{\lambda}$.
\begin{theorem}
\label{thm exp decay}
Assume $\mu_0<0$ and that
there exists $\phi \in C_{per}(\R^N)$, $\phi>0$ satisfying \eqref{spectral-prob}. Then, given any unit vector $e \in \R^N$ and $c\ge \change{c_e^*}$ we have

a) For any positive $\lambda$ so that $\lambda<\lambda(c)$ there exist  $C>0$  such that
\begin{align*}
\psi(s,x)  \le C e^{\lambda s}  \quad  \forall x \in \R^N, \
\forall s\in\R.
\end{align*}

b) There is $\sigma,C>0$ such that
$$
0\le p(x)- \psi(s,x) \le C e^{-\sigma s}
\quad \forall x\in\R^N , \ \forall s\ge0.
$$
\end{theorem}

Equation \eqref{ptw.eq.para-intro} can be related to a class of problems studied by Weinberger in \cite{W2}.
However, as observed in \cite{CD1,shen-zhang}, one of the  main difficulties in dealing with  the nonlocal equation  \eqref{ptw.eq.para-intro} comes  from the lack of regularizing effect of  \eqref{ptw.eq.para-intro}, which makes the framework developed by Weinberger not applicable,
since the compactness assumption required in \cite{W2} does not hold.

Another difficulty in the construction of pulsating fronts is  that the equation satisfied by the function $\psi$ (see \eqref{equ psi} below) involves an integral operator in time and space, which is in some sense degenerate. This  difficulty also appears in the classical reaction diffusion case, and it becomes delicate to proceed using the standard approaches used in  \cite{berestycki-larroutourou-lions,berestycki-nirenberg,KPP}.

Finally, we comment on some of the hypotheses made in the construction.
Regarding smoothness of the data, one can deal with less regularity of $J$ and $f$, but some arguments would have to be modified.
The hypothesis on the support of $J$ in \eqref{hyp J} can be weakened. For example, we believe that  the same results are true  assuming that $J$ satisfies the so called Mollison condition:
$$
\forall\, \lambda>0, \quad \int_{\R^N} J(z) e^{ \lambda|z|} \, d z < +\infty .
$$
Finally, the hypothesis that $\mu_0$ is an eigenvalue seems crucial in our approach. It is an interesting open problem to understand  whether some type of pulsating front exists in the case where $\mu_0$ is not an eigenvalue. We believe that if such solutions exist, they will be qualitatively different from the ones constructed in Theorem~\ref{thm main}.
\change{See also Remark~\ref{rem:pev} for other observations on this hypothesis.}

\changejer{In the preparation of this work, we were informed of a very recent work of Shen and Zhang \cite{shen-zhang2} done independently dealing with the existence and properties of pulsating front for a nonlocal equation like \eqref{ptw.eq.para-intro}.  The construction  of pulsating front proposed by Shen and Zhang  relies on a completely  different method and another definition of pulsating front.  With their method, they are able to construct bounded measurable pulsating fronts for any speed $c>c*(e)$ but fail to construct pulsating front  for the critical speeds $c^*(e)$ due to the lack of good Lipschitz regularity estimates on the fronts.  Some additional properties, such as exact exponential behaviour as $t\to-\infty$, uniqueness of the profile in a appropriate class and some kind of stability of the front are also studied in this work.
The main differences between the results obtained by Shen and Zhang and ours concern essentially the regularity of the fronts. Whereas they obtained bounded measurable front, we obtained  uniform Lipschitz front which is  a significant part of our work. We also have the feeling that our approach is more robust, in the sense that it does not strongly rely on the KPP structure and can be  adapted to other situations such as a monostable or ignition nonlinearity which seems not be the case for the method used in \cite{shen-zhang2}. We have in mind a problem like
\begin{equation*}
\frac{\partial u}{\partial t}=\int_{\R^N}J \left(\frac{x-y}{g(x)g(y)}\right)[ u(y) - u(x)]\,dy + f(u) \qquad t\in \R, \; x \in \R^N,
\end{equation*}
   where $f$ is monostable nonlinearity, $J$ a smooth probability density and $g$ a continuous positive periodic function.
It is worth noticing that in  \cite{shen-zhang2}, the existence of a principal eigenvalue for \eqref{spectral-prob} is also a crucial hypothesis.
}
\section{Scheme of the construction}
The proof of Theorem~\ref{thm 1} is contained in Section~\ref{sec stationary}, and follows by now standard arguments.

To construct a pulsating front solution $u$  of \eqref{ptw.eq.para-intro} in the direction $-e$ with effective speed $c$ connecting  $0$ and a positive periodic stationary solution $p$, we let
$\psi(s,x) = u \left( \frac{s-x\cdot e}{c} , x \right)$.
Then we need to find $\psi$ satisfying
\begin{align}
\label{equ psi}
\left\{
\begin{aligned}
& c \psi_s  =  M [\psi] - \psi + f(x,\psi) \quad \forall s\in \R,\ x\in\R^N
\\
& \psi(s,x + k) = \psi(s,x) \quad \forall s\in \R, \ x \in \R^N, \ k \in \Z^N,
\\
& \lim_{s\to-\infty} \psi(s,x) =0 \quad \text{uniformly in } x,\\
& \lim_{s\to\infty} \psi(s,x) =p(x) \quad \text{uniformly in } x,
\end{aligned}
\right.
\end{align}
where $M$ is the integral operator
$$
M [\psi](s,x) = \int_{\R^N} J(x-y) \psi ( s+(y-x)\cdot e ,y  ) \, d y .
$$


To analyse \eqref{equ psi} we introduce a regularized problem, namely, we consider for $\eps>0$
\begin{align}
\label{approx equ u}
c \psi_s =   M [\psi] - \psi + f(x,\psi) + \eps \Delta \psi  \quad \forall s\in \R, \  x\in \R^N
\end{align}
where $\Delta$ is the Laplacian with respect to the $x$ variables.
The stationary version of this equation is a perturbation of \eqref{stationary}:
\begin{align}
\label{approx stat}
0 = J * u - u + f(x,u) + \eps \Delta u \qquad x \in \R^N.
\end{align}
We will see in Section~\ref{sec stationary} that under the assumption that
\eqref{stationary}
has a positive periodic continuous solution $p$, for small $\eps >0$ the equation \eqref{approx stat} also has a stationary positive solution $p_{\eps} $ and $p_{\eps} \to p$ uniformly
as $\eps \to 0$.



As a step to prove Theorem~\ref{thm main}, for small $\eps>0$ we will find $\change{c_e^*(\eps)}$ such that for $c\ge \change{c_e^*(\eps)}$ there exists a solution $\psi_\eps$ to \eqref{approx equ u} satisfying
\begin{align}
\label{approx cond}
\left\{
\begin{aligned}
&
\lim_{s\to-\infty} \psi(s,x) =0
\\
&
\lim_{s\to+\infty} \psi(s,x) =p_{\eps}(x)
\\
& \hbox{$\psi(s,x)$ is increasing in $s$ and periodic in $x$,}
\end{aligned}
\right.
\end{align}
This is done in Section~\ref{s:approximate pulsating fronts}, following in part the methods developed in \cite{Berestycki-Hamel-Roques-II}.

A substantial part of this article is devoted to obtain  estimates for $\psi_\eps$
that will allow us to prove that $\psi  = \lim_{\eps\to0} \psi_\eps$ exists and solves \eqref{equ psi}.
These estimates are based on the expected exponential decay of $\psi$ as $s\to-\infty$, which we discuss next.
Suppose $\psi$ is a solution of \eqref{equ psi}.
One may expect that for some $\lambda>0$
$$
\psi(s,x)  = e^{\lambda s} w(x) + o(e^{\lambda s})
\quad \text{as } s \to - \infty, \quad x \in \R^N
$$
where $w$ is a positive periodic function, \change{at least when $c>c_e^*$.}
Then at main order the equation in \eqref{equ psi} yields
\begin{align}
\label{main order}
c \lambda w = \int_{\R^N} J(x-y) e^{\lambda(y-x).e} w(y) \, d y - w  + f_u(x,0) w \quad\text{in } \R^N.
\end{align}
%
Define
\begin{align*}
J_\lambda(x) = J(x) e^{-\lambda x\cdot e}
\end{align*}
then \eqref{main order} can be written as the periodic eigenvalue problem
\begin{align}
\label{ev}
\left\{
\begin{aligned}
J_\lambda *w - w  + f_u(x,0) w + \mu_\lambda w =0 \quad \hbox{in } \R^N
\\
w>0 \hbox{ is continuous and periodic,}
\end{aligned}
\right.
\end{align}
which will be studied in Section~\ref{sec p ev}. In particular, under the assumptions of Theorem~\ref{thm main},  we will see that it has a principal eigenvalue $\mu_\lambda$ in the space of continuous periodic functions.
Then the speed of the travelling front should be given by
$c = -\frac{\mu_\lambda}{\lambda}$,
and this leads to the formula for the minimal speed \eqref{def c star}.

For the solutions of \eqref{approx equ u} and \eqref{approx cond} on can guess a similar asymptotic behaviour as $s\to-\infty$ and a formula for the minimal speed
\begin{align}
\label{c star eps}
\change{c_e^*(\eps)} = \min_{\lambda>0}  ( - \frac{\mu_{\eps,\lambda}}{\lambda} )
\end{align}
where $\mu_{\eps,\lambda}$ is the principal eigenvalue of $-L_{\eps,\lambda}$ where
$$
L_{\eps,\lambda} w=  \eps \Delta w + J_\lambda w - w  + f_u(x,0) w
$$
in the space of $C^2$ periodic functions.

Based on the estimates developed in Section~\ref{sec estimates L eps} for the operator $ L_{\eps,\lambda} u $, we prove in Section~\ref{sec Exponential bounds} exponential bounds of the form:
for  $0<\lambda<\lambda_\eps(c)$
\begin{align}
\label{exp est psi eps}
\psi_\eps(s,x)  \le C e^{\lambda s}  \quad  \forall x \in \R^N \quad \forall s\in\R
\end{align}
where
$\lambda_\eps(c)$ is  the smallest positive $\lambda$ such that $c = -\frac{\mu_{\eps,\lambda}}{\lambda} $, and
$C$ does not depend on $\eps>0$.
This exponential bound is obtained by studying the two sided Laplace transform of $\psi_\eps$, an idea present in \cite{carr-chmaj}.

The exponential estimate \eqref{exp est psi eps} allows us in Section~\ref{sec proof main thm} to obtain uniform control of local Sobolev norms $\|\psi_\eps\|_{W^{1,p}}$ with $p>N$, which in turn implies that we obtain a locally uniform limit $\psi = \lim_{\eps\to0} \psi_\eps$ for some subsequence. The final step is to verify that $\psi$ satisfies all the requirements in \eqref{equ psi}.

%
\section{Principal eigenvalue for non-local operators}
\label{sec p ev}
Let us recall the notation
$$
C_{per}(\R^N) =
\{\phi\in C(\R^N) \ | \ \phi \ \hbox{is} \ [0,1]^N- \ \hbox{periodic}\}.
$$
For the rest of the article it is crucial to understand the eigenvalue problem
\eqref{ev}, and the purpose of this section is to study its properties.  We will write
\eqref{ev} in the form
\begin{align}
\label{ev2}
\left\{
\begin{aligned}
& L_\lambda \phi + \mu \phi =0 \quad \hbox{in } \R^N
\\
& \phi \in C_{per}(\R^N), \ \phi>0
\end{aligned}
\right.
\end{align}
where
$$
L_\lambda w = J_\lambda* w + a(x) w
$$
and $a(x)=f_u(x,0)-1  \in C_{per}(\R^N)$.

We say that $L_\lambda$ has a principal eigenfunction if for some $\mu\in \R$ there is a solution in $C_{per}(\R^N)$ of \eqref{ev2}.

As we will see later, it is not true in general that $L_\lambda$ has a principal eigenfunction, but it is convenient to define in all cases
\begin{equation}
\label{defeig}
\mu_\lambda=\sup\{\ \mu\in \R \ | \ \exists \phi \in C_{per}(\R^N), \phi>0,  \ \hbox{such that} \ L_\lambda \phi+\mu \phi\leq 0\}
\end{equation}
and call it the generalized principal eigenvalue of $-L_\lambda$. The name is motivated by the following result.

\begin{proposition}\label{lemma-charact}
Let $\lambda\in \R$. If there is $\mu \in \R$, $\phi \in C_{per}(\R^N)$, $\phi\ge 0$ and nontrivial satisfying $L_\lambda \phi + \mu \phi = 0$, then $\mu$ is given by \eqref{defeig} and it is simple eigenvalue of $L_\lambda$.
\end{proposition}

The proof of this is a direct adaptation of Lemma 3.2 in \cite{Cov6}.

The next proposition characterizes the existence of a principal eigenfunction.
\begin{proposition}\label{lemma1-aux}
\change{If $a\in C_{per}(\R^N)$, then
$\max a(x)+\mu_\lambda \le 0$. Moreover,   $\max a(x)+\mu_\lambda<0$} if and only if $L_\lambda$ admits a principal eigenfunction.
\end{proposition}
\change{For the proof of the above result and the following two (Proposition~\ref{propconvexity} and Corollary~\ref{corollary}) see later in this section.}

\begin{proposition}
\label{propconvexity}
The function $-\mu_\lambda$ is  convex in $\R$ and even. In particular,  \change{$-\mu_\lambda$} is nondecreasing in $[0,\infty)$
and nonincreasing in $(-\infty,0]$.
\end{proposition}

\begin{corollary}\label{corollary}
If $L_0$ has a principal eigenfunction then for all $\lambda \in \R$, $L_\lambda$ has a principal eigenfunction.
\end{corollary}

In general it is difficult to describe precisely in terms of $J$ and $a$ whether $L_\lambda$ has a principal eigenfunction, but we have sufficient and necessary conditions.

\begin{proposition}
\label{prop pri eig}
Assume $a\in C_{per}(\R^N)$ and let $A:=\max_{\R^N}a(x)$. There are constants $C_1, C_2,m>0$
that depend on $J_\lambda$ such that:

a)   if
\begin{align}
\label{integral grande}
\int_{[0,1]^N}\frac{1}{A - a(x)} \, d x \ge C_1 \|a\|_{L^\infty}^m
\end{align}
then $L_\lambda$ admits a principal eigenfunction,

b) if
\begin{align*}
\int_{[0,1]^N}\frac{1}{A - a(x)} \, d x \le C_2
\end{align*}
then $L_\lambda$ has no principal eigenfunction.
\end{proposition}
\change{We give the proof of this Proposition later on inside this section.}

Finally, we need the next proposition to show that the formula \eqref{def c star} is well defined and gives a positive number.
\begin{proposition}
\label{prop propiedades de mu}
The function $\lambda \mapsto \mu_\lambda$ is continuous and for all $\eps>0$ there exists $\sigma>0$ such that
$$
-\mu_\lambda \ge -\mu_0-\eps+ \sigma e^{\sigma|\lambda| } \quad \forall \lambda\in \R.
$$
\end{proposition}
\change{The above Proposition is proved later on inside this section.}

\begin{remark}
Many of the previous results have appeared in similar contexts, or have been proved under slightly different conditions.
Existence of a principal eigenfunction was obtained for symmetric non-local operators  in \cite{HMMV}, and later also in \cite{bates-zhao,Cov6,cdm1,shen-zhang}.
A condition like \eqref{integral grande} is always explicitly or implicitly assumed in these works.
The motivation for  definition \eqref{defeig} is taken from \cite{berestycki-nirenberg-varadhan}. It has been adapted to many elliptic operators, and was first introduced for non-local operators in \cite{Cov6}. In this work the author obtained many of the results described here for an integral operator on a domain in $\R^N$.
A characterization like Proposition~\ref{lemma1-aux}  for $\mu_\lambda$ was first obtained in \cite{Cov6}.
The convexity of $-\mu_\lambda$, Proposition~\ref{propconvexity}, is proved in \cite{shen-zhang} under the assumption that a principal eigenfunction exists. Examples of non-local operators with no principal eigenvalue are also presented in \cite{Cov6,shen-zhang}.

\end{remark}

\change{The rest of this section is devoted to prove Propositions~\ref{lemma1-aux}, \ref{propconvexity}, Corollary~\ref{corollary}, and Propositions~\ref{prop pri eig} and \ref{prop propiedades de mu}.}
We start with some basic facts about the definition \eqref{defeig}.
The following results are simple adaptations from results found in \cite{Cov6}.
\begin{proposition}
\label{propaux1}
(Proposition 1.1 \cite{Cov6})
Given $a\in C_{per}(\R^N)$, and $J:\R^N \to \R$, $J\ge 0$ in $L^1(\R^N)$ define
$$
\mu_p(J,a)=\sup\{\mu\in \R | \exists \phi \in C_{per}(\R^N), \phi>0,  \ \hbox{such that} \ J*\phi + a\phi +\mu \phi\leq 0\},
$$

Then the following hold:
\begin{enumerate}
\item[(i)] If $a_1\ge a_2$, then
$$
\mu_p (J,a_2) \ge \mu_p(J,a_1) .
$$
\item[(ii)]
If $J_1 \ge J_2$ then
$$
\mu_p (J_2,a) \ge \mu_p(J_1,a) .
$$

\item[(iii)] $\mu_p(J,a)$ is Lipschitz  in $a$, more precisely
$$
|\mu_p(J,a_1)-\mu_p(J,a_2)|\leq \|a_1-a_2\|_{\infty} .
$$
\end{enumerate}
\end{proposition}

To prove Proposition~\ref{prop pri eig} we will need a generalization of the Krein-Rutman theorem \cite{krein-rutman} for positive not necessarily compact operators due to Edmunds, Potter and Stuart \cite{EPS}. For this we recall some definitions. A cone in a real Banach space $X$ is a non-empty closed set $K$
such that for all $x, y \in K$ and all $\alpha \ge 0$ one has $x +
\alpha y \in K$, and if $x\in K$, $-x \in K$ then $x=0$. A cone
$K$ is called reproducing if $X = K-K$. A cone $K$ induces a
partial ordering in $X$ by the relation $x\le y$ if and only if
$x-y \in K$. A linear map or operator $T:X\to X$ is called
positive if $T(K)\subseteq K$.

If $T:X\to X$ is a bounded linear map on a complex Banach space X,
its essential spectrum (according to Browder \cite{browder})
consists of those $\lambda$ in the spectrum of $T$ such that at
least one of the following conditions holds : (1) the range of
$\lambda I - T$ is not closed, (2) $\lambda$ is a limit point of
the spectrum of $T$, (3) $\cup_{n=1}^\infty ker(\lambda I -
T)^n$ is infinite dimensional. The radius of the essential
spectrum of $T$, denoted by $r_e(T)$, is the largest value of
$|\lambda|$ with $\lambda$ in the essential spectrum of $T$. For
more properties of $r_e(T)$ see \cite{nussbaum}.

\begin{theorem}
(Edmunds, Potter, Stuart \cite{EPS})
\label{cdm.th.eps}
Let K be a reproducing cone in a real Banach space X, and let
$T\in\L(X)$ be a positive operator such that $T^m(u)\ge cu$ for
some $u\in K$ with $\|u\|=1$, some positive integer $m$ and some
positive number $c$. If $c^{1/m}>r_e(T)$, then $T$ has an
eigenvector $v\in K$ with associated  eigenvalue $\rho\ge
c^{1/m}$.
and  $T^*$ has an  eigenvector $v^*\in K^*$ corresponding to the eigenvalue  $\rho$.
\end{theorem}

If the cone
$K$ has nonempty interior and $T$ is strongly positive, i.e.
$u\ge0$, $u\not=0$ implies $T u \in int(K)$, then $\rho$ is the
unique $\lambda \in \R$ for which there exist nontrivial $v \in K$
such that $T v = \lambda v$ and $\rho$ is simple, see
\cite{zeidler}.

\medskip
\noindent
{\bf Proof of Proposition~\ref{prop pri eig}.}
\

a)
Write
the eigenvalue problem \eqref{ev2}
in the form
\begin{align*}
J_\lambda *u + b(x) u = \nu u
\end{align*}
where
$$ b(x) = a(x) + k, \quad \nu = -\mu + k$$
and $k>0$ is a constant such that $\inf_{} b >0$.
Sometimes we will use the operator notation $J_\lambda[\phi]=J_\lambda * \phi$.
We study this eigenvalue problem in the space $C_{per}(\R^N)$ with uniform norm, where the operator $J_\lambda$ is compact.
Let $u \in C_{per}(\R^N)$, $u\ge 0$ and $m\in \N$   .
Since $u$ and $b$ are non-negative and $J_\lambda$ is a positive operator, we see that
\begin{align}
\label{ineq nu}
(J_\lambda+b(x))^{m}[u] \ge \change{J_\lambda^{m}[u]} +b(x)^{m}u
\end{align}
We observe that there are $m$ and $d>0$ depending on $J$ such that for $\change{u \in C_{per}(\R^N)}$, $u\ge 0$,
$$
\change{J_\lambda^{m}[u]} \ge d \int_{[0,1]^N} u .
$$
\co{proof of inequality}
\change{Indeed,
$$
J_\lambda^m[u] = J_\lambda^{(m)} *u ,
$$
where $J_\lambda^{(m)}$ denotes the $m$-fold convolution $J_\lambda*\ldots*J_\lambda$.
Let $B_R(x_0) $ with $R>0$  be such that $J_\lambda(x)>0$ for points $x\in B_R(x_0)$. Then $J_\lambda*J_\lambda(x) >0$ for $x \in B_{2R}(2x_0)$. Iterating this argument we get  $J_\lambda^{(m)}(x)>0 $ for $x \in B_{mR}(m x_0)$. We choose now $m$ large so that $B_{mR}(m x_0)$ contains some closed cube $Q$ with vertices in $\Z^N$. Let $d = \inf_{x\in Q} J_\lambda^{(m)}(x)>0$. Then, for $u \in C_{per}(\R^N)$, $u\ge 0$,
\begin{align*}
J_\lambda^m [u](x) & = \int_{\R^n}
J_\lambda^{(m)}(x-y) u(y) \, d y
\geq \int_{Q} J_\lambda^{(m)} (z ) u(x-z) \, d z
\\
& \geq d \int_{Q}  u(x-z) \, d z =\int_{[0,1]^N}  u ,
\end{align*}
since $u$ is $[0,1]^N$-periodic.
}

Let $\eps>0$ and define the continuous periodic positive function
$$
u_\eps(x) = \frac{1}{\max b^m - b(x)^m + \eps}  .
$$
We claim that choosing $\eps$ and $C_1$ in \eqref{integral grande} appropriately there is $\delta>0 $ such that
\begin{align}
\label{ineq integral}
J_\lambda^m u_\eps + b(x)^m u_\eps \ge (\max b  + \delta)^m u_\eps
\quad \text{in } \R^N.
\end{align}
Indeed, taking $C_1$ large in \eqref{integral grande} and then $\eps>0$ small, we have
$$
d \int_{[0,1]^N} \frac{1}{\max b^m - b(x)^m + \eps } \, d x >1.
$$
Then to prove \eqref{ineq integral} it is sufficient to have
$$
1 > \frac{(\max b + \delta)^m - b(x)^m}{\max b ^m - b(x)^m + \eps}
\quad\text{in } \R^N.
$$
This last condition holds provided we take $\delta$ sufficiently small.
Therefore, by \eqref{ineq nu} and \eqref{ineq integral} we have
$$
(J_\lambda+b(x))^{m}[u_\eps] \ge (\max b  + \delta)^m u_\eps.
$$
Using the compactness of the operator $J_\lambda$,   we have
$r_e(J_\lambda+b(x)) =\max_{x\in \R^N} b(x)$, and by Theorem
\ref{cdm.th.eps}
we obtain the desired conclusion.
We observe that the principal eigenvalue is simple since
the cone of positive periodic functions has non-empty
interior and, for a sufficiently large $p$, the operator $(J_\lambda +
b)^p $ is strongly positive.
Any point $\nu$ in the spectrum of  $(J_\lambda + b)$ with $|\nu| > r_e(J_\lambda + b)$ is isolated, see \cite{browder}. In particular  the principal eigenvalue is an isolated point in the spectrum.

\medskip
b)
As before, without loss of generality we can assume $a>0$.
Suppose there exists a principal periodic eigenfunction $\phi$ with eigenvalue $\mu$.  Then
$
\max a(x)  + \mu < 0
$.
Let $\q=[0,1]^N$ and note that
$$
J_\lambda * \phi(x) = \int_{\R^N} J(x-y) e^{\lambda (x-y)\cdot e} \phi (y) \, d y
= \int_\q \sum_{k \in \Z^N} J(x-z-k) e^{\lambda (x-z-k) \cdot e} \phi(z) \, d z
$$
$$
\le ( \int_\q \phi ) \, \sup_{x,z \in \q } \sum_{k \in \Z^N} J(x-z-k) e^{\lambda (x-z-k) \cdot e}
$$
But then
$$
 \phi (x)\le  \frac{1}{-(a(x)+\mu) } ( \int_\q \phi ) \, \sup_{x,z \in \q } \sum_{k \in \Z^N} J(x-z-k) e^{\lambda (x-z-k) \cdot e}
 .
$$
Integrating the above inequality we obtain
$$
\int_\q \phi \le  \int_\q \frac{1}{-(a(x)+\mu) } \, d x
\cdot \int_\q \phi \cdot\, \sup_{x,z \in \q } \sum_{k \in \Z^N} J(x-z-k) e^{\lambda (x-z-k) \cdot e},
$$
and hence
$$
1\le \int_\q\frac{1}{-(a(x)+\mu) } \, d x
\cdot\, \sup_{x,z \in \q } \sum_{k \in \Z^N} J(x-z-k) e^{\lambda (x-z-k) \cdot e} .
$$
Since $\mu\le -\max a(\cdot)$
$$
1\le \int_\q\frac{1}{\max a(\cdot) - a(x) } \, d x
\cdot\, \sup_{x,z \in \q } \sum_{k \in \Z^N} J(x-z-k) e^{\lambda (x-z-k) \cdot e}
$$
Let
$$
M=
\sup_{x,z \in \q } \sum_{k \in \Z^N} J(x-z-k) e^{\lambda (x-z-k) \cdot e} .
$$
If
$$
M \int_\q\frac{1}{\max a(\cdot) - a(x) } \, d x   <1
$$
there can not exist a principal eigenfunction.
\qed

\medskip
\noindent
{\bf Proof of Proposition~\ref{lemma1-aux}.}
From the definition we obtain directly
$   \max a(x)+\mu_\lambda\leq 0$ for all $\lambda\in \R$.
If there exists a principal eigenfunction $\phi \in C_{per}(\R^N)$, then clearly
$
\max a(x) + \mu_\lambda < 0
$.

Now suppose that $
\max a(x) + \mu_\lambda < 0
$.
We approximate $a$ by functions $a_\eps \in C_{per}(\R^N)$ such that $\max a=\max a_\eps$, $\|a-a_\eps\|_\infty \to 0$ as $\eps\to 0$, and
\begin{align}
\label{integral infinita}
\int_{[0,1]^N}\frac{1}{\max a_\eps - a_\eps(x)} \, d x =+\infty.
\end{align}
Then, by Proposition~\ref{prop pri eig}  there exists a positive, periodic $\phi_\eps$, with $\|\phi_\eps\|_\infty=1$,  such that
$$J_\lambda*\phi_\eps+(a_\eps(x)+\mu_\lambda^\eps)\phi_\eps=0, \ \hbox{in} \ \R^N.$$
Since by Proposition~\ref{propaux1}, $\mu_\lambda^\eps\to  \mu_\lambda$,  there exists $\delta>0$ such that
$a_\eps(x)+\mu_\lambda^\eps<-\delta$ for all $x$ and $\eps$. Therefore, by a simple compactness argument, we have that $\phi_\eps\to\phi$ uniformly as $\eps\to 0$, with $\phi$ positive satisfying \eqref{eig1}, which concludes the proof.

\qed

\begin{remark}
\label{remark smoothness eigenf}
If  $L_\lambda$ has a principal eigenfunction $\phi \in C_{per}(\R^N)$, and additionally $a \in C^k$, $k\ge 1$ and $J$ is $ C^k$, then $\phi$ is also $C^k$, which follows from
$$
J_\lambda \phi = (-\mu_\lambda - a) \phi
$$
and $ -\mu_\lambda - a \ge \delta$ for some $\delta>0$.
\end{remark}


\noindent{\bf Proof of Proposition~\ref{propconvexity}.}
To prove this result, we will first suppose that $a$ satisfies \eqref{integral infinita}, and then
we proceed by an approximation argument.
We will prove the convexity using an idea from \cite{shen-zhang}.
Let $\lambda_1,\lambda_2\in \R$, and $t\in (0,1)$. If $a$ satisfies \eqref{integral infinita} then by Proposition \ref{prop pri eig}   there exists
$\phi_1, \phi_2$ positive solutions of \eqref{ev2},  with corresponding eigenvalues $\mu_1, \mu_2$, for  $\lambda_1, \lambda_2$ respectively.
Consider $\phi=\phi_1^t\phi_2^{1-t}$. Then by H\"older's inequality we have that
$$J_\lambda*\phi \leq (J_{\lambda_1}*\phi_1)^{t}(J_{\lambda_2}*\phi_2)^{1-t}.$$
Using the inequality above and  that $\phi_1$ and $\phi_2$ are solutions of \eqref{ev2} we obtain that
$$
J_\lambda*\phi \leq ((-a(x)-\mu_1)\phi_1)^{t}((-a(x)-\mu_2)\phi_2)^{1-t}=(-a(x)-\mu_1)^{t}(-a(x)-\mu_2)^{1-t}\phi
$$
and then using Young's inequality we obtain that
$$J_\lambda* \phi\leq (t(-a(x)-\mu_1)+(1-t)(-a(x)-\mu_2))\phi=(-a(x)+t\mu_1+(1-t)\mu_2)\phi,$$
from where
$$\mu_{t\lambda_1+(1-t)\lambda_2}\ge t\mu_1+(1-t)\mu_2,$$
which gives the convexity.

To conclude when  \eqref{integral infinita} does not hold, we just approximate $a$ by $a_\eps$ satisfying \eqref{integral infinita}
and $a_\eps\to a$ uniformly in $\R^N$. Then the result follows by Proposition~ \ref{propaux1} (iii).

Finally, we claim that
the function $\mu_\lambda$ is even. Indeed, suppose first $\mu_\lambda$ is the principal eigenvalue of $L_\lambda$, so $\mu_\lambda + \max a(x) <0$. Considering $L_\lambda$ in the space of $L^2_{loc}(\R^N)$ periodic functions, we have that $L_{-\lambda}$ is its adjoint, and therefore $\mu_\lambda$ is in the spectrum of $L_{-\lambda}$. Using $\mu_\lambda + \max a(x) <0$ it is easy to see that $\mu_\lambda$ is the principal eigenvalue of $L_{-\lambda}$. In the case $L_\lambda$ has no principal eigenfunction, we directly deduce $\mu_{\lambda} = \mu_{-\lambda}$.

Since $-\mu_\lambda$ is even and convex, we obtain, that $\mu$ is nondecreasing in $(0,\infty)$ and noincreasing in $(-\infty,0)$.
$\Box$

\medskip\noindent
{\bf Proof of Proposition~\ref{prop propiedades de mu}.}
For the continuity of $\lambda \mapsto \mu_\lambda$ we argue as follows. Suppose first that  $a$ satisfies \eqref{integral infinita} and $\lambda_j \to \lambda_\infty$.
It is easy to see that $\mu_{\lambda_j}$ is bounded, so up to a subsequence $\mu_{\lambda_j} \to \mu$. Let $\phi_j \in C_{per}(\R^N)$ be the principal eigenfunction associated \change{with} $\mu_{\lambda_j}$ ($j=1,2,\ldots$) normalized so that $\|\phi_j\|_{L^\infty}=1$. Since $\mu + \max a < 0$, we have
$\mu_{\lambda_j} + \max a \le -\delta<0$ for some $\delta>0$ and all $j$ large. Then from
$$
J_{\lambda_j} * \phi_j = ( - \mu_{\lambda_j} - a) \phi_j
$$
we obtain compactness to say that for a subsequence $\phi_j$ converges uniformly to a nontrivial, nonnegative function $\phi \in C_{per}(\R^N)$ satisfying the eigenvalue problem
$$
J_{\lambda_\infty} * \phi = ( - \mu - a) \phi .
$$
Because of the uniqueness of the principal eigenvalue, Proposition~\ref{lemma-charact},  $\mu = \mu_{\lambda_\infty}$.

If  $a$ does not satisfy \eqref{integral infinita} we argue approximating $a$ by $a_\eps$ that satisfy \eqref{integral infinita}. Let $\mu_\lambda^\eps$ denote the principal eigenvalue of $-J_\lambda - a_\eps$.  We note that the convergence $\mu_\lambda^\eps \to \mu_\lambda$ as $\eps\to 0$ is uniform by Proposition~\ref{propaux1} (iii), so continuity of $\mu_\lambda^\eps$ with respect to $\lambda$ for all $\eps$ yields continuity of $\lambda \mapsto u_\lambda$.

Next we  show the exponential growth of $-\mu_\lambda$.
Observe that if $\phi\in C_{per}(\R^N)$ then
\begin{equation*}
J_\lambda *\phi=\int_{[0,1]^N} k_{\lambda}(x,y) e^{-\lambda(x-y) \cdot e} \phi(y)dy,
\end{equation*}
where
$$k_\lambda(x,y)=\sum_{k\in \Z^N} e^{\lambda k\cdot e} J(x-y-k).$$
The function $k_\lambda(\cdot, y)$ is $[0,1]^N$-periodic.
We consider the following eigenvalue problem
\begin{align*}
\hat L_\lambda \phi +( \mu+\eps) \phi = 0
\quad \hbox{with} \ \phi\in C\left([0,1]^N\right),
\end{align*}
where $\eps>0$ and
$$
\hat L_\lambda \phi=\int_{[0,1]^N} k_{\lambda}(x,y) e^{-\lambda(x-y)\cdot e} \phi(y)dy  + a(x) \phi + \mu_0  \phi.
$$

We will assume first that the support of $J$ is large, so that for some constants $b,d>0$:
$$
k_\lambda(x,y) \ge d e^{ b \lambda } \quad\forall x,y \in[0,1]^N.
$$
Let $w(y) = e^{-\lambda y \cdot e }$. Then
$$
\hat L_\lambda w \ge ( d e^{b\lambda} + a(x)+\mu_0 +\eps )  w \ge \delta e^{b\lambda} w
$$
where $\delta>0$ and where we take $\lambda$ large. If $\lambda>0$ is large enough, by Theorem~\ref{cdm.th.eps} we obtain a principal eigenfunction $\hat \phi \in C([0,1]^N)$ of $\hat L_\lambda$, with principal eigenvalue $-\hat \mu_\lambda \ge \delta e^{b\lambda}$.
Since
$k_\lambda(x,y) e^{\lambda(x-y)\cdot e }$ is periodic in $x$,
we see that $\hat \phi$ is periodic. Therefore, extending it periodically to $\R^N$, we find that it is the principal eigenfunction of $L_\lambda$ and $-\mu_\lambda +\mu_0+\eps= -\hat \mu_\lambda \ge \delta e^{b\lambda}$.
Now since $-\mu_\lambda$ is non decreasing in $\lambda$ we have $-\mu_\lambda +\mu_0+\eps\ge \eps$ and by taking $\delta $ smaller if necessary we achieve for all $\lambda$
$$-\mu_\lambda \ge -\mu_0-\eps+ \delta e^{b\lambda}.$$

Without the assumption that the support of $J$ is large, we can assume that $a(x)\ge 0$ and work with $m$ large so that the support of $J^m$ is large. Then
$$
(J_\lambda + a(x))^m \ge J_\lambda^m + a(x)^m .
$$
Notice that
$$
J_\lambda^m (x) = e^{\lambda x\cdot e } J^m(x)
$$
so the previous argument applies and we deduce that the principal eigenvalue of $J_\lambda^m + a(x)^m $ grows exponentially as $\lambda \to+\infty$. Then the same holds for $(J_\lambda + a(x))^m$ and therefore for $J_\lambda + a(x)$.
\qed

\change{
\begin{remark}
\label{rem:pev}
We would like to comment here on the hypothesis in Theorem~\ref{thm main} that there is a principal eigenvalue for problem \eqref{spectral-prob}. In fact, the proof of Theorem~\ref{thm main} reveals that we actually need only that \eqref{ev} has a principal eigenvalue for all $\lambda>0$, which holds  under the stated hypotheses that \eqref{spectral-prob} has a principal eigenvalue (this is a consequence of Propositions~\ref{lemma1-aux} and \ref{propconvexity}). Then it is natural to ask whether  it is always true that \eqref{ev} has a principal eigenfunction, even if \eqref{spectral-prob} does not. Thanks to Proposition~\ref{prop pri eig} one can construct examples where \eqref{ev}  has no principal eigenvalue for  $\lambda$ in some interval around 0.
\end{remark}
}


\section{Convergence of the principal eigenvalue and eigenfunction}
\label{sec conv princip eigen}
Given $\eps\ge0$ we study here the eigenvalue problem:
\begin{align}
\label{eig1}
\left\{
\begin{aligned}
& \eps \Delta w + J_\lambda *w - w  + f_u(x,0) w + \mu w = 0
\quad\text{in } \R^N
\\
& w > 0 \ \hbox{ periodic and $C^2$.}
\end{aligned}
\right.
\end{align}
We will write
\begin{align}
\label{def L eps}
L_{\eps,\lambda} w=  \eps \Delta w + J_\lambda* w - w  + f_u(x,0) w
\end{align}
and $L_{\lambda} = L_{0,\lambda}$.

In this section we will assume that $\mu_0$ is a principal eigenvalue for $-L_0$. Observe that by Corollary \ref{corollary} $\mu_\lambda$ is a principal eigenvalue of $-L_\lambda$.
By the Krein-Rutman theorem,  we know that for $\eps>0$,  $L_{\eps,\lambda}$ has a principal eigenvalue $\mu_{\eps,\lambda}$ and there are principal $C^2$ periodic eigenfunctions $\phi_{\eps,\lambda}>0$ of $L_{\eps,\lambda}$ and $\phi_{\eps,\lambda}^*>0$ of $L_{\eps,\lambda}^*$, that is,
$$
L_{\eps,\lambda} \phi_{\eps,\lambda} + \mu_{\eps,\lambda} \phi_{\eps,\lambda} =0
\quad \hbox{and}\quad
L_{\eps,\lambda}^* \phi_{\eps,\lambda}^* + \mu_{\eps,\lambda} \phi_{\eps,\lambda}^* =0  .
$$

\begin{lemma}
\change{Assume that $\mu_0$ is a principal eigenvalue for $-L_0$.}
For $\eps\ge 0$
\begin{align}
\label{principal ev sup}
\mu_{\eps,\lambda}
 &  = \sup \ \{ \mu\in \R: \exists \phi>0 \quad L_{\eps,\lambda} \phi + \mu \phi\le 0 \}
\\
\label{principal ev inf}
& = \inf \ \{ \mu\in\R: \exists \phi>0 \quad L_{\eps,\lambda} \phi + \mu \phi\ge 0 \} ,
\end{align}
where the $\sup$ and $\inf$ are taken over periodic $C^2$ periodic functions if $\eps>0$ and over continuous periodic functions if $\eps=0$.
\end{lemma}
\noindent
{\bf Proof.}
Let us write:
\begin{align*}
\mu_{\eps,\lambda}^+ & = \sup \ \{ \mu: \exists \phi>0 \quad L_{\eps,\lambda} \phi + \mu \phi\le 0 \}
\\
\mu_{\eps,\lambda}^- &= \inf \ \{ \mu: \exists \phi>0 \quad L_{\eps,\lambda} \phi + \mu \phi\ge 0 \} .
\end{align*}
Using $\phi_\eps$ in the definitions we see that
$$
\mu_{\eps,\lambda}^- \le \mu_{\eps,\lambda} \le \mu_{\eps,\lambda}^+.
$$

Let us prove $\mu_{\eps,\lambda} = \mu_{\eps,\lambda}^-$. Let $\mu \in \R$ be such that there exits $\psi>0$ $C^2$ periodic such that $L_{\eps,\lambda} \psi + \mu \psi \ge 0$. Then
$$
\mu_{\eps,\lambda} \langle \psi, \phi_{\eps,\lambda}^* \rangle
= -\langle \psi , L_{\eps,\lambda}^* \phi_{\eps,\lambda}^* \rangle
= - \langle L_{\eps,\lambda} \psi ,\phi_{\eps,\lambda}^* \rangle
\le  \mu \langle \psi, \phi_{\eps,\lambda}^* \rangle
$$
where $\langle \ , \ \rangle$ denotes $L^2$ inner product on $[0,1]^N$.
Since $\langle \psi, \phi^* \rangle>0$ we deduce that $\mu_{\eps,\lambda} \le \mu$. Hence $\mu_{\eps,\lambda} \le \mu_{\eps,\lambda}^-$.

The proof of $\mu_{\eps,\lambda}^+ \le \mu_{\eps,\lambda} $ is similar.\qed

\medskip

\begin{lemma}
\label{lemma conv mu eps}
\change{Assume that $\mu_0$ is a principal eigenvalue for $-L_0$.}
Let $\mu_{\eps,\lambda}$ be the principal eigenvalue of \eqref{eig1} in the space of \ $C^2$ periodic functions. Then
$$
\mu_{\eps,\lambda} \to\mu_\lambda \quad \hbox{ as } \eps\to 0 ,
$$
\change{and the convergence is uniform for $\lambda$ in bounded intervals.}

Let $\phi_{\eps,\lambda} $ be the principal periodic eigenfunction of $L_{\eps,\lambda}$ normalized so that
\co{cosmetics}
$$
\change{\|\phi_{\eps,\lambda}\|_{L^2([0,1]^N)} = 1.}
$$
Then
$$
\phi_{\eps,\lambda} \to \phi_\lambda
\quad
\text{ in } C(\R^N)
\text{ as } \eps\to0
$$
where  $\phi_\lambda $ is the principal periodic eigenfunction of $L_\lambda$.
\end{lemma}

\noindent
{\bf Proof.}
Under the stated hypotheses  \eqref{hyp J}, \eqref{hyp f1} on $J$ and $f$, $\phi_\lambda$ is $C^2$ by
Proposition~\ref{prop pri eig}.
Let $\mu>\mu_\lambda$.
Then
$$
L_{\eps,\lambda} \phi_\lambda + \mu \phi_\lambda = \eps \Delta \phi_\lambda + (\mu- \mu_\lambda) \phi_\lambda \ge 0
$$
if $\eps $ is small. Using formula \eqref{principal ev inf} we see that for small $\eps$, $ \mu_{\eps,\lambda} \le \mu $.
Thus
$$
\limsup_{\eps\to 0} \mu_{\eps,\lambda} \le \mu_\lambda.
$$
Using \eqref{principal ev sup} we can prove
$$
\liminf_{\eps\to 0} \mu_{\eps,\lambda} \ge \mu_\lambda.
$$

Next we prove the uniform convergence of $\phi_{\eps,\lambda}$
and for this we derive \textit{a priori} estimates. Since  $\phi_{\eps,\lambda}$ satisfies \eqref{eig1} and $f_u(x,0)$ is $C^2$ we see that $\phi_{\eps,\lambda}$ is in $C^{3,\alpha}(\R^N) $ for any $\alpha \in (0,1)$. Fix $i\in \{1,\ldots,N\}$ and differentiate \eqref{eig1}
with respect to $x_i$. Let us write
$w_i = \partial_{x_i } \phi_{\eps,\lambda}$. Then
\begin{align}
\label{eq partial eig f}
\eps \Delta w_i + g_i - \ w_i
+ f_u(x,0) w_i
+ \mu_{\eps,\lambda} w_i = 0
\quad
\text{in }\R^N,
\end{align}
where
$$
g_i(x)
=  \int_{\R^N}
\left(
\partial_{x_i} J(x-y )
- \lambda e_i
\right) e^{\lambda (y-x)\cdot e} \phi_{\eps,\lambda}(y)
\, d y + \partial^2_{x_i u }f(x,0)  \phi_{\eps,\lambda} .
$$
Let $p\ge 1$.
Multiplying \eqref{eq partial eig f} by $|w_i|^{p-2} w_i$ and integrating on the period $[0,1]^N$ we get
$$
\eps
\int_{[0,1]^N}
\Delta w_i
|w_i|^{p-2} w_i
\, d x
+
\int_{[0,1]^N}
g_i |w_i|^{p-2} w_i
\, d x
$$
$$
 +
\int_{[0,1]^N}
( -  1  + f_u(x,0)  + \mu_{\eps,\lambda} )
|w_i|^p
\, d x
= 0 .
$$
Integrating by parts
\co{formula seems correct as it is}
\begin{align*}
\eps
(p-1) \int_{[0,1]^N}
|w_i|^{p-2} |\nabla w_i|^2
&  +
\int_{[0,1]^N}
( 1  - f_u(x,0)  - \mu_{\eps,\lambda} )
|w_i|^p
\, d x
\\
& =
\int_{[0,1]^N}
g_i |w_i|^{p-2} w_i
\, d x
\end{align*}
and therefore
$$
\int_{[0,1]^N}
( 1  - f_u(x,0)  - \mu_{\eps,\lambda} )
|w_i|^p
\, d x
\le
\int_{[0,1]^N}
g_i |w_i|^{p-1}
\, d x  .
$$
By H\"older's inequality
\begin{align}
\label{ineq eigen holder}
\int_{[0,1]^N}
( 1  - f_u(x,0)  - \mu_{\eps,\lambda} )
|w_i|^p
\, d x
\le
\left(
\int_{[0,1]^N}
|w_i|^p
\right)^{1-1/p}
\left(
\int_{[0,1]^N}
|g_i|^p
\right)^{1/p} .
\end{align}
Since the operator $L_\lambda$ has a principal eigenfunction $\phi_\lambda>0$ from the relation
$$
J_\lambda* \phi_\lambda = (1 - f_u(x,0) - \mu_\lambda) \phi_\lambda
$$
we see that
$$
\inf_{x\in\R^N}
(1 - f_u(x,0) - \mu_\lambda) > 0.
$$
Since $\mu_{\eps,\lambda} \to \mu_\lambda$ as $\eps\to 0$, for sufficiently small $\eps>0$ we have
$$
(1 - f_u(x,0) - \mu_{\eps,\lambda} ) \ge c > 0 \quad
\text{for all } x\in \R^N.
$$
We deduce from this and \eqref{ineq eigen holder} that
$$
\|w_i\|_{L^p([0,1]^N)} \le C \|g_i\|_{L^p([0,1]^N)}
$$
with $C$ independent of $\eps$.
But
$$
\|g_i\|_{L^p([0,1]^N)} \le C \|\phi_{\eps,\lambda} \|_{L^p([0,1]^N)}
$$
and therefore, recalling the definition of $w_i$, we obtain
\begin{align}
\label{grad estimate}
\|\nabla \phi_{\eps,\lambda}\|_{L^p([0,1]^N)} \le  C
\| \phi_{\eps,\lambda} \|_{L^p([0,1]^N)}
\end{align}
with $C$ independent of $\eps$. Since we have normalized $
\| \phi_{\eps,\lambda} \|_{L^2([0,1]^N)}=1$, using \eqref{grad estimate} repeatedly and Sobolev's inequality we deduce that for any $p>1$
$$
\|\nabla \phi_{\eps,\lambda}\|_{L^p([0,1]^N)} \le  C
$$
for some constant $C$.  By Morrey's inequality we deduce that $\phi_{\eps,\lambda}$ is bounded in $C^\alpha([0,1]^N)$ for any $0<\alpha<1$. Therefore, for a subsequence we have that $\phi_{\eps,\lambda} \to \phi$ uniformly on $[0,1]^N$ to some continuous function $\phi$. Then, multiplying \eqref{eig1} by a periodic smooth function and integrating by parts twice we deduce that $\phi \ge 0$  is a periodic eigenfunction of $L_\lambda$ with eigenvalue $\mu_\lambda$. Then $\phi$ is a multiple of $\phi_\lambda$ and since both have $L^2$ norm equal to 1, we conclude that $\phi=\phi_\lambda$. We also deduce that the whole family $\phi_{\eps,\lambda} $ converges to $\phi_\lambda$ as $\eps\to 0$.
\qed


\section{The stationary problem}
\label{sec stationary}

In this section  we give the proof of Theorem \ref{thm 1}.
\co{section completely rewritten}
The same result for Dirichlet boundary condition appears in \cite{Cov6}.

First we state a result analogous to Theorem \ref{thm 1}
for the perturbed problem.
\begin{proposition}
\label{prop stat}
Assume \eqref{hyp f1}.
Let $\mu_{\eps}$ denote the principal periodic eigenvalue of $-L_{\eps}$ where for $\eps> 0$
$$
L_{\eps} \phi = \eps \Delta \phi + J*\phi-\phi + f_u(x,0) \phi .
$$
\co{erased "unique up to translation"}
The perturbed stationary equation \eqref{approx stat} has a
positive periodic solution if and only if $\mu_{\eps} < 0$ and
this solution is unique.
\end{proposition}
We will omit the proof, since it is very similar to
\cite{Berestycki-Hamel-Roques-I,cdm1}.

\begin{lemma}
\label{lemma bounds}
Assume $\mu_0<0$, so for $\eps>0$ small $\mu_\eps<0$ and there exists a positive solution $p_\eps$ of \eqref{approx stat}.
\co{this lemma is new (same material as before)}
Then there is a constant $C>0$ such that for $\eps>0$ small
$$
\frac {1}{C} \le p_\eps(x) \le C \quad \forall x\in \R^N.
$$
Also, $p_\eps$ is uniformly Lipschitz for $\eps>0$ small, i.e., there is $C$ such that
$$
|p_\eps(x)-p_\eps(x')|\leq C|x-x'|
\quad\text{for all }x,x'\in\R^N
$$
and for all $\eps>0$ small.
\end{lemma}

\noindent{\bf Proof.}
For the proof of upper and lower bounds, it suffices to exhibit super and subsolutions which are bounded and bounded away from zero, uniformly for $\eps>0$ small. As a super solution we just take a large fixed constant.

Let us proceed with the construction of a sub solution.
We follow an argument developed in \cite{Cov6}.
Let  $a(x):=f_u(x,0)-1$ and $\sigma:=\sup_{\R^N} a(x)$.
Since $a(x)$ is smooth and periodic  there exists a point $x_0$ such that  $\sigma=a(x_0)$. By continuity of  $a(x)$, for each $n$ there exists $\eta_n$ such that for all $x \in B_{\eta_n}(x_0)$ we have  $ |\sigma -a(x)|\le \frac{2}{n}$.

Now let us consider a sequence of real numbers $(\eps_n)_{n\in\N}$ which converges to zero such that  $\eps_n\le\frac{ \eta_n}{2}$.
 Next, let $(\chi_n)_{n\in \N}$ be  the following sequence of  cut-off functions:
\co{comment of the referee on periodic function}
\change{ $\tilde \chi_n(x):=\chi(\frac{\|x-x_0\|}{\eps_n})$} where $\chi$ is a smooth function such that $0\le \chi \le 1$, $\chi(x) = 0$ for $|x| \ge 2$ and $\chi(x) = 1$ for $|x| \le 1$.
\change{Next, we let
$$
\chi_n(x) = \sum_{k\in \Z^N} \tilde \chi_n(x-k)
$$
so that for $n$ large, $\chi_n$ is well defined, smooth, and $[0,1]^N$ periodic.}

Let us consider the following sequence of continuous periodic
functions $(a_n)_{n\in \N}$, defined by
\begin{align*}
a_n(x):= \max \{a(x),\sigma\chi_n\} .
\end{align*}
Then $\|a_n-a\|_\infty\to 0$ as $n\to\infty$.
Now consider a  $C^\infty$ regularization
$b_n(x):=\rho_{n}* a_n(x)$ where $\rho_n$ is an adequate sequence of mollifiers with support in $B_{\frac{\eps_n}{4}}(0)$, such that
$\|b_n - a_n\|_\infty\le \|a_n-a\|_\infty$.
Let $\phi_{\eps,n} > 0$ be the principal eigenfunction of the following eigenvalue problem
$$
\eps \Delta \phi_{\eps,n} +  J*\phi_{\eps,n} +b_n(x)\phi_{\eps,n} + \mu_{\eps,n} \phi_{\eps,n} = 0 \quad \hbox{in } \R^N .
$$
Since $b_n$ is constant in a small neighborhood of $x_0$, which is a point where it attains its maximum,
by Proposition~\ref{prop pri eig}, there is a  principal eigenvalue $\mu_{n}$ and eigenfunction $\phi_n>0$ for the problem
$$
J*\phi_{n} +b_n(x)\phi_{n} + \mu_{n} \phi_{n} = 0 \quad \hbox{in } \R^N.
$$
We normalize $\|\phi_n\|_{L^\infty([0,1]^N)}=1$.

Using that $\|b_n(x) -a(x)\|_{\infty} \to 0$ as $n \to \infty$, from the Lipschitz continuity with respect to $a(x)$  (Proposition~\ref{propaux1}) it follows that for $n$ big enough,  say $n\ge n_0$, we have
\co{Tried to simply proofs. Before the same argument appears twice.}
$$
\mu_n \leq \frac{\mu_{0}}{2}<0.
$$
We fix $n_0$ large so that
$$
\|b_{n_0} -a \|_\infty \le \frac{|\mu_0|}{8}
$$
Having fixed $n_0$, we work with $\eps_0>0$ small so that
$$
\mu_{\eps,n_0}  \leq \frac{\mu_0}{4}<0,
\quad\text{for all } 0<\eps\leq \eps_0,
$$
which is possible since $\mu_{\eps,n_0} \to \mu_{n_0}$ as $\eps\to 0$ by Lemma~\ref{lemma conv mu eps}.

Now for $\sigma>0$ we have
\begin{align*}
\eps \sigma \Delta \phi_{\eps,n_0}+J*\sigma\phi_{\eps,n_0} -\sigma \phi_{\eps,n_0} +f(x,\sigma\phi_{\eps,n_0})
&\ge -\big(\| a(x)- b_{n_0}(x)\|_{\infty}+\mu_{\eps,n_0}\big)\sigma\phi_{\eps,n_0} \\
 & \ \ \ +o(\sigma\phi_{\eps,n_0})\\
& \ge -\frac{\mu_{0}}{8}\sigma\phi_{\eps,n_0}+o(\sigma\phi_{\eps,n_0})>0.
\end{align*}
Therefore, for $\sigma>0$ sufficiently small, $\sigma\phi_{\eps,n_0}$ is a
subsolution of \eqref{stationary}.
By Lemma~\ref{lemma conv mu eps}, $\phi_{\eps,n_0} \to \phi_{n_0}$ uniformly in $\R^N$
as $\eps\to0$.
Since $\phi_{n_0}>0$ we find the lower
bound $p_\eps \ge 1/C$ for some $C>0$ and all $\eps>0$ small.

Let us prove now that $p_\eps$ is uniformly Lipschitz.
Let $v = \frac{\partial p_\eps}{\partial x_j}$ for some $j \in \{
1,\ldots,N\}$. Then $v$ satisfies
\begin{align*}
J * v - v  + \eps \Delta v   + f_u(x,p_\eps) v
+ f_{x_j}(x,p_\eps)= 0  \qquad x \in \R^N.
\end{align*}
We use that $f(x,u)/u$ is a decreasing
function for $u>0$. This implies that $f(x,u) - f_u(x,u) u>0$ for
all $x\in \R^N$ and all $u>0$. Since there is a fixed lower bound
for $p_\eps \ge \frac 1C$ ($\eps>0$ small) we find a fixed lower
bound for the quantity
$$
f(x,p_\eps) - f_u(x,p_\eps) p_\eps \ge \delta_0 >0 \quad \forall
x\in \R^N
$$
and all $\eps>0$ small. Then $p_\eps$ satisfies
$$
\eps\Delta p_\eps + J* p_\eps-p_\eps+f_u(x,p_\eps) p_\eps =
f_u(x,p_\eps) p_\eps - f(x,p_\eps) \le - \delta_0.
$$
By the maximum principle we conclude that
$$
|v|\le \frac{\|f_{x_j}\|_{\infty}}{\delta_0} p_\eps \le C \quad \hbox{in } \R^N .
$$
Thus $p_\eps$ is uniformly Lipschitz.
\qed

\medskip

\medskip

\noindent
{\bf Proof of Theorem~\ref{thm 1}.}
Uniqueness is proved as in \cite{cdm1,Cov6}.
Also the proof that $\mu_0<0$ is necessary for existence is very similar to \cite{cdm1,Cov6}, so we omit the details.

Assume now $\mu_0<0$ and let us prove that there exists a continuous solution.
Let $p_\eps$ be the positive solution of \eqref{approx stat}, which exists since $\mu_\eps<0$ for $\eps>0$ small.
By Lemma~\ref{lemma bounds}, $p_\eps$ is uniformly Lipschitz and therefore, up to subsequence $p_\eps$, converges uniformly in $[0,1]^N$ as $\eps\to 0$ to a continuous function $p>0$ which is periodic and solves \eqref{stationary}. By the uniqueness of the positive periodic solution of \eqref{stationary}, we have convergence of the whole family $p_\eps$.
\qed

\medskip

Directly from the previous proof we get the following result.
\begin{corollary}
Assume $\mu_{0}<0$, so $\mu_{\eps} < 0$ for $\eps>0$
small . Let $p$ be the positive continuous
periodic solution of \eqref{stationary} and $p_{\eps}$ be the
positive periodic solution of \eqref{approx stat} for $\eps>0$
small. Then
$$
p_\eps \to p \quad \hbox{uniformly as } \eps\to0 .
$$
\end{corollary}

\section{Construction of approximate pulsating fronts}
\label{s:approximate pulsating fronts}

Let $\eps>0$ be small enough so that
$$
0 = J * p - p + \eps\Delta p + f(x,p) \qquad x \in \R^N
$$
has a positive  periodic solution $p_{\eps}$, which is unique.

Here the main result is the following.
\begin{proposition}
\label{pro aprox}
Let $ \change{c_e^*(\eps)}$ be defined by \eqref{c star eps}.
For $c\ge \change{c_e^*(\eps)}$ there is a solution to
\begin{align}
\label{eq aprox}
c \partial_s \psi = M \psi - \psi + \eps \Delta \psi  + f(x,\psi) \quad \hbox{in } \R \times \R^N
\end{align}
such that
\begin{align}
\label{cond aprox}
\left\{
\begin{aligned}
& \lim_{s\to-\infty} \psi(s,x)=0
\\
&
\lim_{s\to+\infty} \psi(s,x)=p_{\eps}(x)
\\
&
\psi(s,x) \hbox{ is increasing in $s$ and periodic in $x$.}
\end{aligned}
\right.
\end{align}
\end{proposition}
\co{\changejer{I have change the notation for the regularisation parameter, keep $\eps$ for the x and introduce $\changejer{\kappa}$ for the s}}

To prove this result, we first work with an elliptic regularization  $\L_{\changejer{\kappa}}$ of  the operator $M-Id+\eps\Delta_x-c\partial_s $ as it is done  in  \cite{berestycki-hamel-cpam,Cov4,CDM2} and introduce a truncated problem as follows.
Given $\changejer{\kappa},r,R>0$, $\sigma\ge 0$ and $c\in \R$ consider the problem
\begin{align}
\label{tr 1}
\left\{
\begin{aligned}
& \L_{\changejer{\kappa}} \psi + f(x,\psi) +H(s,x) =0 \quad \hbox{in }  (-r, R) \times \R^N
\\
& \psi(s,\cdot) = \sigma \phi \quad \hbox{for } s \le - r
\\
& \psi(s,\cdot) = p_{\eps} \quad \hbox{for } s \ge R
\\
& \hbox{$\psi(s,\cdot)$ is $[0,1]^N$-periodic for all $s$}
\end{aligned}
\right.
\end{align}
where $$\L_{\changejer{\kappa}}\psi:= \int_{[-r\le s+(y-x)\cdot e \le R]} J(x-y) \psi ( s+(y-x)\cdot e ,y  ) \, d y - \psi  + \eps\Delta_x \psi + \changejer{\changejer{\kappa}} \partial_{ss}\psi - c \partial_s\psi,$$
$\phi_\eps$ is the principal periodic eigenfunction
associated \change{with} the principal eigenvalue $\mu_{\eps}$ of the following problem
$$
\eps \Delta \phi + J* \phi - \phi  + f_u(x,0) \phi  + \mu_{\eps} \phi = 0 ,
$$
and
$$
H(s,x) = \sigma \int_{[s+(y-x)\cdot e \le -r]} J(x-y)\phi_\eps(y)  \, d y
+
 \int_{[s+(y-x)\cdot e \ge R]} J(x-y) p_{\eps}(y) \, d y .
$$

\begin{proposition}
\label{pro truncated}
 There exists $\sigma_0$ such that for all $0\le \sigma \le \sigma_0$ and for
 any $c\in\R$ there exists a unique  solution of
\eqref{tr 1}.
Moreover, the corresponding  solution is  increasing in $s$,  and continuous with respect to $\sigma$ with values in $C^2( [-r,R]\times \R^N )$.
\end{proposition}

\medskip
\noindent{\bf Proof.}
Note that by construction, since $J$ is smooth then $H(s,x)$ is also smooth and the problem \eqref{tr 1} can be solved by super and sub-solutions techniques. We call a function  $\psi \in C^2(\R^N \times [-r,R])$ a super-solution of \eqref{tr 1} if
\begin{align*}
&  \L_{\changejer{\changejer{\kappa}}} \psi + f(x,\psi) +H(s,x) \le 0
, \quad -r < s < R
\\
& \psi(-r,x) \ge \sigma \phi_\eps, \quad \psi(R,x) \ge p_{\eps}(x) \quad \forall x \in \R^N
\\
& \hbox{$\psi$ is periodic in $x$} .
\end{align*}
Subsolutions are defined similarly reversing the inequalities. If there exists a subsolution $ \Psi_1 \in C^2([-r,R]  \times  \R^N)$ and  a supersolution $\Psi_2 \in C^2([-r,R]  \times  \R^N)$ such that $\Psi_1 \le \Psi_2$, then using monotone iterations one can construct a minimal solution $\underline \psi$ and a maximal solution $\overline \psi$ of \eqref{tr 1} such that $\Psi_1 \le \underline \psi \le \overline \psi \le \Psi_2$. The monotone iterations can be taken for instance of the form
$$
\psi_0 = \Psi_1
$$
and $\psi_n$ defined recursively as
\begin{align}
\label{eq iter}
\left\{
\begin{aligned}
& -  \eps \Delta_x \psi_{n+1} - \changejer{\changejer{\kappa}} \partial_{ss}\psi_{n+1}
+ c \partial_s\psi_{n+1}
+ (A+1)\psi_{n+1}
\\
&
\qquad \qquad \qquad \qquad \qquad
=
\tilde M \psi_n  + f(x,\psi_n)   + A \psi_n +H(x,s)\quad \hbox{in }  (-r, R) \times \R^N
\\
& \psi_{n+1}(-r,x) = \sigma \phi_\eps , \quad \psi_{n+1}(R,x) = p_{\eps}(x) \quad \forall x \in \R^N
\\&
 \hbox{$\psi_{n+1}$ is periodic in $x$} .
\end{aligned}
\right.
\end{align}
where $\tilde M$ denotes the operator
$$\tilde M \psi(s,x) =\int_{[-r\le s+(y-x)\cdot e \le R]} J(x-y) \psi ( s+(y-x)\cdot e ,y  ) \, d y.$$
Here $A>0$ is a large constant such that $u \mapsto f(x,u) + A u $ is increasing for all $u\in[0,\max p_{\eps}]$ and all $x$. Then the right hand side of \eqref{eq iter} is a monotone operator.

Now since, $p_{\eps}$ and $w$ are bounded and strictly positive functions,  the following quantity  $\sigma^*$ is well defined  $$\sigma^*:=\sup\{\sigma >0\, |\, \sigma \phi_\eps\le p_{\eps}\}.$$
 Take now $0\le\sigma \le \sigma^*$. Then from the definition of $H(s,x)$ we  see that
 $ p_{\eps}$ is a supersolution of \eqref{tr 1}. Indeed,  a short computation shows that
$$
 \L_{\changejer{\kappa}}[p_\eps]  + f(x,p_\eps) + H(x,s)\le  (J* p_{\eps} - p_{\eps}) + f(x,p_{\eps}) + \eps \Delta_x p_{\eps}=0.
$$

Working with $\eps>0$ sufficiently small we have  that $\mu_{\eps}<0$.
Let us now observe that when $0\le\sigma \le \sigma^*$ and $\sigma$ is small enough the function $\sigma \phi_\eps$ is a subsolution of \eqref{tr 1}.  Indeed,  as above using that $\sigma \phi_\eps\le p_{\eps}$ a short computation shows that
\begin{align*}
\L_{\changejer{\kappa}}[\sigma \phi_\eps ]  + f(x,\sigma \phi_\eps) + H(x,s)&\ge  \sigma(J* \phi_\eps - \phi_\eps) + f(x,\sigma\phi_\eps) + \eps \sigma\Delta_x \phi_\eps\\
&\ge  \sigma \phi_\eps\left( - \mu_{\eps} + \frac{f(x,\sigma \phi_\eps)}{\sigma \phi_\eps}-f_u(x,0)\right).
\end{align*}
Since $\phi_\eps$ is uniformly bounded, using the regularity of $f(x,s)$    we have  for $\sigma\ge 0$ small enough say $\sigma\le \sigma_1$
$$
\left( -\mu_{\eps} + \frac{f(x,\sigma \phi_\eps)}{\sigma \phi_\eps}-f_u(x,0)\right)\ge -\frac{\mu_{\eps}}{2}\ge 0.
$$
Thus for $\sigma\le \sigma_0:= \inf\{ \sigma_1, \sigma^*\}$, $\sigma \phi_\eps$ is a subsolution to \eqref{tr 1} with $\sigma \phi_\eps\le p_{\eps}$.

We prove now that for all  $\sigma\le \sigma_0$ the corresponding \changejer{ problem \eqref{tr 1} has a unique positive solution denoted $\psi_\sigma$.}
To this end we use  a standard   sliding method.
\co{\changejer{ I have remove the subscript $\sigma$ in the proof. I have used $\psi_\sigma$ to denote the unique solution}}
First observe that for any $0\le \sigma \le \sigma_0$,  then any bounded solution $\psi$ of the corresponding problem  \eqref{tr 1} satisfies
$$\sigma \phi_\eps<\psi <p_{\eps}.$$
Indeed, let us start with the proof of the inequality $\psi\le p_{\eps}$. Since $p_{\eps}$ is bounded away from $0$ the following quantity is well defined
$$\gamma^*:=\inf \{\gamma >0\, |\, \psi\le \gamma  p_{\eps}\}.$$
To prove the inequality, we are reduced to show that $\gamma^*\le 1$. Assume by contradiction that $\gamma^*>1$.
From the definition of $\gamma^*$, using the periodicity of the functions $\psi$, $p_{\eps}$ and a standard argument we see that there exists a point $(s_0,x_0)\in (-r,R)\times \R^N$ such that $\gamma^*p_{\eps}(s_0,x_0)=\psi(s_0,x_0)$.

Observe that  since $\frac{f(x,s)}{s}$ is a  decreasing function of $s$, the function $\gamma^*p_{\eps}$ is a supersolution of  \eqref{tr 1}.  Moreover, for some positive constant $A$ big enough, the function $\gamma^*p_{\eps}-\psi$ satisfies
\begin{align*}
&  \L_{\changejer{\kappa}}(\gamma^*p_{\eps}- \psi) -A (\gamma^*p_{\eps}- \psi)  \le 0
, \quad \text{ in }\quad (-r,R)\times \R^N
\\
& (\gamma^*p_{\eps}- \psi)(-r,x) \ge 0, \quad (\gamma^*p_{\eps}- \psi)(R,x) \ge 0 \quad \forall x \in \R^N.
\end{align*}
Since $ \L_{\changejer{\kappa}}$ is elliptic in $(-r,R)\times \R^N$ and $\gamma^*p_{\eps}(s_0,x_0)=\psi(s_0,x_0)$,   from  the strong maximum principle it follows that
$$\gamma^*p_{\eps}\equiv \psi \qquad \text{ in }  \qquad (-r,R)\times \R^N,$$
which is impossible since $ \gamma^*p_{\eps}(x)>p_{\eps}(x)\ge \sigma \phi_\eps(x)=\psi(-r,x)$.
Therefore we have $\gamma^*\le 1$ and $\psi \le p_{\eps}$.
The strict inequality comes from the strong maximum principle.
 Now observe that to obtain  the other inequality $\sigma \phi_\eps<\psi$ we can just reproduce the above argumentation with $\sigma \phi_\eps$ in the role of $\psi$ and $\psi$ in the role of $p_{\eps}$.

We are now in position to prove the uniqueness of the solution of  \eqref{tr 1}.
Suppose $\psi_1$, $\psi_2$ are 2 solutions of \eqref{tr 1}. Define the following  continuous functions
$$
\bar \psi_1(s,x):=\begin{cases} \sigma \phi_\eps (x)  & \hbox{if } s < -r \; \text{ and } x\in \R^N\\
 \psi_1(s,x)  & \hbox{if } -r\le s \le R \; \text{ and } x\in \R^N\\
 p_{\eps}(x)  & \hbox{if } s > R \; \text{ and } x\in \R^N\\
\end{cases}
$$
and
$$
\bar \psi_2(s,x):=\begin{cases} \sigma \phi_\eps (x)  & \hbox{if } s < -r \; \text{ and } x\in \R^N\\
 \psi_2(s,x)  & \hbox{if } -r\le s \le R \; \text{ and } x\in \R^N\\
 p_{\eps}(x)  & \hbox{if } s > R \; \text{ and } x\in \R^N.\\
\end{cases}
$$
 Note that with this notation the equation  \eqref{tr 1} satisfied by $\psi_1$ and $\psi_2$ can be  rewritten
\begin{equation}
  \eps \Delta \psi_{i} +\changejer{\kappa} \partial_{ss}\psi_i - c\partial_s \psi_i  -\psi_i +f(x,\psi_i)=-M\bar \psi_i \quad \text{ in }\quad (-r,R)\times \R^N,
\end{equation}
with $i\in \{1,2\}$.

Let us define
$$
\bar \psi_1^\tau(s,x) :=  \bar \psi_1(s+\tau,x)
$$
with $\tau \in \R$. Obviously, we have
$$
 \bar \psi_1^\tau(s,x) :=   \psi_1(s+\tau,x)  \qquad \text{ in }\qquad (-r,R-\tau)\times \R^N
$$

We claim that for all $\tau \in [0, R+r]$
\begin{align}
\label{ineq phi lambda}
\bar \psi_1^\tau(s,x) > \bar \psi_2(s,x) \quad \hbox{for }  (s, x) \in \R\times\R^N.
\end{align}
By construction we easily see that $\bar\psi_1 ^{R+r}\ge \bar\psi_2$ in $\R\times\R^N$  since we know that   $$\sigma \phi_\eps\le \psi_i \le p_{\eps} \quad \hbox{for }  (s, x) \in \R\times\R^N.$$

Moreover, using that  we have a strict inequality in $(-r,R)$, that is to say   $$\sigma \phi_\eps<\psi_i < p_{\eps} \quad \hbox{for }  (s, x) \in (-r,R)\times\R^N,$$
 we can find a positive $\eps$ such that  for any $\tau \in [R+r-\eps,R+r]$ we have
 $$
\bar \psi_1^\tau(s,x) > \bar \psi_2(s,x) \quad \hbox{for }  (s, x) \in \R\times\R^N.
$$

Note also that  by construction for all $\tau\ge 0$ we  have
\begin{equation} \label{ptw.eq.slid.esti1}
\bar\psi_1^{\tau}\ge \bar\psi_2 \qquad \text{ in }\qquad \left((-\infty,-r]\cup [R-\tau,+\infty)\right)\times \R^N.
\end{equation}

Now let us define
$$
\tau^* = \inf\{ \tau \in [0,R] :
\bar\psi_1^{\tau'} \ge \bar \psi_2 \quad \hbox{for }\tau' \in [\tau,R+r]
\}
$$
then $0\le \tau^* <R+r$. Assume that $\tau^*>0$.
in this case
$$
\bar \psi_1^{\tau^*} \ge \bar \psi_2 \quad\hbox{in  } \R\times \R^N
$$
and since $J\ge 0$ we have
$$M(\bar \psi_1^{\tau^*}-\bar \psi_2)\ge 0.$$

Now, fix $A>0$ large so that $f(x,u) + A u$ is monotone increasing in $[0,\max p_{\eps}]$.
Let us denote $z:=\bar\psi_1^{\tau^*}-\bar\psi_2$. Then using the definition of $\bar \psi_1^\tau$ and $\bar \psi_2$ in $(-r,R-\tau^*)\times \R^N$, we have
\begin{align*}
&\eps \Delta z +\changejer{\changejer{\kappa}} \partial_{ss}z - c\partial_s z  -(A+1)z \le -M(\bar\psi_1^{\tau^*}-\bar\psi_2) \le 0,\\
&z(-r,x)>0 \qquad \text{ for all } \; x\in \R^N,\\
&z(R-\tau^*,x) >0 \qquad \text{ for all } \; x\in \R^N.
\end{align*}
By the strong maximum principle, it follows that $z>0$ in $(-r,R-\tau^*)\times \R^N$.
Therefore, we have $\bar\psi_1^{\tau^*}-\bar\psi_2>0$ in $[-r,R-\tau^*]\times \R^N$ and by continuity for $\delta$ small
we have for any $\tau$ in $(\tau^*-\delta,\tau^*)$
\begin{equation}
\bar\psi_1^{\tau}-\bar\psi_2 \ge 0 \qquad \text{ in }\qquad [-r,R-\tau]\times \R^N.
\end{equation}
Combining the later with \eqref{ptw.eq.slid.esti1} it follows that
for any positive $\tau$ in $(\tau^*-\delta,\tau^*)$ we have
\begin{equation*}
\bar\psi_1^{\tau}-\bar\psi_2 \ge 0 \qquad \text{ in }\qquad \R\times \R^N,
\end{equation*}
which contradicts the definition of $\tau^*$. Therefore, $\tau^* =0$ and $\bar \psi_1 \ge \bar \psi_2$.
By interchanging the role of  $\psi_1$ and $\psi_2$ in the above argument we end up with $\bar \psi_1\ge \bar \psi_2 \ge \bar \psi_1$, which prove the uniqueness of the solution of  \eqref{tr 1}.

Taking $\psi_2 = \psi$ in \eqref{ineq phi lambda} shows that $\psi$ is increasing in $s$.
Finally, \changejer{denoting $\psi_\sigma$ the unique solution of the corresponding problem \eqref{tr 1} one can see that  the map} $\sigma \mapsto \psi_\sigma$ is continuous, thanks to the uniqueness of the solution to \eqref{tr 1} and standard elliptic estimates.
\qed

\begin{proposition}\label{prop norm}
Suppose $c > \change{c_e^*(\eps)}$. Then there exists $r_0>0$, $\changejer{\changejer{\kappa}}(c)>0$ and $k>0$ such that for $r\ge r_0$, $R\ge r_0$, $\changejer{\changejer{\kappa}}\le \changejer{\changejer{\kappa}}(c)$ there is $\sigma \in (0, \sigma_0)$ for which the unique increasing solution $\psi$ of \eqref{tr 1} satisfies
$$
\max_{x\in [0,1]^N} \psi(0,x) = \frac{1}{k} \min_{\R^N} p_{\eps}.
$$
\end{proposition}
\noindent{\bf Proof.}
Let $\psi_\sigma$ denote the unique solution of \eqref{tr 1} constructed in Proposition~\ref{pro truncated}.

Choose $k>0$, so that $$\sigma_0\max_{\R^N}\phi_\eps>\frac{1}{k} \min_{\R^N} p_{\eps},$$
where $\phi_\eps$ denote the positive periodic principal eigenfunction associated \change{with} the eigenvalue problem
$$ J*  \phi - \phi + \eps \Delta \phi + f_u(x,0) \phi + \mu_{\eps} \phi=0 .$$

Observe that since $\psi_\sigma$ is increasing in $s$, we have $\max_{\R^N}\psi_{\sigma_0}(0,x)>\frac{1}{k} \min_{\R^N} p_{\eps}$.
Next we prove that for  $\sigma= 0$, we have $\max_{x\in \R^N} \psi_0(0,x) < \frac{1}{k} \min_{\R^N} p_{\eps}$.

Recall that $$\change{c_e^*(\eps)}:=\inf_{\lambda>0}\left(-\frac{\mu_{\eps,\lambda}}{\lambda}\right),$$
where \changejer{$\mu_{\eps,\lambda}$} is the principal periodic eigenvalue of the problem
$$
J_\lambda * \phi - \phi + \eps \Delta \phi + f_u(x,0) \phi + \mu_{\eps,\lambda} \phi =0.
$$

Since $c> \change{c_e^*(\eps)}$ there is $\bar\lambda >0$ such that $c \bar\lambda  +\mu_{\eps,\bar\lambda}>0$.
Let us denote $\changejer{\phi_{\eps,\bar\lambda}}$ the principal periodic eigenfunction associated with $\changejer{\mu_{\eps,\bar\lambda}}$ and consider the function
$$w:=e^{\bar\lambda(s-s_0)}\changejer{\phi_{\eps,\bar\lambda}},$$
where \changejer{$s_0\in \R$ is chosen so that $$e^{-\bar\lambda s_0}\max_{\R^N}\phi_{\eps,\bar \lambda}<\frac{1}{k} \min_{\R^N} p_{\eps},$$}
and take $R>0$ large so that
$$
e^{ \bar\lambda (R-s_0)} \changejer{\min_{\R^N}\phi_{\eps,\bar\lambda}} \ge p_{\eps}(x).
$$
Since $w$ is monotone increasing  in $s$ we have
\co{\changejer{I have remove the extra $\phi$ and remove the subscript  on  $\phi$} }
$$
w(s,x) \ge p_{\eps}(x)\quad \text{ for  any } \; (s,x) \in [R,+\infty)\times \R^N.
$$
Finally,  observe that
$$
e^{\bar\lambda(-r-s_0) } \changejer{\phi_{\eps,\bar \lambda}}(x) \ge 0 \quad \text{ for  any } \; (s,x) \in \R\times \R^N.
$$

We claim that the function $w$ is a supersolution of \eqref{tr 1} with $\sigma=0$ for $\changejer{\changejer{\kappa}}$ small enough. Indeed,  in $(-r,R)$ we have
\changejer{
\begin{align*}
\L_{\eps}w +f(x,w)+H(s,x)&\le
(J_{\bar\lambda} * \phi_{\eps,\bar\lambda} - \phi_{\eps,\bar\lambda} + \eps \Delta \phi_{\eps,\bar\lambda} +f_u(x,0)\phi_{\eps,\bar\lambda}- c\bar\lambda \phi_{\eps,\bar\lambda} +\changejer{\kappa}\bar\lambda^2 \phi_{\eps,\bar\lambda})e^{\bar\lambda s}\\
&\le -( \mu_{\eps,\bar\lambda} +c\bar \lambda -\changejer{\changejer{\kappa}}\bar\lambda^2)w
\end{align*}
}
Therefore, for $\changejer{\kappa}\le \frac{c+ \mu_{\bar\lambda} }{\bar\lambda^2}=:\changejer{\kappa}(c)$ we have

\begin{align*}
&\L_{\eps}w +f(x,w)+H(s,x)\le 0\qquad \text{ for all } \; (s,x)\in (-r,R)\times\R^N,\\
&w(-r,x)>0 \qquad \text{ for all } \; x\in \R^N,\\
&w(R,x) >p_{\eps} \qquad \text{ for all } \; x\in \R^N.
\end{align*}
Since $0$ is a subsolution of \eqref{tr 1} with $\sigma=0$ and $w\ge 0$
using the  uniqueness of the solution of \eqref{tr 1} we must have
$\psi_0(s,x)\le \changejer{w(s,x)}$. Therefore $$\max_{\R^N}\psi_0(0,x)\le \max_{\R^N}w(0,x)<\min_{\R^N}\frac{p_{\eps}}{k}.$$

With $R>0$ fixed, we see that the map $\sigma \in [0,\sigma_0] \mapsto \psi_\sigma$ is continuous, and at $\sigma_0$ satisfies $ \max\psi_{\sigma_0}(0,x) > \min \frac{p_{\eps}}{k} $ and $\max\psi_{0}(0,x) < \min \frac{p_{\eps}}{k} $. By  continuity there is $\sigma \in [0,\sigma_0]$ such that $\max\psi_\sigma(0,x) = \min \frac{p_{\eps}}{k}$.
\qed

\begin{proposition}
For $c> \change{c_e^*(\eps)}$ and $\changejer{\kappa}\le \changejer{\kappa}(c)$ there is a solution to
\begin{align}
\label{eq aprox-reg}
c \partial_s \psi = M \psi - \psi + \eps \Delta \psi +  \changejer{\kappa}\partial_{ss}\psi  + f(x,\psi) \quad \hbox{in } \R \times \R^N
\end{align}
such that
\begin{align*}
& \lim_{s\to-\infty} \psi(s,x)=0
\\
&
\lim_{s\to+\infty} \psi(s,x)=p_{\eps}(x)
\\
&
\psi(s,x) \hbox{ is increasing in $s$ and periodic in $x$.}
\end{align*}
\end{proposition}

\noindent{\bf Proof.}
For $r>0$ large, let $\psi_r$ be the solution of \eqref{tr 1} with $R=r$  obtained in Proposition \ref{prop norm} where $\sigma = \sigma(r) \in (0,\sigma_0)$ is such that
\begin{align}
\label{normalization}
\max_{x\in  \R^N} \psi_r(0,x) =  \min_{x\in \R^N} \frac{p_{\eps}(x)}{k} .
\end{align}
We let $r\to\infty$. Since $\psi_r$ is locally bounded in $C^{1,\alpha}$, there is a subsequence such that $\psi_r$ converges locally in $C^{1,\alpha}$ to a function $\psi :\R \times \R^N$ which satisfies \eqref{eq aprox-reg} with the speed $c$, is increasing in $s$ and periodic in $x$.

The limit $w(x) =\lim_{s\to-\infty} \psi(s,x)$ exists and is a solution of the stationary problem. By Proposition \ref{prop stat} part b)
this solution is either 0 or the unique positive stationary solution $p_{\eps}$. By
\eqref{normalization} we conclude that $w\equiv 0$. Similarly $\lim_{s\to+\infty} \psi(s,x) = p_{\eps}(x)$.

\qed

In the next proposition we establish some \textit{a priori} estimates satisfied by the solutions of \eqref{eq aprox-reg}. Namely, we have

\begin{proposition}\label{pro aprox esti}
Let $c> \change{c_e^*(\eps)}$ and $\changejer{\changejer{\kappa}}\le \changejer{\changejer{\kappa}}(c)$ then the solution $(\psi_{\kappa,\eps},c)$ of \eqref{eq aprox-reg} satisfies
\begin{itemize}
\item[(i)]
\begin{equation*}
c\int_{\R\times \q}|\partial_s \psi_{\kappa,\eps}|^2 =-\frac{\eps}{2}\int_{\q}|\nabla_x p_{\eps}|^2 - \changejer{\frac{1}{4}}\int_{\q^2}\tilde \j(x,y)(p_{\eps}(x)-p_{\eps}(y))^2+\int_{\q}F(x,p_{\eps})
\end{equation*}
where $\q=[0,1]^N$ and $\tilde\j=\sum_{k\in \Z^N} J(x-y-k)$ is a symmetric positive kernel.
\item[(ii)] For all compact set $\K\subset \R\times\R^N$, there exists $R>0$, a constant $\gamma(R)$ and $n\in\N$ so that
\begin{equation*}
\int_{\K}|\nabla_x \psi_{\kappa,\eps}|^2 \le \gamma(R)(2n)^N.
\end{equation*}
\item[(iii)] Given $R>0$, let
$$
Q_R = \{ (s,x) \in \R\times\R^N : |x|<R, \, |s|<R \} .
$$
Then there exists  positive constant $M,M'$ independent of $\eps$ such that
 \begin{align*}
&\sup_{Q_{R/4}}|\nabla_x \psi_{\kappa,\eps}| \le M
\left(|c|+\frac{1}{R} + R \left(2\sup_{Q_R}|p_{\eps}(x)|+\sup_{Q_R}|f_u(x,0)|\right)
\right)\sup_{Q_R}|\psi_{\kappa,\eps}|\\
 &\sup_{Q_{R/4}}
\frac{|\psi_{\kappa,\eps}(t_1,x)-\psi_{\kappa,\eps}(t_2,x)|}{|t_1-t_2|^{\frac{1}{2N}}} \le M' \sup_{Q_{R}}|\nabla_x \psi_{\kappa,\eps}|.
\end{align*}
\end{itemize}
\end{proposition}
We give the proof of this Proposition in Appendidx~\ref{appendix}.
%
We are now in a position to prove the Proposition ~\ref{pro aprox}

\medskip
\noindent{\bf Proof of Proposition~\ref{pro aprox}.}
Let us first assume that $c>\change{c_e^*(\eps)}$. Then from the above construction, for any $\changejer{\kappa}\le \changejer{\kappa}(c)$, there exists a function $\changejer{\psi_{\kappa,\eps}(s,x)}$ increasing in $s$ and periodic in $x\in \R^N$ that is solution of \eqref{eq  aprox-reg}. Without loss of generality, we can assume that $\psi_{\kappa,\eps}$ is normalized as follows
\begin{equation*}
\max_{\R^N}\psi_{\kappa,\eps}(0,x)=\min_{\R^N}\frac{p_{\eps}}{k}.
\end{equation*}
We let $\changejer{\kappa}\to 0$ along a sequence. Thanks to the apriori estimates of Proposition \ref{pro aprox esti}, we can extract a subsequence  of $(\psi_{\kappa_n,\eps})_{n\in\N}$ which converges locally uniformly in $\R\times\R^N$ to a function $\psi_\eps\in H_{loc}^1(\R^N) \cap C^\alpha(\R \times \R^N)$ for some $\alpha\in(0,1)$, that satisfies \eqref{eq aprox} in the sense of distributions.
Since $\psi_{\kappa_n,\eps}$ is periodic in $x$, monotone increasing in  $s$, and $0\le \psi_{\kappa_n,\eps}\le p_\eps$, we also have that $\psi_\eps$ is periodic in $x$,  monotone non decreasing in $s$, and
$0\le \psi_\eps\le p_{\eps}$. Note also that from the normalization condition, since $\psi_{\kappa_n,\eps}\to \psi_\eps$ locally uniformly, we also deduce that
\begin{equation}
\label{eq-lim-normali}
\max_{\R^N}\psi_{\eps}(0,x)=\min_{\R^N}\frac{p_{\eps}}{k}.
\end{equation}
Furthermore, using standard parabolic estimate, on can show that $\psi_\eps$ is a classical solution of \eqref{eq aprox}. Thus $\psi_\eps$ satisfies
\begin{align*}
\left\{
\begin{aligned}
& \eps \Delta \psi_\eps -c\partial_s\psi_\eps +M[\psi_\eps] -\psi_\eps + f(x,\psi_\eps) =0 \quad \hbox{in }  \R \times \R^N,
\\
& 0\le \psi \le p_{\eps},  \quad \partial_s\psi\ge 0 \quad \hbox{in }  \R \times \R^N
\\
& \hbox{$\psi_\eps(s,\cdot)$ is $[0,1]^N$-periodic for all $s$.}
\end{aligned}
\right.
\end{align*}
By standard estimates the limit $w(x) =\lim_{s\to-\infty} \psi_\eps(s,x)$ exists and is a solution of the stationary problem. By Proposition \ref{prop stat} part b)
this solution is either 0 or the unique positive stationary solution $p_{\eps}$. By
\eqref{eq-lim-normali} we conclude that $w\equiv 0$. Similarly $\lim_{s\to+\infty} \psi_\eps(s,x) = p_{\eps}(x)$.
\qed

\section{Estimates for $L_{\eps,\lambda}$}
\label{sec estimates L eps}
Recall the notation from \eqref{def L eps}:
$$
L_{\eps,\lambda} u = \eps \Delta u + J_\lambda * u - u + f_u(x,0)u .
$$
\co{\changejer{I have switch the parameter in the notation of  the operator, }}
\co{\changejer{I have recall the definition  of $\mu_\eps?$}}
\begin{lemma}
\label{lemma L infty epsilon}
Let $\lambda$ be such that  $0<\lambda c < - \mu_{\eps,\lambda}$, \changejer{ where $\mu_{\eps,\lambda}$ is the principal periodic eigenvalue of the operator $-L_{\eps,\lambda}$ defined in section \ref{sec conv princip eigen}}.
If $u \in C^2(\R^N)$, $u \ge 0 $ is a periodic solution to
$$
L_{\eps,\lambda} u - \lambda c u = h \quad \hbox{in }\R^N
$$
then
$$
\|u\|_{L^\infty([0,1]^N)} \le \change{C_{\eps,\lambda}} \|h\|_{L^\infty([0,1]^N)}.
$$
\end{lemma}
Note that
for any $\eps>0$ and $0<\lambda_0<\lambda_1<-\mu_{\eps,\lambda}/c$ we have
$$
\sup_{\lambda_0 \le \lambda \le \lambda_1} C_{\eps,\lambda} < \infty,
$$
but the constant depends on $\eps$.

\noindent
{\bf Proof.}
Let $\phi_{\eps,\lambda}^*$ be the principal eigenfunction of the adjoint operator $L_{\eps,\lambda}^*$. Then multiplying the equation by $\phi_{\eps,\lambda}^*$ and integrating we find
$$
(-\mu_{\eps,\lambda} - \lambda c) \int_{[0,1]^N} u \phi_{\eps,\lambda}^* =  \int_{[0,1]^N} h \phi_{\eps,\lambda}^* .
$$
Since $\lambda c < - \mu_{\eps,\lambda}$, $u\ge 0$ and $\phi_{\eps,\lambda}^*$ is strictly positive and bounded, we obtain
$$
\|u\|_{L^1([0,1]^N)} \le \change{C_{\eps,\lambda}} \|h\|_{L^1([0,1]^N)}.
$$
The uniform norm follows because of standard elliptic estimates for the operator $L_{\eps,\lambda}$.
\qed

\begin{proposition}
\label{prop uniform bound}
There is $\rho>0$, such that for any $0<\rho'<\rho$ there is $\eps_0>0$ and $\change{\overline C} $ such that for any $0<\eps\le \eps_0$, any $\lambda$
that satisfies  $(-\change{\mu_{\eps,\lambda}}-\rho)/c \le \lambda \le (-\change{\mu_{\eps,\lambda}}-\rho')/c$ and
any $u \ge 0 $ that is a periodic solution to
\begin{align}
\label{eq u}
L_{\eps,\lambda} u  - \lambda c u = h \quad \hbox{in }\R^N
\end{align}
for some $h \in L^\infty$ we  have
$$
\|u\|_{L^\infty([0,1]^N)} \le \change{\overline C}  \|h\|_{L^\infty([0,1]^N)}.
$$
The constant $\rho>0$ does not depend on $\eps$ or $\lambda$.
\end{proposition}
\noindent
{\bf Proof.}
\change{Let $\mu_\lambda$ be the principal eigenvalue of $-L_\lambda$.}
Recall that $ \inf_{x \in [0,1]^N } (1 - f_u(x,0) - \change{\mu_0}) > 0$, so we can fix $\rho>0$ such that
$ \inf_x ( 1 - f_u(x,0) - \change{\mu_0} - \rho > 0 ) $.
\change{Since $\mu_\lambda\leq \mu_0$, see Proposition~\ref{propconvexity}, also
$ \inf_x ( 1 - f_u(x,0) - \mu_\lambda- \rho > 0 ) $.}
Let $0<\rho'<\rho$ and let us proceed by contradiction. Assume that there exist sequences $\eps_n \to 0$,
\change{
$\lambda_n \in\R$,
periodic functions $(h_n)$ in $L^\infty$, $(u_n)$ in $C^2$,
such that:
$\lambda_n$ satisfies
$(-\mu_{n}-\rho)/c \le \lambda_n \le (-\mu_{n}-\rho')/c$,
\change{where $\mu_n = \mu_{\eps_n,\lambda_n}$},
$u_n$ solves \eqref{eq u} and}
$$
\| h_n \|_{L^\infty} \to 0 \quad \hbox{ and } \quad
\| u_n \|_{L^\infty} =1.
$$
We write equation \eqref{eq u} as
\begin{align}
\label{eq un}
\eps_n \Delta u_n - a_n(x) u_n = -g_n
\end{align}
where
$$
a_n(x) = 1 - f_u(x,0) + \lambda_n c
\quad
\hbox{and}
\quad
g_n = J_{\lambda_n} u_n - h_n.
$$
After extracting a subsequence we may assume that $\lambda_n\to \lambda$,  $u_n \to u$ weakly-* in $L^\infty([0,1]^N)$ and then $J_{\lambda_n}  u_n \to J_{\lambda}  u $ uniformly. Hence $g_n \to g = J_\lambda u$ uniformly, and $g$ is continuous. By Lemma~\ref{lemma conv mu eps} we have \change{$\mu_{n} = \mu_{\eps_n,\lambda_n} \to \mu_\lambda$} as $n\to\infty$.
Since
$$
a_n(x) = 1 - f_u(x,0) + \lambda_n c \ge 1 - f_u(x,0) - \change{\mu_{n}} - \rho
$$
and $ 1 - f_u(x,0) - \change{\mu_\lambda} - \rho > 0$, by working with $n$ large we may assume that
$$
\inf_{x} a_n(x) \ge a_0 > 0 \quad \hbox{for all } n.
$$
Note that $a_n \to a = 1 - f_u(x,0) + \lambda c$, which is a continuous positive function, and the convergence is uniform.
We claim that $u_n \to g/a$ uniformly. For the next argument we will assume that $g_n > 0$, which we can achieve by replacing $u_n$ by $u_n + M$ and $g_n$ by $g_n+ a_n M$ where $M>0$ is large. Note that \eqref{eq un} and  $g_n \to g$ uniformly still hold. Let $0<\sigma<1/2$ and $x_0 \in \R^N$.
By uniform convergence $g_n\to g$, $a_n \to a$ and the continuity of $g$ and $a$, we have
$$
\inf_{x\in B_r(x_0)} \frac{g_n(x)}{\beta + a_n(x)} \ge (1-\sigma)  \frac{g(x_0)}{a(x_0)}
\quad \hbox{in } B_r(x_0)
$$
provided we choose $r>0$, $\beta>0$ small and $n\ge n_0$ with $n_0$ large, and this is uniform in $x_0$.
Let $z$ be the principal eigenfunction for $-\Delta$ in $B_r(x_0)$ such that $\max_{B_r(x_0)} z = 1$ and let $\nu_r = C /r^2$ be the corresponding principal eigenvalue, that is,
\begin{align*}
\left\{
\begin{aligned}
& \Delta z + \nu_r z = 0 , \quad z>0 \qquad \hbox{in } B_r(x_0)
\\
& z=0 \qquad  \hbox{on } \partial B_r(x_0) .
\end{aligned}
\right.
\end{align*}
Define
$$
v_n = u_n - z d_n
\quad
\hbox{where}
\quad
d_n = \inf_{B_r(x_0)} \frac{g_n(x)}{\nu_r \eps_n + a_n(x)}.
$$
Then
$$
\eps_n \Delta v_n - a_n v_n = - g_n + d_n ( \eps_n \nu_r + a_n) z \le 0
$$
by the choice of $d_n$ and $z\le 1$. Since $v_n = u_n \ge 0$ on $\partial B_r(x_0)$ by the maximum principle we deduce that
$$
u_n \ge \left( \inf_{B_r(x_0)} \frac{g_n(x)}{\nu_r \eps_n + a_n(x)}\right) z \quad \hbox{in } B_r(x_0).
$$
In particular, if $n\ge n_0$ is large enough so that $\nu_r \eps_n \le \beta$ we obtain
$$
u_n(x_0) \ge (1-\sigma)  \frac{g(x_0)}{a(x_0)} .
$$
This proves that
$$
\liminf_{n\to\infty}  \inf_{x} (u_n - g/a) \ge 0.
$$
A similar argument shows that
$$
\limsup_{n\to\infty}  \sup_{x} (u_n - g/a) \le 0
$$
which proves the uniform convergence $u_n \to g/a$. We deduce that $u = g/a$, and therefore $u$ solves the equation
$$
J_\lambda u - u + f_u(x,0) u - \lambda c u = 0 .
$$
But since $\|u_n\|_{L^\infty}=1$ and $u_n$ converges uniformly we also deduce that $\|u\|_{L^\infty}=1$. Moreover $u\ge 0$. Then necessarily $\lambda c$ is the principal eigenvalue
\change{$-\mu_\lambda$ of $L_\lambda$}.
\change{This not possible because we assumed $\lambda_n c \le -\change{\mu_n} - \rho'$, so $\lambda c \leq -\mu_\lambda-\rho'$, a contradiction.}
\qed

\section{Exponential bounds}
\label{sec Exponential bounds}
Suppose we have a solution of
\begin{align}
\label{eq a}
\left\{
\begin{aligned}
& c \psi_s =  \eps \Delta\psi + M[\psi] - \psi + f(x,\psi) \quad \forall s\in \R, x\in \R^N
\\
& \hbox{$\psi(\cdot,x)$ is nondecreasing for all $x$}
\\
&  \hbox{$\psi(s,\cdot)$ is $[0,1]^N$ periodic for all $s$}
\\
& \psi(s,x) \to 0 \quad \hbox{ as } s\to-\infty
\\
&
\psi(s,x) \to p_\eps(x) \quad \hbox{ as } s\to\infty.
\end{aligned}
\right.
\end{align}
Let $\delta>0$ be fixed. We assume the following normalization on $\psi$:
\begin{align}
\label{normalization phi}
\max_{x\in [0,1]^N} \psi(0,x) = \delta.
\end{align}

Let $\lambda_\eps(c)$ be the smallest positive $\lambda$ such that $c = -\frac{\mu_{\eps,\lambda}}{\lambda} $.
The main result in this section is the following.
\begin{proposition}
\label{prop exp decay}
For any $0<\lambda<\lambda_\eps(c)$ there are $\delta>0$,  $C>0$  such that if $\psi$ satisfies \eqref{eq a} and \eqref{normalization phi}, then
\begin{align}
\label{decay phi}
\psi(s,x)  \le C e^{\lambda s}  \quad  \forall x \in \R^N , \
\quad \change{\forall s\leq 0,}
\end{align}
where $C$ does not depend on $\eps>0$.
\end{proposition}

As a corollary we have:
\begin{proposition}
\change{For all $\eps>0$ small and any fixed $\lambda$ such that
$0<\lambda < \lambda_\eps(c)$ there exists $C_\lambda$ independent of $\eps$ such that  if $\psi$ satisfies \eqref{eq a} and \eqref{normalization phi}, then}
\begin{align}
\label{decay phis}
|\psi_s(s,x)| \le C_\lambda e^{\lambda s} \quad
\forall s\le 0, \ \forall x \in \R^N
\\
\label{decay phix}
\eps^{1/2} |\nabla_x \psi(s,x)| \le C_\lambda e^{\lambda s} \quad
\forall s\le 0, \ \forall x \in \R^N
\\
\label{decay phixx}
\eps |\nabla^2_x\psi(s,x)| \le C_\lambda e^{\lambda s} \quad
\forall s\le 0, \ \forall x \in \R^N .
\end{align}
\end{proposition}
The proof of this proposition is based on scaling in the $x$ variable and applying Schauder estimates for parabolic equations. We omit the proof.

The proof has several steps.

\begin{lemma}
\label{lemma a}
There exists $\lambda_0>0$ and $C>0$ such that if $\delta>0$ is sufficiently small and $\psi$ satisfies \eqref{eq a} and \eqref{normalization phi}, then
\begin{align}
\label{eq lemma a}
\int_{[0,1]^N} \int_{-\infty}^\infty \psi(s,x) e^{-\lambda s} \, d s \, d x  \le C \quad  \forall 0<\lambda \le \lambda_0
\end{align}
where the constants do not depend on $\eps>0$. Moreover,
$$
 \int_{-\infty}^\infty \psi(s,x) e^{-\lambda s} \, d s   \le C_\eps \quad  \forall 0<\lambda \le \lambda_0
$$
where $C_\eps$ depends on $\epsilon$.
\end{lemma}
\noindent
{\bf Proof.}
Let $\eta_n:\R\to\R$ be a s smooth function such that $\eta_n(s) = 1$ for all $s\ge-n$, $\eta_n(s) =0$ for all $s\le - 2 n$, $\eta_n'\ge0$. Let $\lambda>0$ and define
$$
U_n(x,\lambda)  = \int_{-\infty}^\infty \psi(s,x) e^{-\lambda s} \eta_n(s) \, d s.
$$
We multiply \eqref{eq a} by $\eta_n(s) e^{-\lambda s}$ and integrate on $(-\infty,\infty)$. The term involving $M\psi$ yields
\begin{align*}
\int_{-\infty}^\infty M \psi(s,x) \eta_n(s) e^{-\lambda s} \, d s & =
\int_{-\infty}^\infty \int_{\R^N} J(x-y) \psi(s + (y-x)\cdot e,y) \eta_n(s) e^{-\lambda s} \, d y \, d s
\end{align*}
\begin{align*}
=\int_{\R^N} J(x-y) e^{- \lambda(x-y)\cdot e} \int_{-\infty}^\infty \psi(s + (y-x)\cdot e,y) \eta_n(s) e^{-\lambda( s + (y-x)\cdot e )} \, d s \, d y
\end{align*}
\begin{align*}
= \int_{\R^N} J(x-y) e^{- \lambda(x-y)\cdot e} \int_{-\infty}^\infty \psi(\tau,y) e^{-\lambda \tau} \eta_n(\tau-(y-x)\cdot e) \, d \tau \, d y
\end{align*}
and we write this term as
$$
J_\lambda U_n(\cdot,\lambda) + \int_{\R^N} J(x-y) e^{- \lambda(x-y)\cdot e} \int_{-\infty}^\infty \psi(\tau,y) e^{-\lambda \tau} \left[ \eta_n(\tau-(y-x)\cdot e) - \eta_n(\tau) \right] \, d \tau \, d y
$$
Hence
\begin{align}
\label{eq 2}
\eps \Delta U_n + J_\lambda U_n - U_n + f_u(x,0) U_n - c \lambda U_n = D_n+E_n+F_n
\end{align}
where
$$
D_n = \int_{\R^N} J(x-y) e^{- \lambda(x-y)\cdot e} \int_{-\infty}^\infty \psi(\tau,y) e^{-\lambda \tau} \left[  \eta_n(\tau) -\eta_n(\tau-(y-x)\cdot e)  \right] \, d \tau \, d y
$$
$$
E_n = \int_{-\infty}^\infty ( f(x,\psi(s,x)) - f_u(x,0) \psi(s,x) ) e^{-\lambda s}  \eta_n(s) \, d s
$$
$$
F_n = - c \int_{-\infty}^\infty \psi(s,x) \eta_n'(s) e^{-\lambda s} \, d s.
$$
Observe that in $D_n$, we can assume that the integral in $y$ ranges on $|y-x|\le 1$ (because we assume that $J$ has support contained in the unit ball). Then $|(y-x)\cdot e| \le 1$ and since $\eta$ is nondecreasing
$$
\int_{\R^N} J(x-y) e^{- \lambda(x-y)\cdot e} \int_{-\infty}^\infty \psi(\tau,y) e^{-\lambda \tau} \eta_n(\tau-(y-x)\cdot e)  \, d \tau \, d y
$$
$$
\ge
\int_{\R^N} J(x-y) e^{- \lambda(x-y)\cdot e} \int_{-\infty}^\infty \psi(\tau,y) e^{-\lambda \tau} \eta_n(\tau-1)  \, d \tau \, d y
$$
$$
=
\int_{\R^N} J(x-y) e^{- \lambda(x-y)\cdot e} \int_{-\infty}^\infty \psi(\tau+1,y) e^{-\lambda (\tau+1)} \eta_n(\tau)  \, d \tau \, d y
$$
$$
\ge
e^{-\lambda}
\int_{\R^N} J(x-y) e^{- \lambda(x-y)\cdot e} \int_{-\infty}^\infty \psi(\tau,y) e^{-\lambda \tau} \eta_n(\tau)  \, d \tau \, d y
$$
because $\psi(\cdot,x)$ is nondecreasing.
It follows that
$$
D_n \le (1 - e^{-\lambda}) J_\lambda U_n(\cdot,\lambda).
$$
Thus, from \eqref{eq 2} and since $F_n\le 0$
$$
\eps \Delta U_n + J_\lambda U_n - U_n + f_u(x,0) U_n  - c \lambda U_n \le   (1 - e^{-\lambda}) J_\lambda U_n(\cdot,\lambda)+E_n
$$
Write
$$
E_n =\int_{-\infty}^0 \ldots \, d s
+  \int_{0}^\infty \ldots \, d s
$$
and note that
$$
\int_{0}^\infty \change {\left|( f(x,\psi(s,x)) - f_u(x,0) \psi(s,x) ) e^{-\lambda s} \eta_n(s) \right| }\, d s \le C_1
$$
\co{put absolute value}
with $C_1 \sim 1/\lambda$ as $\lambda\to0^+$.
We estimate the other integral as follows:
\begin{align*}
\int_{-\infty}^0 & ( f(x,\psi(s,x)) - f_u(x,0) \psi(s,x) ) e^{-\lambda s} \, d s
 \le C_f \int_{-\infty}^0   \psi(s,x)^2 e^{-\lambda s} \eta_n(s) \, d s
\\
&
\le C_f \delta \int_{-\infty}^0   \psi(s,x) e^{-\lambda s} \eta_n(s) \, d s \le C_f \delta U_n(x,\lambda)
\end{align*}
where $C_f$ is a constant that depends only on $f$.

In this way we obtain
\begin{align}
\label{eq U}
\eps \Delta U_n + J_\lambda U_n - U_n + f_u(x,0) U_n  - c \lambda U_n \le   (1 - e^{-\lambda}) J_\lambda U_n(\cdot,\lambda)+ C_f \delta U_n + C_1
\end{align}
Let $\mu_{\eps,\lambda}$ be the principal eigenvalue of the operator $-(\eps \Delta \phi + J_\lambda \phi - \phi + f_u(x,0) \phi )$, $\phi_{\eps,\lambda}$, the principal eigenfunction and $\phi_{\eps,\lambda}^*$ be the principal eigenfunction for the adjoint operator.
Since $\mu_{\eps,\lambda}\to \mu_\lambda$ as $\eps\to0$ and $\mu_\lambda<0$, we can assume that $\mu_{\eps,\lambda}<0$.
Multiplying \eqref{eq U} by $\phi_{\eps,\lambda}^*$ and integrating over the period $[0,1]^N$ we find
$$
(-\mu_{\eps,\lambda} - c \lambda) \int_{[0,1]^N} U_n(x,\lambda) \phi_{\eps,\lambda}^*(x) \, d x
\le  (1 - e^{-\lambda}) \int_{[0,1]^N} J_\lambda U_n(x,\lambda) \phi_{\eps,\lambda}^* (x) \, d x +
$$
$$
+  C_f \delta \int_{[0,1]^N} U_n(x,\lambda) \phi_{\eps,\lambda}^*(x) \, d x
 + C_1 \int_{[0,1]^N}  \phi_{\eps,\lambda}^*(x) \, d x
$$
But
$$
\int_{[0,1]^N} J_\lambda U_n(x,\lambda) \phi_{\eps,\lambda}^* (x) \, d x =
\int_{[0,1]^N} (J_\lambda)^* \phi_{\eps,\lambda}^* (x) U_n(x,\lambda)  \, d x
$$
$$
= \int_{[0,1]^N} \left[ -\mu_{\eps,\lambda} \phi_{\eps,\lambda}^* + \phi_{\eps,\lambda}^* - f_u(x,0) \phi_{\eps,\lambda}^* - \eps \Delta \phi_{\eps,\lambda}^* \right] U_n(x,\lambda)  \, d x
$$
\change{Note that $\phi_{\eps,\lambda}^*$ is uniformly bounded in $C^2([0,1]^N)$ as $\eps\to0$, see Remark~\ref{remark smoothness eigenf}, a property where use that $f$ is $C^3$.
Using the uniform smoothness of $\phi_{\eps,\lambda}^*$} and the fact that it is uniformly bounded below $\phi_{\eps,\lambda}^*(x) \ge c>0$ as $\eps\to 0$ with $\lambda>0$ fixed, we see that
\begin{align*}
\int_{[0,1]^N} J_\lambda U_n(x,\lambda) \phi_{\eps,\lambda}^* (x) \, d x \le
C  \int_{[0,1]^N} U_n(x,\lambda) \phi_{\eps,\lambda}^*(x) \, d x.
\end{align*}
Therefore
$$
(-\mu_{\eps,\lambda} - c \lambda) \int_{[0,1]^N} U_n(x,\lambda) \phi_{\eps,\lambda}^*(x) \, d x \le  (  (1-e^{-\lambda}) C + C_f \delta )\int_{[0,1]^N} U_n(x,\lambda) \phi_{\eps,\lambda}^*(x) \, d x
$$
$$ + C_1 \int_{[0,1]^N}  \phi_{\eps,\lambda}^*(x) \, d x
$$
Choosing $\delta>0$ and $\lambda>0$ sufficiently small we deduce that
$$
\int_{[0,1]^N} U_n(x,\lambda) \phi_{\eps,\lambda}^*(x) \, d x \le C
$$
and again using that $\phi_{\eps,\lambda}^*$ is uniformly bounded below, we find
\begin{align}
\label{eq lemma a2}
\int_{[0,1]^N} U_n(x,\lambda) \, d x \le C
\end{align}
where $C$ is independent of $\eps$ and $n$. Now letting $n\to\infty$, we obtain the conclusion \eqref{eq lemma a}.

To prove the last part we observe that
$$
\lim_{n \to \infty} U_n(x,\lambda) = U(x,\lambda) $$
by monotone convergence where
$$
U(x,\lambda) = \int_{-\infty}^\infty \psi(s,x) e^{-\lambda s} \, d s  .
$$
By \eqref{eq lemma a2}, $U(\cdot,\lambda)$ is in $L^1([0,1]^N)$ and is a weak solution of

$$
\eps \Delta U + J_\lambda U - U - c \lambda U = \tilde E  \quad \hbox{in } \R^N
$$
where
$$
\tilde E  = \int_{-\infty}^\infty  f(x,\psi(s,x))  e^{-\lambda s}  \, d s .
$$
Note that
$$
\| \tilde E \|_{L^p([0,1]^N)} \le C \| U (\cdot,\lambda)\|_{L^p([0,1]^N)}
$$
for all $p\ge 1$. Then, using standard elliptic $L^p$ estimates we deduce that $U(\cdot,\lambda)\in L^\infty$ for $0<\lambda\le \lambda_0$.
\qed

\medskip

\begin{lemma}
\label{trivial lemma}
Suppose $\psi:(-\infty,0] \to [0,\infty)$ is nondecreasing and let $\lambda \in \R$.
Then
\begin{align}
\label{eq c}
\psi(s)  \le \lambda \frac{e^{\lambda s}}{1 - e^{\lambda s}} \int_{-\infty}^0 \psi(\tau) e^{-\lambda \tau} \, d \tau \quad \forall s \le 0.
\end{align}
\end{lemma}
\noindent
{\bf Proof.}
Let $t\le 0$. Then
$$
\psi(t) \int_t^0 e^{-\lambda s} \, d s  \le \int_t^0 \psi(s) e^{-\lambda s } \, d s .
$$
\qed

We prove first the exponential decay of $\psi$ for some constant that depends on $\eps$.
\begin{lemma}
\label{lemma b}
For any $\lambda<\lambda_\eps(c)$ there is  $C_\eps>0$  such that if $\psi$ is a solution of \eqref{eq a} then
\begin{align}
\label{eq lemma b}
\psi(s,x)  \le C_\eps e^{\lambda s}  \quad  \forall x \in \R^N , \
\forall s \in \R.
\end{align}
\end{lemma}
\noindent
{\bf Proof.}
In this proof $\eps>0$ is fixed and we find $\delta_\eps>0$ such that if $\psi$ satisfies
\begin{align}
\label{norm2}
\max_{x\in [0,1]^N} \psi(0,x) \le  \delta_\eps
\end{align}
then the conclusion \eqref{eq lemma b} holds. Given any solution of \eqref{eq a} we know already by Lemma~\ref{lemma a} that $\psi(s,x) \to 0$ as to $-\infty$ uniformly in $x$, even at an exponential rate, so that \eqref{norm2} holds provided we replace $\psi(x,s)$ by $\psi(x,s-\tau)$ with $\tau$ sufficiently large.

Let $\eta \in C^\infty(\R)$ be such that $\eta(t) =1$ for $t\le 1$ and $\eta(t)=0$ for $t\ge2 $.
For $\lambda \in \R$, $x\in [0,1]^N$, let $U$ be defined by
\begin{align}
\label{def U}
U(x,\lambda) = \int_{-\infty}^\infty \psi(s,x) e^{-\lambda s} \eta(s) \, d s
\end{align}
with values in $[0,\infty]$.
At this moment we  know from Lemma~\ref{lemma a} that $U(x,\lambda)<+\infty$ if we take $0<\lambda\le\lambda_0$ where $\lambda_0>0$ is \change{a} small fixed number. The objective is to prove that for any $\lambda$ such that $0<\lambda c < - \change{\mu_{\eps,\lambda}}$
$$
\|U(\cdot,\lambda)\|_{L^\infty([0,1]^N)} < + \infty.
$$
Then from \eqref{eq c} we obtain the desired conclusion.

Assume that $\lambda$ is such that $
\|U(\cdot,\lambda)\|_{L^\infty([0,1]^N)} < + \infty
$.
We multiply \eqref{eq a} by $\eta(s) e^{-\lambda s}$ and integrate on $(-\infty,\infty)$. We obtain
\begin{align*}
\eps \Delta U  + J_\lambda U - U + f_u(x,0) U - c \lambda U =
D_\lambda(x)+E_\lambda(x)+F_\lambda(x)
\end{align*}
where
\begin{align*}
D_\lambda(x) = \int_{\R^N} J(x-y) e^{- \lambda(x-y)\cdot e} \int_{-\infty}^\infty \psi(\tau,y) e^{-\lambda \tau} \left[  \eta(\tau) -\eta(\tau-(y-x)\cdot e)  \right] \, d \tau \, d y
\end{align*}
\begin{align*}
E_\lambda(x) = \int_{-\infty}^\infty ( f(x,\psi(s,x)) - f_u(x,0) \psi(s,x) ) e^{-\lambda s}  \eta(s) \, d s
\end{align*}
\begin{align*}
F_\lambda(x) = - c \int_{-\infty}^\infty \psi(s,x) \eta'(s) e^{-\lambda s} \, d s.
\end{align*}
Thus
$$
( L_{\eps,\lambda} - \lambda c) U  = D_\lambda +E_\lambda +F_\lambda .
$$
Since $U$ is nonnegative, we may apply Lemma~\ref{lemma L infty epsilon} and deduce
$$
\| U(\cdot,\lambda) \|_{L^\infty} \le C_{\eps,\lambda} (  \| D_\lambda +E_\lambda +F_\lambda \|_{L^\infty})
$$
Write $U= U_1 + U_2$ where
\begin{align}
\label{U1 U2}
U_1 = \int_{-\infty}^0 \psi(s,x) e^{-\lambda s} \eta(s) \, d s ,
\qquad
U_2 = \int_0^{\infty} \psi(s,x) e^{-\lambda s} \eta(s) \, d s.
\end{align}
Since $U_2\ge 0$, we also have
$$
\|U_1\|_{L^\infty([0,1]^N)} \le C_{\eps,\lambda}   \
\| D_\lambda+E_\lambda+F_\lambda\|_{L^\infty([0,1]^N)} .
$$
In $D_\lambda(x)$ one can restrict $\tau$ to \change{$[-1,4]$}.
\co{it does not matter much}
Hence
$$
\| D_\lambda \|_{L^\infty([0,1]^N)} \le C
$$
and the constant remains bounded as $\lambda$ varies in a bounded interval of $\R$. Similarly the integral in $F_\lambda(x)$ is restricted to $1\le \tau \le 2$ and hence
$$
\| F_\lambda \|_{L^\infty([0,1]^N)} \le C
$$
with $C$ as before.
We estimate
$$
|E_\lambda(x)| =  \left|\int_{-\infty}^\infty ( f(x,\psi(s,x)) - f_u(x,0) \psi(s,x) ) e^{-\lambda s} \eta(s) \, d s \right|
$$
$$
\le C\int_{-\infty}^{-1} |\psi(s,x)|^2 e^{-\lambda s} \, d s  + C
$$
By \eqref{eq c}
$$
|\psi(s,x)| \le C_0 e^{\lambda s}\|U_1(\cdot,\lambda)\|_{L^\infty} \quad \forall x \in [0,1]^N, \ \forall s \le -1.
$$
Hence, using \eqref{norm2},
$$
|E_\lambda(x)| \le C \delta_\eps^{1/2} \int_{-\infty}^{-1} |\psi(s,x)|^{3/2} e^{-\lambda s} \, d s  + C
$$
$$
\le C \delta_\eps^{1/2} \|U_1(\cdot,\lambda)\|_{L^\infty}^{3/2} \int_{-\infty}^{-1} e^{\lambda s/2} \, d s  + C = C_{\lambda_0} \delta_\eps^{1/2} \|U_1(\cdot,\lambda)\|_{L^\infty}^{3/2}  + C .
$$
where $C_{\lambda_0} \sim 1/\lambda_0$.
Therefore
\begin{align}
\label{est a13}
\|U_1(\cdot,\lambda)\|_{L^\infty([0,1]^N)} \le  \delta_\eps^{1/2} C_{\lambda_0}  C_{\eps,\lambda} \|U_1(\cdot,\lambda)\|_{L^\infty}^{3/2}  + C_1.
\end{align}
If we choose $\delta_\eps>0$ small this implies that there is a gap for $\|U_1(\cdot,\lambda)\|_{L^\infty([0,1]^N)}$. For example we can achieve
\begin{align*}
\hbox{either $\|U_1(\cdot,\lambda)\|_{L^\infty([0,1]^N)} \le 2 C_1$ or  $\|U_1(\cdot,\lambda)\|_{L^\infty([0,1]^N)} \ge 3 C_1$}.
\end{align*}
Indeed, first fix $0<\lambda_0<\lambda_1< \lambda_\eps(c)$. Then we know from Lemma~\ref{lemma L infty epsilon} that
$$
\sup_{\lambda_0\le\lambda\le\lambda_1} C_{\eps,\lambda}  <\infty.
$$
Choose $\delta_\eps>0$ such that
$$
\delta_\eps^{1/2}  (3 C_1)^{1/2} C_{\lambda_0} \left( \sup_{\lambda_0\le\lambda\le\lambda_1} C_{\eps,\lambda}  \right) \le \frac13.
$$
Suppose that $\|U_1(\cdot,\lambda)\|_{L^\infty([0,1]^N)} \le 3 C_1$. Then by \eqref{est a13}
\begin{align*}
\|U_1(\cdot,\lambda)\|_{L^\infty([0,1]^N)}
& \le  \delta_\eps^{1/2} C_{\lambda_0}   C_{\eps,\lambda}   \|U_1(\cdot,\lambda)\|_{L^\infty}^{3/2}  + C_1
\\
& \le \delta_\eps^{1/2} C_{\lambda_0}   C_{\eps,\lambda}  (3 C_1)^{1/2 }\|U_1(\cdot,\lambda)\|_{L^\infty}  + C_1
\\
& \le \frac13\|U_1(\cdot,\lambda)\|_{L^\infty}  + C_1 \le 2 C_1 .
\end{align*}
Using Lemma~\ref{lemma a} and increasing $C_1$ and decreasing $\delta_\eps$ if necessary, we can assume that
$$
\| U_1(\cdot,\lambda_0) \|_{L^\infty} \le 2C_1.
$$
Since $\lambda \mapsto \| U_1(\cdot,\lambda) \|_{L^\infty}$ is continuous we see that
$$
\| U_1(\cdot,\lambda) \|_{L^\infty} \le 2C_1 \quad \forall \lambda_0\le  \lambda \le \lambda_1.
$$
\qed

\medskip
\noindent
{\bf Proof of Proposition~\ref{prop exp decay}.}
We argue as in Lemma~\ref{lemma b}.
In this proof we take $\rho>0$ as in as in Proposition~\ref{prop uniform bound} and let $0<\rho'<\rho$. We restrict $\lambda$
\change{so that it satisfies $(-\mu_{\eps,\lambda}-\rho)/c \le \lambda \le (-\mu_{\eps,\lambda}-\rho')/c$} and take $0<\eps\le\eps_0$.

Let $U$ be defined by \eqref{def U}, and $U_1$, $U_2$ defined in \eqref{U1 U2}. Following the proof of Lemma~\ref{lemma b}, if $\psi$ satisfies \eqref{eq a} and \eqref{normalization phi} then,
\change{using Proposition~\ref{prop uniform bound},}
\begin{align*}
\|U_1(\cdot,\lambda)\|_{L^\infty([0,1]^N)} \le  \delta^{1/2} \change{\overline C}  \|U_1(\cdot,\lambda)\|_{L^\infty}^{3/2}  + C_1 ,
\end{align*}
where $\change{\overline C} $ now remains bounded
\change{for any $0<\eps\le \eps_0$
if $\lambda$ satisfies
$(-\mu_{\eps,\lambda}-\rho)/c \le \lambda \le (-\mu_{\eps,\lambda}-\rho')/c$.}
Again, choosing  $\delta>0$ small such that
$$
\delta^{1/2}  (3 C_1)^{1/2} \change{\overline C}  \le \frac13
$$
we obtain
\begin{align*}
\hbox{either $\|U_1(\cdot,\lambda)\|_{L^\infty([0,1]^N)} \le 2 C_1$ or  $\|U_1(\cdot,\lambda)\|_{L^\infty([0,1]^N)} \ge 3 C_1$}.
\end{align*}
Let $\psi_\tau(s,x)=\psi(s-\tau,x)$ where $\tau>0$ and $U_{1,\tau}$ denote the corresponding {\em Laplace transform} as in \eqref{def U}, \eqref{U1 U2}. By Lemma~\ref{lemma b}
$$
\| U_{1,\tau}(\cdot,\lambda) \|_{L^\infty} \to 0 \quad \hbox{as } \tau \to+\infty.
$$
Since $\tau \mapsto \| U_{1,\tau}(\cdot,\lambda) \|_{L^\infty}$ is continuous we see that
$$
\| U_{1,0}(\cdot,\lambda) \|_{L^\infty} \le 2C_1 .
$$
\change{Then by Lemma~\ref{trivial lemma} we obtain \eqref{decay phi}.}
\qed

\section{Proof of the main theorem}
\label{sec proof main thm}

In this section we prove Theorem~\ref{thm main}, by establishing a uniform estimate in $W^{1,p}_{loc}$ of $\psi_\eps$, the convergence of $\psi_\eps$ to a function $\psi$ satisfying the equation, and finally establishing that $\psi$ solves the full problem.

\begin{proposition}
\label{prop sobolev est}
There is $\delta>0$ such that if $ \psi_\eps$ is a solution of \eqref{eq a} satisfying the normalization condition
\eqref{normalization phi}, then for any
for any $1\le p < \infty$ and bounded open set $D$ in $ \R \times \R^N$ there is a constant $C$ independent of $\eps$ as $\eps\to 0$ such that:
\begin{align}
\label{W 1p estimate}
\| \psi_\eps\|_{W^{1,p}(D)} \le C.
\end{align}
\end{proposition}
\noindent
{\bf Proof.}
For simplicity we write $\psi=\psi_\eps$
\change{and we use the notation $\psi_{x_i} = \frac{\partial\psi}{\partial x_i}$.}
We differentiate the equation in \eqref{eq a} with respect to $x_i$
and get
\begin{align}
\label{derivada en x}
c \psi_{s x_i} =  \eps \Delta \psi_{x_i} + M_{x_i} [\psi] - e_i M [\psi_s] - \psi_{x_i} + f_u(x,\psi) \psi_{x_i} + f_{x_i}(x,\psi)
\end{align}
where
$$
M_{x_i} [\psi](s,x) = \int_{\R^N} J_{x_i}(x-y) \psi(s+(y-x)\cdot e,y) \, d y
$$
$e=(e_1,\ldots,e_N)$.
We write this as
\begin{align}
\nonumber
c \psi_{s x_i}  + (1- f_u(x,0) ) \psi_{x_i} & =
\eps \Delta \psi_{x_i} + M_{x_i} [\psi] - e_i M [\psi_s]
\\
\label{eq d ds}
& \qquad + ( f_u(x,\psi) - f_u(x,0) )\psi_{x_i} + f_{x_i}(x,\psi) .
\end{align}
Let $1\le p < +\infty$ and $\theta>0$ to be fixed later on.
Then
$$
\frac{\partial }{\partial s} \left( e^{s p (1-f_u(x,0) -\theta)/c} |\psi_{x_i}|^p\right)
=
\frac{p}{c} e^{s p(1-f_u(x,0) -\theta)/c} \left( c \psi_{s x_i}  + (1-f_u(x,0) - \theta) \psi_{x_i} \right)
|\psi_{x_i}|^{p-2} \psi_{x_i}
$$
Using \eqref{eq d ds} we obtain
\begin{align*}
\frac{\partial }{\partial s} \left( e^{s p (1-f_u(x,0) -\theta)/c} |\psi_{x_i}|^p
\right)
& = \frac{p}{c} e^{s p(1-f_u(x,0) -\theta)/c} \Big( \eps \Delta \psi_{x_i} + M_{x_i} [\psi] - e_i M [\psi_s]
\\
& \quad   + ( f_u(x,\psi) - f_u(x,0) )\psi_{x_i} + f_{x_i}(x,\psi)   - \theta \psi_{x_i} \Big)
|\psi_{x_i}|^{p-2} \psi_{x_i} .
\end{align*}
We integrate now with respect to $x$ over the period $[0,1]^N$ and estimate the terms on the right hand side.
\begin{align*}
\frac{c}{p}  \frac{\partial }{\partial s}  &
\int_{[0,1]^N} e^{s p(1-f_u(x,0) -\theta)/c}
|\psi_{x_i}|^p
  \, d x
= I_1 + I_2 + I_3 + I_4 + I_5 + I_6
\end{align*}
where
\begin{align*}
I_1 & = \eps \int_{[0,1]^N} e^{s p(1-f_u(x,0) -\theta)/c}
\Delta \psi_{x_i}
|\psi_{x_i}|^{p-2} \psi_{x_i}
\, d x
\\
I_2 & = \int_{[0,1]^N} e^{s p(1-f_u(x,0) -\theta)/c}
M_{x_i}[ \psi ] |\psi_{x_i}|^{p-2} \psi_{x_i}
\, d x
\\
I_3 & = - e_i \int_{[0,1]^N} e^{s p(1-f_u(x,0) -\theta)/c}
M[ \psi_s ] |\psi_{x_i}|^{p-2} \psi_{x_i}
\, d x
\\
I_4 & =  \int_{[0,1]^N} e^{s p(1-f_u(x,0) -\theta)/c}   ( f_u(x,\psi) - f_u(x,0) )
|\psi_{x_i}|^{p}
\, d x
\\
I_5 & = \int_{[0,1]^N} e^{s p(1-f_u(x,0) -\theta)/c}
f_{x_i}(x,\psi)
|\psi_{x_i}|^{p-2} \psi_{x_i}
\, d x
\\
I_6 & = - \theta \int_{[0,1]^N} e^{s p(1-f_u(x,0) -\theta)/c}  |\psi_{x_i}|^{p}
\, d x .
\end{align*}
Integrating by parts we can estimate
\begin{align*}
I_1 & =
- \eps (p-1)
\int_{[0,1]^N} e^{s p(1-f_u(x,0) -\theta)/c} |\psi_{x_i}|^{p-2}
|\nabla \psi_{x_i}|^2  \, d x
\\
& \quad
- \eps \int_{[0,1]^N}
\nabla
\left(  e^{s p(1-f_u(x,0) -\theta)/c} \right)
\nabla \psi_{x_i}
|\psi_{x_i}|^{p-2} \psi_{x_i}
\, d x
\\
& \le  \frac{\eps |s| p }{c}
\int_{[0,1]^N} e^{s p(1-f_u(x,0) -\theta)/c}
|\nabla_x f_{u}(x,0)| \,
\change{|\nabla \psi_{x_i} | \,
|\psi_{x_i}|^{p-1} }
\, d x .
\end{align*}
By Young's inequality
\begin{align*}
I_1
& \le \frac{\theta}{5}
\int_{[0,1]^N}
e^{s p(1-f_u(x,0) -\theta)/c}
|\psi_{x_i}|^{p}
\, d x
\\
& \quad
+
C \eps^p |s|^p
\int_{[0,1]^N}  e^{s p(1-f_u(x,0) -\theta)/c}
| \nabla \psi_{x_i} |^p \, d x
\end{align*}
where $C$ depends on $\theta $ and $\| f\|_{C^2}$.
In a similar way
\begin{align*}
I_2
& \le
\frac{\theta}{5}
\int_{[0,1]^N} e^{s p(1-f_u(x,0) -\theta)/c}
|\psi_{x_i}|^{p}
\, d x +
C
\int_{[0,1]^N} e^{s p(1-f_u(x,0) -\theta)/c}
|M_{x_i}[ \psi ]|^p
\, d x
\\
I_3
& \le
\frac{\theta}{5}
\int_{[0,1]^N} e^{s p(1-f_u(x,0) -\theta)/c}
|\psi_{x_i}|^{p}
\, d x +
C
\int_{[0,1]^N} e^{s p(1-f_u(x,0) -\theta)/c}
|M[\psi_s]|^p
\, d x
\\
I_5
&
\le \frac{\theta}{5}
\int_{[0,1]^N} e^{s p(1-f_u(x,0) -\theta)/c}
|\psi_{x_i}|^{p} \, d x  +
C
\int_{[0,1]^N} e^{s p(1-f_u(x,0) -\theta)/c}
|f_{x_i}(x,\psi)|^p \, d x
\end{align*}
To estimate $I_4$ we write
$$
I_4
 \le
\sup_{y}  | f_u(y,\psi(s,y))) - f_u(y,0) |
\int_{[0,1]^N} e^{s p(1-f_u(x,0) -\theta)/c}
|\psi_{x_i}|^{p}
\, d x .
$$
We work with $\delta>0$ small so that from the normalization condition \eqref{normalization phi} we get
$$
\sup_{y}  | f_u(y,\psi(s,y))) - f_u(y,0) |
\le \frac{\theta}{5}\quad \hbox{for all } s \le 0.
$$
Then
\begin{align*}
I_4 \le\frac{\theta}{5}
\int_{[0,1]^N} e^{s p(1-f_u(x,0) -\theta)/c}
|\psi_{x_i}|^{p}
\, d x
\end{align*}
Combining the previous estimates we obtain
\begin{align}
\label{ineq derivative}
\frac{c}{p}  \frac{\partial }{\partial s}
&\int_{[0,1]^N} e^{s p(1-f_u(x,0) -\theta)/c}
|\psi_{x_i}|^p
\, d x
\\
\nonumber
& \le
C \eps^p |s|^p
\int_{[0,1]^N}  e^{s p(1-f_u(x,0) -\theta)/c}
| \nabla \psi_{x_i} |^p \, d x
\\
\nonumber
& \quad
+ C
\int_{[0,1]^N} e^{s p(1-f_u(x,0) -\theta)/c}
|M_{x_i}[ \psi ]|^p
\, d x
\\
\nonumber
& \quad +
C
\int_{[0,1]^N} e^{s p(1-f_u(x,0) -\theta)/c}
|M[\psi_s]|^p
\, d x
\\
\nonumber
& \quad +
C
\int_{[0,1]^N} e^{s p(1-f_u(x,0) -\theta)/c}
|f_{x_i}(x,\psi)|^p \, d x
\end{align}
Let $t_0\le t \le 0$. We integrate with respect to $s$ over $[t_0,t]$ and then let $t_0 \to - \infty$. By \eqref{decay phix}, given any $0<\lambda < \lambda_\eps(c)$ there is $C$ such that
\begin{align}
\label{exp dec ineq 1}
\nonumber
&
\int_{[0,1]^N} e^{s p(1-f_u(x,0) -\theta)/c}
|\psi_{x_i}(s,x)|^p
\, d x
\\ & \qquad
\le \frac{C}{\eps^{p/2}}
\int_{[0,1]^N}
\exp ( s p(1-f_u(x,0) -\theta + \lambda c )/c )
\, d x .
\end{align}
We choose now $\lambda$ and $\theta$ as follows.
\change{We fix a large $\Lambda_0>0$.}
We note that
since there is a principal periodic eigenfunction $\phi_\lambda \in \change{C_{per}(\R^N)}$, $\phi_\lambda>0$ for
$$
J_\lambda*\phi_\lambda - \phi_\lambda
+f_u(x,0 ) \phi_\lambda + \mu_\lambda \phi_\lambda = 0
\quad \text{in } \R^N
$$
we must have
$$
\gamma \equiv
\change{\inf_{\lambda \in [0,\Lambda_0]}}
\inf_{x\in\R^N } ( 1-f_u(x,0) -\mu_\lambda)
=
\change{\inf_{\lambda \in [0,\Lambda_0]}}
\inf_{x\in\R^N }  \change{\frac{J_\lambda*\phi_\lambda  (x)}{\phi_\lambda(x)}}
>0.
$$
Since $\mu_{\eps,\lambda} \to \mu_\lambda $ as $\eps \to 0$, for $\eps>0$ sufficiently small
$$
\inf_{x\in\R^N } ( 1-f_u(x,0) -\mu_{\eps,\lambda})  \ge \gamma/2 >0.
$$
and since for $\lambda = \lambda_\eps(c)$ we have $\lambda c = -\mu_{\eps,\lambda}$ we get
$$
\lambda_\eps(c) \ge \frac{\gamma}{2c } + \sup_{x\in \R^N} \frac{f_u(x,0) - 1}{c}  .
$$
Take $\lambda>0$ such that
\begin{align}
\label{choice lambda}
\sup_{x\in \R^N} \frac{f_u(x,0) - 1}{c}
+\frac{\gamma}{4c}
\change{ \leq } \lambda \le \lambda_\eps(c) - \frac{\gamma}{4c}.
\end{align}
Then choose $\theta= \gamma/8>0$ and get
\begin{align}
\label{def sigma}
\sigma \equiv \inf_{x\in\R^N} \left(  \frac{1-f_u(x,0) -\theta}{c} + \lambda \right) >0 .
\end{align}
Then from \eqref{exp dec ineq 1} we obtain
$$
\int_{[0,1]^N} e^{s p(1-f_u(x,0) -\theta)/c}
|\psi_{x_i}(s,x)|^p
\, d x
\le \frac{C}{\eps^{p/2}} e^{p \sigma s },
\quad
\change{\quad \forall s\leq 0,}
$$
and therefore
\begin{align}
\label{limit int}
\lim_{s\to-\infty}
\int_{[0,1]^N} e^{s p(1-f_u(x,0) -\theta)/c}
|\psi_{x_i}(s,x)|^p
\, d x
=0.
\end{align}
Integrating \eqref{ineq derivative} in $[t_0,t]$ with $t_0\le t \le 0$ and using \eqref{limit int}
we obtain
\begin{align}
\label{after int}
\frac{c}{p}
\int_{[0,1]^N} e^{s p(1-f_u(x,0) -\theta)/c}
|\psi_{x_i}|^p
\, d x \le  K_1 + K_2 + K_3 + K_4
\end{align}
where
\begin{align*}
K_1 & = C \eps^p
\int_{-\infty}^t |s|^p
\int_{[0,1]^N}
e^{s p(1-f_u(x,0) -\theta)/c}
| \nabla \psi_{x_i} |^p \, d x
\, d s
\\
K_2 & = C
\int_{-\infty}^t
\int_{[0,1]^N} e^{s p(1-f_u(x,0) -\theta)/c}
|M_{x_i}[ \psi ]|^p
\, d x
\, d s
\\
K_3 & = C
\int_{-\infty}^t
\int_{[0,1]^N} e^{s p(1-f_u(x,0) -\theta)/c}
|M[\psi_s]|^p
\, d x
\, d s
\\
K_4 & =  C
\int_{-\infty}^t
\int_{[0,1]^N} e^{s p(1-f_u(x,0) -\theta)/c}
|f_{x_i}(x,\psi)|^p \, d x
\, d s
\end{align*}
Next we claim that $K_1,K_2,K_3,K_4$ remain bounded as $\eps \to0$.
Indeed, by \eqref{decay phixx} and \eqref{def sigma},
\begin{align*}
e^{s p(1-f_u(x,0) -\theta)/c}
| \nabla \psi_{x_i} |^p
& \le
\change{ e^{s p(1-f_u(x,0) -\theta)/c}
| \nabla^2_x \psi |^p }
\\
& \le
\frac{C}{\eps^p}
e^{s p(1-f_u(x,0) -\theta + \lambda c )/c}
\le \frac{C}{\eps^p}
e^{s p\sigma } ,
\end{align*}
for \change{$s\leq 0$}, $x \in \R^N$
with $C$ independent of $\eps$
\change{(note  that $\nabla \psi_{x_i}$ is a second order derivative  of $\psi$)}.
Therefore $K_1$ is bounded as $\eps\to0$. The other ones can be bounded similarly, using \eqref{decay phi}, \eqref{decay phis} and the hypotheses $f(x,0)=0$, $f\in C^3$ which imply
$$
|f_{x_i}(x,u)|\le  C u
\quad\text{for } 0\le u \le \delta
$$
for some $C$. Thus from \eqref{after int} we deduce that there exists $C$ independent of $\eps$ for $\eps$ small such that for all $s\le 0$
\begin{align}
\label{uniform s le 0}
\int_{[0,1]^N} e^{s p(1-f_u(x,0) -\theta)/c}
|\psi_{x_i}(s,x)|^p
\, d x \le  C
\end{align}
This together with \eqref{decay phis} proves the estimate \eqref{W 1p estimate} for any bounded open set $D\subset (-\infty,0) \times \R^N$. To obtain \eqref{W 1p estimate} for any bounded open set $D\subset \R \times \R^N$ we proceed similarly as before. We multiply \eqref{derivada en x} by $|\psi_{x_i}|^{p-2} \psi_{x_i} $ and integrate over $[0,1]^N$.  Using that $\psi$ has a uniform upper bound we obtain
$$
\frac{d}{d s}
\int_{[0,1]^N}
|\psi_{x_i}|^p \, d x \le
C \int_{[0,1]^N} |\psi_{x_i}|^p \, d x.
$$
Then, using Gronwall's inequality we deduce for $s\ge 0$
$$
\int_{[0,1]^N}
|\psi_{x_i}(s,x)|^p
\, d x \le
e^{C s}
\int_{[0,1]^N}
|\psi_{x_i}(0,x)|^p
\, d x
+ C.
$$
Since by  \eqref{uniform s le 0} we have  a uniform control of the form
$\int_{[0,1]^N}
|\psi_{x_i}(0,x)|^p
\, d x \le C $,
we obtain that for all $R>0$ there exists $C>0$ independent of $\eps$ such that
$$
\int_{[0,1]^N}
|\psi_{x_i}(s,x)|^p
\, d x \le   C
\quad \text{for all } |s|\le R.
$$
Using this and  \eqref{decay phis} we obtain the estimate \eqref{W 1p estimate} for any bounded open set $D\subset \R \times \R^N$.
\qed

\begin{lemma}
\label{lemma convergence}
If $c\ge \change{c_e^*}$  there exists a function $\psi:\R \times \R^N$ which is $C^1$ in $s$ and Lipschitz continuous and satisfies
\begin{align}
\label{eq psi}
c \psi_s  =  M [\psi] - \psi + f(x,\psi) \quad \forall s\in \R,\ x\in\R^N
\end{align}
and
$$
\lim_{s\to-\infty} \psi(s,x) =0.
$$
Furthermore $\psi>0$ is periodic in $x$ and non-decreasing in $s$.
\end{lemma}
\noindent{\bf Proof.}
Let $c \ge \change{c_e^*}$. If $c>\change{c_e^*}$ then $c>\change{c_e^*(\eps)}$ for $\eps>0$ small and we let, for small $\eps>0$, $\psi_\eps$ be the solution constructed in
Proposition~\ref{pro aprox} with speed $c$. If  $c = \change{c_e^*}$  we let
 $\psi_\eps$ be the solution constructed in
Proposition~\ref{pro aprox} with speed $c_\eps=\change{c_e^*(\eps)}$. In any case we have a solution of
\eqref{pro aprox} with speed $c_\eps \to c$, satisfying also \eqref{cond aprox}.

Let $\delta >0$ be from Lemma~\ref{prop sobolev est} and shift in $s$ so that $\psi_\eps$ satisfies
$$
\max_{x\in [0,1]^N} \psi_\eps(0,x) = \delta .
$$
Then, choosing $p>N$ in Lemma~\ref{prop sobolev est}  we can find a sequence $\eps_n \to 0$ such that $\psi_{\eps_n} \to \psi$ uniformly \change{on} compact sets.  Using this local uniform convergence we see that the function  $\psi $ satisfies \eqref{eq psi} in the following weak form
$$
- c \int_{-\infty}^\infty
\int_{[0,1]^N} \psi \varphi_s \, d x d s
=
 \int_{-\infty}^\infty
\int_{[0,1]^N}
( M [\psi] - \psi + f(x,\psi) ) \varphi \, d x d s
$$
for all $\varphi :\R \times \R^N \to \R$  smooth periodic function with compact support. This implies that $\psi$ is $C^1$ in $s$ and satisfies \eqref{eq psi} classically. Since $\psi_\eps$ is non-decreasing in $s$ and periodic in $x$ we deduce
that $\psi$ is also
non-decreasing in $s$ and periodic in $x$.
Moreover, by Proposition~\ref{prop exp decay}, if we take $0<\lambda<\lambda_c$ we have $\psi_\eps(s,x) \le C e^{\lambda s}$ with $C$ independent of $\eps$. Letting $\eps\to0$ we find the same inequality for $\psi$ and hence $\lim_{s\to-\infty} \psi(s,x) =0 $.

\change{
Finally, we prove that $\psi$ is Lipschitz continuous, which follows the same lines of Proposition~\ref{pro aprox}, so we point out the main steps.
\co{added proof}
Let $b_i$, $i=1,\dots,N$ denote the canonical basis in $\R^N$. Given $h\in \R$ we define
$$
D_i^h \psi(s,x) = \frac{\psi(s,x+b_i h) - \psi(s,x)}{h} .
$$
We choose $\lambda,\theta,\sigma>0$ as in \eqref{choice lambda}, \eqref{def sigma} so that
\begin{align}
\label{decay exp}
e^{2 s (1-f_u(x,0) -\theta)/c}
\le e^{2 s (\sigma-\lambda)}
\quad\forall x\in \R^N, \ s \le 0.
\end{align}
Then we compute
\begin{align*}
&
\frac{\partial}{\partial s}
\left(
e^{2 s (1-f_u(x,0) -\theta)/c} (D_i^h \psi)^2
\right) =
\\
&=\frac{2}{c}
e^{2 s (1-f_u(x,0) -\theta)/c}
\Big(
M_i[\psi^h] -e_i  M[D_s^{-h e_i} \psi] + (f_u(x,\tilde \psi) - f_u(x,0)) D_i^h \psi
\\
& \qquad\qquad
+ D_i^h f(\cdot,\psi(s,x+b_i h)) - \theta D_i^h \psi
\Big) D_i^h \psi
\end{align*}
where $e = (e_1,\ldots,e_N)$,
\begin{align*}
M_i[g](s,x)
&= \int_{\R^N}
\frac{J(x+b_i h - y)- J(x-y)}{h} g(s + (y-x)\cdot e,y) \, d y
\\
\psi^h(s,x)
&= \psi(s- e_i h,x)
\\
D_s^\tau \psi(s,x)
&= \frac{\psi(s+\tau,x) - \psi(s,x)}{\tau} ,
\end{align*}
and $\tilde \psi(s,x)$ lies between $\psi(s,x)$ and $\psi(s,x+b_i h)$. From here we deduce
\begin{align*}
&
\frac{\partial}{\partial s}
\left(
e^{2 s (1-f_u(x,0) -\theta)/c} (D_i^h \psi)^2
\right) \le
\\
&
\qquad\qquad\qquad
e^{2 s (1-f_u(x,0) -\theta)/c}
\Big(
M_i[\psi^h]^2
+M[D_s^{-e_i h} \psi]^2
+( D_i^h f(\cdot,\psi(s,x+b_i h)) )^2
\Big) .
\end{align*}
Using the exponential decay $\psi(s,x) \le C e^{\lambda s}$ for all $s\le 0$ and all $x\in\R^N$, and a similar one for $\psi_s$ (c.f. \eqref{decay phis}), we deduce from this and \eqref{decay exp} that
$$
\frac{\partial}{\partial s}
\left(
e^{2 s (1-f_u(x,0) -\theta)/c} (D_i^h \psi)^2
\right)
\le C e^{2 \sigma s}  .
$$
Integrating from $-\infty$ to $s\le 0$, we conclude that
there exists $C$ independent of $h$ such that
$$
|D_i^h \psi (s,x)) |\le C e^{ \lambda s}, \quad \forall x\in\R^N, \
\forall s\leq 0.
$$
This proves that $\psi(s,\cdot)$ is Lipschitz continuous for all $s\leq 0$. An argument similar to the one at the end of Proposition~\ref{pro aprox}  shows that it is also Lipschitz continuous for all $s\in \R$.}
\qed
\bigskip

We now prove the exponential convergence $\psi(s,x) \to p(x)$ as $s\to+\infty$, uniformly in $x$, by constructing appropriate subsolutions.
\begin{lemma}
Let $\psi$ be the function constructed in Lemma~\ref{lemma convergence}.
Then there exists $C$, $\sigma>0$ such that
$$
0\le p(x) - \psi(s,x)  \le C e^{-\sigma s}
\quad \text{for all } s\ge 0.
$$
In particular
$$
\lim_{s\to+\infty} \psi(s,x) = p(x) \quad
\text{uniformly for } x \in \R^N.
$$
\end{lemma}
\noindent
{\bf Proof.}
First we note that
$$
\psi(s,x) \le p(x) \quad
\text{for all } s \in \R, \ x \in \R^N.
$$

Next we show that $\psi(s,x) \to p(x)$ as $s\to+\infty$ uniformly for $x\in \R^N$. For this we will prove that there exists $\eps_0>0$ such that for any $0<m_0<1$ there is $s_0 \in \R$ such that
\begin{align}
\label{prop m0}
\psi_\eps(s,x) \ge m_0 p_\eps(x)
\quad
\text{for all } x\in \R^N,\  s \ge s_0, \ 0<\eps\le \eps_0.
\end{align}
The value $s_0$ depends on $m_0$ but not on $\eps$.

Recall that we have normalized $\psi_\eps$ by
$$
\max_{x\in [0,1]^N} \psi_\eps(0,x) = \delta
$$
where $\delta >0$ is from Proposition~\ref{prop sobolev est}.
By Lemma~\ref{lemma convergence}
$$
\psi_\eps \to \psi \quad \text{as } \eps \to 0
$$
uniformly on compact sets of $ \R \times \R^N $.
Since $\psi>0$ in $\R^N \times \R$  and is continuous
we see that that there is $\eps_0>0$ and $a>0$ such that for $0<\eps\le\eps_0$
$$
\psi_\eps (0,x) \ge 2 a p_\eps(x)
\quad \forall x\in \R^N .
$$
Note that  $a<1$.
Then we also have
\begin{align*}
\psi_\eps (s,x) \ge 2 a p_\eps(x)
\quad \forall x\in \R^N,\ s \ge 0,
\end{align*}
because $\psi_\eps(\cdot,x)$ is non-decreasing.

Given $a\le m\le 1$, $R\ge 1$, we construct a family of functions
$$
v_m(s,x) = \lambda_m(s) p_\eps(x) \quad
s \in \R, \ x \in \R^N
$$
where
$$
\lambda_m(s) = a + \frac{(m-a) s }{R+1} (1-\eta(s-R) ) + (m-a) \eta(s-R)
$$
and $\eta \in C^\infty(\R)$ is a cut-off function such that $\eta(s)= 0$ for $s\le 0$, $\eta(s) = 1$ for $s \ge 1$, $0\le\eta\le 1$ and $0\le \eta' \le 2$.
Note that $a\le \lambda_m(s) \le m $
\change{for all $s\geq0$}.

Fix $0< m_0 <1$ and let $a\le m \le m_0$.
It can be shown that we can choose $R>0$ large enough, independently of $\eps$, so that $v_m$ satisfies
$$
\eps \Delta v_m + M [v_m] - v_m + f(x,v_m) - c (v_m)_s  \ge 0
$$
for $s \ge 1$ and $x\in \R^N$.

Using a sliding argument we obtain that $a \le m \le m_0$
\begin{align*}
\psi_\eps \ge v_{m}
\quad
\text{for all } s \ge 1, \ x\in [0,1]^N.
\end{align*}

Using this inequality with $m=m_0$ we establish \eqref{prop m0}. Letting $\eps\to0$ we the deduce that
$$
\lim_{s\to+\infty}\psi(s,x) = p(x)
\quad\text{uniformly for } x\in \R^N.
$$

Finally, let us show that there is exponential convergence. For this we construct a subsolution $w_m$ with this property. Indeed,
let $\sigma>0$ to be fixed shortly and  $0\le m \le 1$.
We set
$$
w_m(s,x) = m ( 1- e ^{-\sigma s} ) p(x).
$$

Choosing  $S_0$ large and  $\sigma >0$ small we obtain that
$$
M[w_m] -  w_m+  f(x,w_m) - c (w_m)_s \ge
0
\quad
\text{in } [S_0,+\infty) \times \R^N.
$$
Let $S_1$ be  such that
$$
\psi(s,x) \ge (1-e^{-\sigma (S_0+1) } ) p(x)
 \quad\forall s \ge S_1 , \ x \in \R^N.
$$
This can be done because we know that $\psi(s,x) \to p(x)$ as $s\to+\infty$ uniformly for $x \in \R^N$.

Using again a  sliding argument we can prove that
\begin{align*}
\psi(s,x) \ge w_m(s+S_0-S_1,x)
\quad\forall s \ge S_1 , \ x \in \R^N
\end{align*}
and all $0\le m < 1$. Letting $m \to 1$ we find
$$
\psi(s,x) \ge (1 - e^{-\sigma(s + S_0 - S_1)} ) p(x)
\quad\text{for all } s \ge s_0, \  x\in \R^N,
$$
which finishes the proof of the lemma.
\qed

\begin{remark}
\change{
The limit
$
\tilde p(x) = \lim_{s\to\infty} \psi(s,x)
$
exists by monotonicity, but we cannot assert that it defines a continuous function (we have not proved uniform continuity of $\psi(s,x)$ as $s\to\infty$).
One could then argue that $\tilde p$ is a bounded measurable solution of the stationary problem and that Theorem~\ref{thm 1} also asserts the uniqueness of this solution. This would yield pointwise convergence
$\lim_{s\to+\infty} \psi(s,x) = p(x) $ for all $x \in \R^N$.
}
\end{remark}

Lastly, to finish the proof of  Theorem \ref{thm main} we prove the non-existence of front for speed $c< \change{c_e^*}$.
\begin{lemma}
Let $J$ and $f$ satisfy \eqref{hyp J} and \eqref{hyp f1} and let  $e \in \R^N$ be a unit vector. Assume $\mu_0<0$ and that there exists $\phi \in C_{per}(\R^N)$, $\phi>0$ satisfying \eqref{spectral-prob}. Then there exists no pulsating front $(\psi,c)$  connecting $0$ and $p(x)$ in the direction $e$ so that $c<\change{c_e^*}$.
\end{lemma}
\dem{Proof.}
Assume by contradiction that there exists a pulsating front $\psi$ with speed $c<\change{c_e^*}$. Then up to a shift $\psi$ is a supersolution of the parabolic problem \eqref{ptw.eq.para-intro} for any initial data $u_0\ge 0$ so that
$$ \sup_{\R^N}u_0<\min_{\R^N}p(x),\; \liminf_{r\to +\infty} \inf_{x.e\le r} u_0>0, \; u_0=0 \;\text{ for } x.e <<-1  $$
Let $u$ be the solution of the parabolic problem \eqref{ptw.eq.para-intro} with initial data $u_0$ satisfying the above condition then by the maximum principle, we have for all $(t,x) \in \R^+\times \R^N$,
 $$u(t,x)\le \psi(x.e+ct+t_0,x)$$ for some fixed $t_0$.
From  Shen and Zhang results, Theorem C in  \cite{shen-zhang}, since $c<\change{c_e^*}$ we have
$$\liminf_{t\to +\infty} \inf_{x.e+ct\ge 0} (u(x,t) -p(x)) =0.$$
Thus we get the following contradiction
\begin{align*}
0=\liminf_{t\to +\infty} \inf_{x.e+ct\ge 0} (u(x,t) -p(x)) &\le \liminf_{t\to +\infty} \inf_{x.e+ct\ge 0} (\psi(x.e+ct+t_0,x) -p(x)) \\
 &\le (\psi(t_0,x) -p(x))<0.
 \end{align*}
\fdem

\appendix
\section{Uniform estimates for solutions some regularized problems}
\label{appendix}

In this section we prove Proposition~\ref{pro aprox esti}. The estimates in this proposition divide naturally in 2 parts, one consisting in energy type estimates, and the other one are Schauder type estimates.

\medskip
\noindent{\bf Proof of Proposition~\ref{pro aprox esti} i).}
We proceed as in Lemma 2.5 in \cite{Berestycki-Hamel-Roques-II}.
Let us denote $\phi_{\kappa,\eps}$ the solution of \eqref{eq aprox-reg}. Then multiply equation \eqref{eq aprox-reg} by $\partial_s \psi_{\kappa,\eps}$ and integrate over $[-R,R]\times \q$ where $\q:=[0,1]^N$. Then it follows that

\begin{multline*}
c\int_{[-R,R]\times \q} |\partial_s \psi_{\kappa,\eps}|^2 = \changejer{\changejer{\kappa}}\int_{[-R,R]\times \q} \partial_s \psi_{\kappa,\eps}\partial_{ss}\psi_{\kappa,\eps} + \eps \int_{[-R,R]\times \q} \partial_s \psi_{\kappa,\eps}\Delta_x\psi_{\kappa,\eps} \\ +\int_{[-R,R]\times \q} \partial_s \psi_{\kappa,\eps} (M\psi_{\kappa,\eps} -\psi_{\kappa,\eps})+\int_{[-R,R]\times \q} \partial_s \psi_{\kappa,\eps} f(s,\psi_{\kappa,\eps})
\end{multline*}

Excepted the term $\mathcal{I}:=\int_{[-R,R]\times \q} \partial_s \psi_{\kappa,\eps} (M\psi_{\kappa,\eps} -\psi_{\kappa,\eps})$, all the term can be estimated  as in the proof of Lemma 2.5 in \cite{Berestycki-Hamel-Roques-II}, so we only deal with $\mathcal{I}$.


A simple computation shows that
$$ \int_{[-R,R]\times \q} \partial_s \psi_{\kappa,\eps} \psi_{\kappa,\eps} =\frac{1}{2}\int_{[-R,R]\times \q}\partial_s(\psi_{\kappa,\eps})^2=\frac{1}{2}\int_{ \q}[(\psi_{\kappa,\eps})^2]^R_{-R}.$$
So it remains to compute $$ I:=\int_{[-R,R]\times \q} \partial_s \psi_{\kappa,\eps} M\psi_{\kappa,\eps}.$$

Let us denote $\q_k:=k+\q$ where $k\in \Z^N$. With this notation, using the periodicity in $x$ of the function $\psi_{\kappa,\eps}$  we have
\begin{align*}
M\psi_{\kappa,\eps} 
&=\sum_{k\in\Z^N}\int_{k+\q}J(x-y)\psi_{\kappa,\eps}(s+(y-x).e,y)\, dy\\
&=\sum_{k\in\Z^N}\int_{\q}J(x-k-y)\psi_{\kappa,\eps}(s+(y-x).e+k.e,y)\, dy.
\end{align*}

Now using  integration by parts it follows that
 \begin{multline*}
 I=\int_{\q\times \q}\sum_{k\in\Z^N}J(x-y-k) [\psi_{\kappa,\eps}(s,x) \psi_{\kappa,\eps}(s+(y-x).e+k.e,y)]^R_{-R}\\-
 \int_{\q\times \q}  \sum_{k\in\Z^N}J(x-y-k) \int_{-R}^R\psi_{\kappa,\eps}(s,x) \partial_s\psi_{\kappa,\eps}(s+(y-x).e+k.e,y).
\end{multline*}

Let us make the change of variable $\tau=s+(y-x).e+k.e$ in the last term of the right hand side. Then we have
 \begin{multline*}
 \int_{\q\times \q} \int_{-R}^R \sum_{k\in\Z^N}J(x-y-k) \psi_{\kappa,\eps}(s,x) \partial_s\psi_{\kappa,\eps}(s+(y-x).e+k.e,y)\\=  \int_{\q\times \q} \sum_{k\in\Z^N}J(x-y-k)  \int_{-R+(y-x).e+k.e}^{R+(y-x).e+k.e}  \psi_{\kappa,\eps}(\tau + (x-y).e-k.e,x) \partial_s\psi_{\kappa,\eps}(\tau,y)
\end{multline*}

Let $R\to \infty$. Using that  $\psi_{\kappa,\eps} \to p_{\eps}$ respectively $0$ as $s\to \pm \infty$, $\psi_{\kappa,\eps}\ge 0, \partial_s\psi_{\kappa,\eps}\ge 0$ we obtain
\begin{align}
& \int_{\R\times \q} \partial_s \psi_{\kappa,\eps} \psi_{\kappa,\eps}=\frac{1}{2}\int_{\q}p_\eps^2 \label{ptw.eq.esti.ener1}
\end{align}
and
\begin{multline*}
\int_{\R\times \q}\partial_s\psi_{\kappa,\eps} M\psi_{\kappa,\eps} =\int_{\q\times \q}\sum_{k\in\Z^N}J(x-y-k) p_{\eps}(x) p_{\eps}(y) \\-\int_{\q\times \q}  \sum_{k\in\Z^N}J(x-y-k) \int_{-\infty}^{+\infty}\psi_{\kappa,\eps}(\tau+(x-y).e-k.e,x) \partial_s\psi_{\kappa,\eps}(\tau,y).
\end{multline*}
Going back to the definition of $M\psi_{\kappa,\eps}$ and using the symmetry of $J$ we can rewrite the above equality the following way
$$
\int_{\R\times \q}\partial_s\psi_{\kappa,\eps} M\psi_{\kappa,\eps} =\int_{\q}J* p_{\eps}(x) p_{\eps}(x)\,dx-\int_{\R\times \q}  M\psi_{\kappa,\eps}(\tau,y) \partial_\tau\psi_{\kappa,\eps}(\tau,y)\,d\tau dy.
$$
Thus we have
$$\int_{\R\times \q}\partial_s\psi_{\kappa,\eps} M\psi_{\kappa,\eps} =\frac{1}{2}\int_{\q}J* p_{\eps}(x) p_{\eps}(x)\,dx.$$

Set $\tilde \j(x,y):=\sum_{k\in\Z^N}J(x-y+k)$, the the above equality rewrites as follows

\begin{equation}\label{ptw.eq.esti.ener2}
\int_{\R\times \q}\partial_s\psi_{\kappa,\eps} M\psi_{\kappa,\eps} =\frac{1}{2}\int_{\q}\int_C \tilde \j(x,y)p_{\eps}(y) p_{\eps}(x)\,dydx
\end{equation}
Finally, combining \eqref{ptw.eq.esti.ener1} and \eqref{ptw.eq.esti.ener2}, we obtain
$$
\int_{\R\times \q}\partial_s\psi_{\kappa,\eps} (M\psi_{\kappa,\eps}-\psi_{\kappa,\eps})=-\frac{1}{4}\int_{\q\times \q}\tilde \j(x,y)(p_{\eps}(x)-p_{\eps}(y))^2\, dxdy.
$$

Hence,
$$ c\int_{\R\times \q}|\partial_s \psi_{\kappa,\eps}|^2 =-\frac{\eps}{2}\int_{\q}|\nabla_x p_{\eps}|^2 - \frac{1}{4}\int_{\q^2}\tilde \j(x,y)(p_{\eps}(x)-p_{\eps}(y))^2+\int_{\q}F(x,p_{\eps}) $$
which proves (i).
\qed

\medskip
\noindent{\bf Proof of Proposition~\ref{pro aprox esti} ii).}
Let $\K$ be a compact set of $\R\times\R^N$. Then since $\K$ is bounded, there exists $n\in\N$ and $R>0$ so that
$\K\subset (-R_0,R_0)\times n\tilde Q$ where $\tilde Q:=[-1,1]^N$.

Let us denote $\e(u)$ the following energy on the set of periodic function
$$\e(u):=-\frac{\eps}{2}\int_{\q}|\nabla_x u|^2 - \frac{1}{4}\int_{\q^2}\tilde \j(x,y)(u(x)-u(y))^2+\int_{\q}F(x,u). $$
From $(i)$, there exists  $R\in [R_0,R_0+1]$ so that
\begin{equation}\label{ptw.eq.R}
 c\int_{\q}|\partial_s\psi_{\kappa,\eps}|^2(R)\le \e(p_{\eps})
 \end{equation}

Let us now multiply \eqref{eq aprox-reg} by $\psi_{\kappa,\eps}$ and integrate over $(-R,R)\times \tilde Q$. Then we have
\begin{multline*}
\frac{c}{2}\int_{\tilde Q}[\psi^2_{\kappa,\eps}]^{R}_{-R}=\changejer{\kappa}\int_{\tilde Q}[\psi_{\kappa,\eps}\partial_s\psi_{\kappa,\eps}]^{R}_{-R}-\changejer{\kappa}\int_{(-R,R)\times \tilde Q}|\partial_s\psi_{\kappa,\eps}|^2-\eps\int_{(-R,R)\times \tilde Q}|\nabla_x\psi_{\kappa,\eps}|^2\\+\int_{(-R,R)\times \tilde Q}(M\psi_{\kappa,\eps}-\psi_{\kappa,\eps})\psi_{\kappa,\eps}+\int_{(-R,R)\times \tilde Q}f(x,\psi_{\kappa,\eps})\psi_{\kappa,\eps}
\end{multline*}

Therefore since $\psi_{\kappa,\eps}$ is uniformly bounded and periodic in $x$ we have,
$$
 \eps\int_{(-R,R)\times \tilde Q}|\nabla_x\psi_{\kappa,\eps}|^2=2\gamma(R)
$$
where
\begin{multline*}
\gamma(R):=-\frac{c}{2}\int_{\q}[\psi^2_{\kappa,\eps}]^{R}_{-R}-\changejer{\kappa}\int_{(-R,R)\times  \q}|\partial_s\psi_{\kappa,\eps}|^2 +\changejer{\kappa}\int_{ \q}[\psi_{\kappa,\eps}\partial_s\psi_{\kappa,\eps}]^{R}_{-R}\\+\int_{(-R,R)\times  \q}(M\psi_{\kappa,\eps}-\psi_{\kappa,\eps})\psi_{\kappa,\eps}+\int_{(-R,R)\times   \q}f(x,\psi_{\kappa,\eps})\psi_{\kappa,\eps}.
\end{multline*}

Since $0\le \psi_{\kappa,\eps}\le p_{\eps}$, $\partial_s\psi_{\kappa,\eps}\ge 0$ and  $f$ is uniformly bounded,  using Cauchy-Schwartz inequality   it follows that
\begin{align*}
\gamma(R)\le |c|\int_{\q}p_{\eps}^2 +\changejer{\kappa}\int_{\q}p^2_{\eps}\int_{\q}|\partial_s\psi_{\kappa,\eps}|^2(R,x)+2R\int_{\q}(J*p_{\eps}) p_{\eps}+2R \|f\|_{\infty}\int_{\q}p_{\eps}.
\end{align*}
Thus, since $c> 0$ by \eqref{ptw.eq.R} we have
\begin{align*}
\gamma(R)\le |c|\int_{\q}p_{\eps}^2 +\frac{\changejer{\kappa}\e(p_{\eps})}{|c|}\int_{\q}p^2_{\eps} +2R\int_{\q}(J*p_{\eps}) p_{\eps}+2R \|f\|_{\infty}\int_{\q}p_{\eps}.
\end{align*}
Hence the estimate (ii) follows by periodicity.
\medskip

\fdem

\medskip
The proof  of Proposition~\ref{pro aprox esti} iii) is based
on the next 2 lemmas.
The first one is a version of a result of \cite{berestycki-hamel-compde}, on gradient estimates for elliptic regularizations of semilinear parabolic equations.
The result in \cite{berestycki-hamel-compde} is based on Bernstein type estimates and is nonlinear in nature, while the estimates below have a linear character, and are based on a  technique of Brandt \cite{brandt1} (see also \cite{brandt2,knerr} and \cite{gilbarg-trudinger} Chap.\@ 3).

Given $R>0$ let
$$
Q_R = \{ (t,x) \in \R\times\R^N: |t|<R, \ |x_i|<R \quad
\forall i=1,\ldots,N\}.
$$
\begin{lemma}
\label{lemma grad est}
Suppose $u \in C^2(Q_R)$ satisfies
$$
\Delta_x u + \eps u_{tt} +  u_t = f(x,t)
\quad\text{in } Q_R
$$
where $0<\eps\le 1$, $f \in L^\infty(Q_R)$.
Then
\begin{align}
\label{bd partial der}
|\partial_{x_i} u(0,0)| \le
\left(
\frac{2(N+1)}{R} +2
\right)
\sup_{Q_R}|u| + \frac{R}{2} \sup_{Q_R}|f|
\end{align}
for all $i=1,\ldots,N$,
where $C$ is independent of $R$, $\eps$.
\end{lemma}
\noindent
{\bf Proof.}
Let us write $x = (x_1,x') \in \R^N$ with $x_1 \in \R$, $x'\in \R^{N-1}$.
Define
$$
\tilde Q = \{ (t,x_1,x') \in  \R \times \R \times \R^{N-1}
:  0<x_1<R, \ |x_i|<1 \quad \forall i=2,\ldots,N, \ |t|<1\}
$$
and
$$
v(t,x_1,x') =
\frac{1}{2} \left(
u(t,x_1,x') - u(t,-x_1,x')
\right)
$$
for $(t,x_1,x') \in \tilde Q$.
Let us write
$$
L v = \Delta_x v +\eps v_{tt} +  v_t .
$$
Then $L$ is an elliptic operator and satisfies the maximum principle. We have
$$
L v(t,x_1,x')=\frac{1}{2} \left(
f(t,x_1,x') - f(t,-x_1,x')
\right)
\quad\text{for } (t,x_1,x') \in \tilde Q
$$
and
$$
|v| \le \sup_{Q_R}|u|
\quad\text{in } \tilde Q.
$$
Let
$$
\bar v (t,x_1,x') = A x_1(R-x_1) + B ( x_1^2 + |x'|^2 + t^2)
$$
where
$$
B = \frac{1}{R^2} \sup_{Q_R}|u|
$$
and
$$
A = \frac{1}{2}
\left(
\sup_{Q_R}|f| + B(2N+2\eps +2R)
\right) .
$$
With these choices we see that
$$
|v|\le \bar v \quad \hbox{on } \partial \tilde Q.
$$
and
$$
L \bar v \le -\sup_{Q_R}|f|
\quad\text{in }  \tilde Q.
$$
By the maximum principle $\bar v - v \ge 0$ in $\tilde Q$. Similarly $\bar v+v\ge 0$ in $\tilde Q$ and therefore
$$
|v| \le \bar v \quad\text{in } \tilde Q.
$$
This implies
$$
|\partial_{ x_1} v(0,0)| \le A R
$$
and gives \eqref{bd partial der} for $i=1$.
The same proof replacing $x_1$ by any of the other variables $x_2,\ldots,x_n$ yields \eqref{bd partial der}.
\qed

\begin{lemma}
Suppose $u \in C^2(Q_2)$ satisfies
$$
u_t - \Delta_x u - \eps u_{tt}   = f(x,t)
\quad \text{in } Q_2
$$
where $\eps>0$ and $f\in L^\infty(Q_2)$.
Then for some $0<\alpha<1$ there is a constant $C$ independent of $\eps$ such that
$$
\sup_{|x|\le1, t_1,t_2 \in [-1,1]}
\frac{|u(x,t_1)-u(x,t_2)|}{|t_1-t_2|^\alpha} \le C
\left(
\sup_{Q_2}|f| + \sup_{Q_2} |u|
\right)
$$
\end{lemma}
\noindent
{\bf Proof.}
Let us write
$$
M = \sup_{Q_2}|f| + \sup_{Q_2} |u| .
$$
By Lemma~\ref{lemma grad est}
\begin{align}
\label{grad est}
\sup_{Q_{1}} |\nabla_x u| \le C M.
\end{align}
Let $\varphi \in C^1(\R^N)$ have support in the ball closed ball $\bar B_1 $ of $\R^N$. Multiplying the equation by $u \varphi $ and integrating in $B_2$ we find
$$
\frac{1}{2}\frac{d}{d t} \int_{B_2} u^2 \varphi \, d x
- \eps \frac{d}{d t} \int_{B_2} u u_t \varphi \, d x
+ \eps \int_{Q_1} u_t^2 \varphi \, d x
+ \int_{B_2} |\nabla u|^2 \varphi \, d x
+ \int_{B_2} \nabla u \nabla \varphi  u \, d x
$$
$$=
 \int_{B_2} f u \varphi \, d x.
$$
Integrating this from $t_0$ to $t_1$ with $-1 \le t_0 < t_1 \le 1$ and using \eqref{grad est}
gives
$$
 - \frac{\eps}{2}
\frac{d}{dt}\int_{B_2} u^2 \varphi \, d x \Big|_{t=t_1}
 + \frac{\eps}{2}
\frac{d}{dt}\int_{B_2} u^2 \varphi \, d x \Big|_{t=t_0}
+ \eps \int_{t_0}^{t_1 }\int_{Q_1} u_t^2 \varphi \, d x
= O(M^2)
$$
where $O(M^2)$ is uniform in $\eps$.
Integrate now with respect to $t_0 \in [1/2,2/3]$ and $t_1 \in [5/6,1]$. We obtain
$$
\eps \int_{1/2}^1  \int_{B_2} g(t)  u_t^2 \varphi \, d x \, d t = O(M^2)
$$
where $g(t)$ is a continuous function which is positive in $[1/2,1]$.
Therefore one can always select $t_0 \in [1/2,1]$, possibly depending on $\eps$, such that
\begin{align}
\label{int ut2 t0}
\eps\int_{B_2}   u_t(t_0)^2 \varphi \, d x = O(M^2) .
\end{align}
Now multiply the equation by $u_t \varphi$ and integrate in $B_2$, to obtain
$$
\int_{B_2} u_t^2 \varphi \, d x
-\frac{\eps}{2} \frac{d}{dt}\int_{B_2} u_t^2 \varphi \, d x
+ \frac{1}{2} \frac{d}{dt} \int_{B_2} |\nabla u|^2 \varphi
\, d x +
\int_{B_2}\nabla u \nabla \varphi u_t \, d x =
\frac{d}{d t} \int_{B_2} f u \varphi
$$
Integrating with respect to $t \in [-1/2,t_0]$ with $t_0$ as above yields
$$
\int_{-1/2}^{t_0}
\int_{B_2} u_t^2 \varphi \, d x \, d t
-
\frac{\eps}{2} \int_{B_2} u_t^2 \varphi \, d x
\Big|_{-1/2}^{t_0}
+
\frac{1}{2}  \int_{B_2} |\nabla u|^2 \varphi
\, d x \Big|_{-1/2}^{t_0}
+
\int_{B_2}\nabla u \nabla \varphi u_t \, d x
$$
$$=
\int_{B_2} f u \varphi\Big|_{-1/2}^{t_0}
$$
Using \eqref{grad est} and \eqref{int ut2 t0} we find
\begin{align}
\label{est a12}
\int_{-1/2}^{t_0}
\int_{B_2} u_t^2 \varphi \, d x \, d t
+
\int_{B_2}\nabla u \nabla \varphi u_t \, d x
= O(M^2)
\end{align}
But
$$
\left|
\int_{B_2}\nabla u \nabla \varphi u_t \, d x
\right|\le
\frac{1}{2}
\int_{B_2}|\nabla u|^2 \frac{|\nabla \varphi|^2}{\varphi} \, d x
+\frac{1}{2}
\int_{B_2} \varphi u_t^2 \, d x
$$
One can select a function $\varphi\ge 0$ with support the ball $|x|\le 1$ and positive in $|x|<1$ such that $ \frac{|\nabla \varphi|^2}{\varphi}  $ is bounded. So by \eqref{grad est}
$$
\left|
\int_{B_2}\nabla u \nabla \varphi u_t \, d x
\right|\le
O(M^2)+\frac{1}{2}
\int_{B_2} \varphi u_t^2 \, d x
$$
and integrating on $[-1/2,t_0]$ we have
$$
\left|
\int_{-1/2}^{t_0}
\int_{B_2}\nabla u \nabla \varphi u_t \, d x  \, d t
\right|\le
O(M^2)+\frac{1}{2}
\int_{-1/2}^{t_0}
\int_{B_2} \varphi u_t^2 \, d x
\, d t.
$$
This combined with \eqref{est a12}
gives
$$
\int_{-1/2}^{t_0}
\int_{B_2} \varphi u_t^2 \, d x
\, d t \le C M^2.
$$
We may further restrict $\varphi$ such that $\varphi\ge 1$ in the ball $|x|\le 1/2$ ad deduce
\begin{align}
\label{int est ut2}
\int_{Q_{1/2}} u_t^2 \, d x
\, d t \le C M^2
\end{align}
Let $t_1,t_2 \in [-1/4,1/4]$, with $t_1\le t_2$. Let $x\in\R^N$ with $|x|\le 1$. Then
$$
u(x,t_2)- u(x,t_1) = \int_{t_1}^{t_2} u_t(x,t) \, d t.
$$
Now integrate this with respect to $x$ in the ball of center $x_0$, $|x_0|\le 1/4$ and radius $r=(t_2-t_1)^{1/(2N)}$:
$$
\int_{B(x_0,r)} ( u(x,t_2)- u(x,t_1) ) \, d x
= \int_{t_1}^{t_2} \int_{B(x_0,r)}  u_t(x,t) \, d x \, d t.
$$
By the mean value theorem there is some $\bar x \in B(x_0,r)$ such that
$$
 u(\bar x,t_2)- u(\bar x,t_1)  =
\frac{C}{r^N}\int_{B(x_0,r)} ( u(x,t_2)- u(x,t_1) ) \, d x
$$
and therefore, using \eqref{int est ut2}
\begin{align*}
| u(\bar x,t_2)- u(\bar x,t_1)  |
& \le
\frac{C}{r^N}
\int_{t_1}^{t_2} \int_{B(x_0,r)} | u_t(x,t)| \, d x \, d t
\\
& \le \frac{C  (t_2-t_1)^{1/2}}{r^{N/2}}
\left(
\int_{t_1}^{t_2} \int_{B(x_0,r)} u_t(x,t)^2 \, d x \, d t
\right)^{1/2 }
\\
& \le C M (t_2-t_1)^{1/4} .
\end{align*}
Since \eqref{grad est} holds we deduce
$$
|u(x_0,t_2) - u(x_0,t_1) |\le C M (t_2-t_1)^{1/(2N)} .
$$
\qed

\noindent
{\bf Acknowledgments.}
J.D. was partially supported by   Fondecyt
1090167,
S.M. was partially supported by   Fondecyt  1090183, and both
acknoledge CAPDE-Anillo ACT-125 and Fondo Basal CMM.
This work is
also part of the MathAmSud NAPDE project (08MATH01)
and ECOS contract no. C09E06.

We are grateful to the reviewers of this article for very useful comments.


\begin{thebibliography}{00}

\bibitem{AB}
G. Alberti, G. Bellettini,
\emph{A nonlocal anisotropic model for phase transitions. {I}. {T}he optimal profile problem}, Math. Ann., 310, (1998), 3, 527--560.

\bibitem{BFRW}
P. W. Bates, P. C. Fife, X. Ren, X. Wang, \emph{Traveling Waves in a convolution model for phase transition.}  Arch. Rational Mech. Anal.  138  (1997),  no. 2, 105--136.

\bibitem{bates-zhao}
P.W. Bates, G. Zhao, Existence, uniqueness and stability of the stationary solution to a nonlocal evolution equation arising in population dispersal. J. Math. Anal. Appl. 332 (2007), no. 1, 428--440

\bibitem{berestycki-hamel-compde}
H. Berestycki, F. Hamel, {\em Gradient estimates for elliptic regularizations of semilinear parabolic and degenerate elliptic equations.} Comm. Partial Differential Equations 30 (2005), no. 1-3, 139--156.

\bibitem{berestycki-hamel-cpam}
H.~Berestycki, F.~Hamel, {\em Front propagation in periodic
excitable media.} Comm. Pure Appl. Math. 55 (2002), no. 8,
949--1032.

\bibitem{berestycki-hamel-contempmath}
H.~Berestycki, F.~Hamel, {\em Generalized travelling waves for reaction-diffusion equations.}  Perspectives in nonlinear partial differential equations,  101--123, Contemp. Math., 446, Amer. Math. Soc., Providence, RI, 2007.

\bibitem{Berestycki-Hamel-Nadin}
H. Berestycki, F. Hamel, G. Nadin, {\em Asymptotic spreading in heterogeneous diffusive excitable media.} J. Funct. Anal. 255 (2008), no. 9, 2146--2189

\bibitem{Berestycki-Hamel-Roques-I}
H.~Berestycki, F.~Hamel, L.~Roques, {\em Analysis of the
periodically fragmented environment model. I. Species persistence.}
J. Math. Biol. 51 (2005), no. 1, 75--113.

\bibitem{Berestycki-Hamel-Roques-II}
H.~Berestycki, F.~Hamel, L.~Roques, {\em Analysis of the
periodically fragmented environment model. II. Biological invasions
and pulsating travelling fronts.} J. Math. Pures Appl. (9) 84
(2005), no. 8, 1101--1146.

\bibitem{berestycki-larroutourou-lions}
H.~Berestycki, B.~Larrouturou, P.-L. ~Lions,
\emph{Multi-dimensional travelling-wave solutions of a flame propagation model}, Arch. Rational Mech. Anal., 111 (1990), 1, 33--49.

\bibitem{berestycki-nirenberg}
H.~Berestycki,  L. ~Nirenberg \emph{Travelling fronts in cylinders}, Ann. Inst. H. Poincar\'e Anal. Non Lin\'eaire, 9 (1992), 5, 497--572.

\bibitem{berestycki-nirenberg-varadhan}
H.~Berestycki, L.~Nirenberg, S.R.S.~Varadhan,  {\em The principal eigenvalue and maximum principle for second-order elliptic operators in general domains.} Comm. Pure Appl. Math. 47 (1994), 47--92.

\bibitem{brandt1}
A.~Brandt, {\em Interior estimates for second-order elliptic differential (or finite-difference) equations via the maximum principle.} Israel J. Math.  7 (1969), 95--121.

\bibitem{brandt2}
A.~Brandt, {\em Interior {S}chauder estimates for parabolic differential- (or difference-) equations via the maximum principle.} Israel J. Math. 7 (1969), 254--262.


\bibitem{browder}
F.E. Browder, {\em On the spectral theory of elliptic differential
operators. I.}  Math. Ann.  142  (1960/1961), 22--130.

\bibitem{CMS}
M. L. Cain,  B.G. Milligan, A.E. Strand,
\emph{Long-distance seed dispersal in plant populations}, Am. J. Bot., 87 (2000), 9, ,1217-1227.


\bibitem{carr-chmaj}
J.~Carr, A.~Chmaj,
{\em Uniqueness of travelling waves for nonlocal monostable equations.}
Proc. Amer. Math. Soc. 132 (2004), no. 8, 2433--2439.


\bibitem{Chen}
X. Chen, \emph{Existence, uniqueness, and asymptotic stability of traveling waves in nonlocal evolution equations}.   Adv. Differential Equations  2  (1997),  no. 1, 125--160.

\bibitem{Clark}
J.S. Clark,
\emph{Why Trees Migrate So Fast: Confronting Theory with Dispersal Biology and the Paleorecord}, The American Naturalist,
152 (1998), 2, 204-224.


\bibitem{Cov2}
{J. Coville},
{\em On uniqueness and monotonicity of solutions of non-local reaction diffusion equation},
{Ann. Mat. Pura Appl. (4)}
{185 (2006)},{3},{461--485.}


\bibitem{Cov4}
{J. Coville},
{\em Travelling fronts in asymmetric nonlocal reaction diffusion equation: The bistable and ignition case}, {Preprint of the CMM}.

\bibitem{Cov6}
{J. Coville},
{\em On a simple criterion for the existence of a principal eigenfunction of some nonlocal operators},
J. Differential Equations 249 (2010) 2921--2953.


\bibitem{CD1}
J. Coville, L. Dupaigne,
\emph{On a non-local reaction diffusion equation arising in population dynamics}.   Proc. Roy. Soc. Edinburgh Sect. A  137  (2007),  no. 4, 727--755.

\bibitem{cdm1}
J.~Coville, J.~D{\'a}vila, S.~Mart{\'{\i}}nez, {\em Existence and uniqueness of solutions to a nonlocal equation with monostable nonlinearity.} SIAM J. Math. Anal. { 39} (2008), no. 5, 1693--1709.

\bibitem{CDM2}
J.~Coville, J.~D{\'a}vila, S.~Mart{\'{\i}}nez, {\em Nonlocal anisotropic dispersal with monostable nonlinearity.} J. Differential Equations 244 (2008), no. 12, 3080--3118.

\bibitem{DGP}
A. De Masi, T.  Gobron, E.  Presutti, \emph{Travelling fronts in non-local evolution equations}, Arch. Rational Mech. Anal., 132 (1995), 2, 143--205.

\bibitem{DOPT1}
A. De Masi, E. Orlandi, E.  Presutti, L. Triolo, \emph{Glauber evolution with Kac potentials. I. Mesoscopic and macroscopic limits, interface dynamics}, Nonlinearity 7 (1994), 633-696.


\bibitem{DOPT2}
A. De Masi, E. Orlandi, E.  Presutti, L. Triolo,
\emph{Uniqueness and global stability of the instanton in nonlocal evolution equations}, Rend. Mat. Appl. (7) , 14 (1994), 4, 693--723.

\bibitem{DK}
  C. Deveaux,  E. Klein, \emph{Estimation de la dispersion de pollen \`a longue distance \`a l'echelle d'un paysage agicole : une approche exp\'erimentale}, Publication du Laboratoire Ecologie, Syst\`ematique et Evolution, 2004.



\bibitem{EPS}
D.E. Edmunds, A.J.B. Potter, C.A. Stuart, {\em Non-compact
positive operators}.  Proc. Roy. Soc. London Ser. A  328  (1972),
no. 1572, 67--81.

\bibitem{EMc}
G.B. Ermentrout, J. B. McLeod, \emph{Existence and uniqueness of travelling waves for a neural network}, Proc. Roy. Soc. Edinburgh Sect. A, 123 (1993), 3, 461--478.

\bibitem{F1}
P.C. Fife,
\emph{Mathematical aspects of reacting and diffusing systems}, Lecture Notes in Biomathematics, 28, Springer-Verlag, Berlin, 1979.


\bibitem{F2}
P.C. Fife,  \emph{An integrodifferential analog of semilinear parabolic {PDE}s}, Partial differential equations and applications, Lecture Notes in Pure and Appl. Math., 177, 137--145, Dekker, New York, 1996.

\bibitem{Freidlin} M.I. Freidlin, \emph{On wavefront propagation in periodic media}, Stochastic analysis and applications, Adv. Probab. Related Topics, 7, 147--166, Dekker, New York, 1984.


\bibitem{FreidlinGartner} M.I. Fre{\u\i}dlin,  Ju. Gertner, \emph{The propagation of concentration waves in periodic and random media}, Dokl. Akad. Nauk SSSR, 249, (1979), 3, 521--525.


\bibitem{GR1}
J. Garc{\'{\i}}a-Meli{\'a}n, J. D. Rossi, \emph{A logistic equation with refuge and nonlocal diffusion}, Commun. Pure Appl. Anal., 8, (2009), 6, 2037--2053.

\bibitem{gilbarg-trudinger}
D.~Gilbarg, N.S.~Trudinger,  Elliptic partial differential equations of second order. Reprint of the 1998 edition. Classics in Mathematics. Springer-Verlag, Berlin, 2001.


\bibitem{Hamel-Roques}
F. Hamel, L. Roques, \emph{Uniqueness and stability properties of monostable pulsating fronts.} Journal of the European Mathematical Society. In Press.

\bibitem{Heinze}
S. Heinze, \emph{Wave solutions to reaction-diffusion systems in perforated domains}, Z. Anal. Anwendungen, 20, (2001), 3, 661--676.


\bibitem{HPS}
S. Heinze,  G. Papanicolaou, A. Stevens, \emph{Variational principles for propagation speeds in inhomogeneous media}, SIAM J. Appl. Math., 62, (2001), 1, 129--148.

\bibitem{HZ}
W. Hudson, B. Zinner, \emph{Existence of traveling waves for reaction diffusion equations of {F}isher type in periodic media}, Boundary value problems for functional-differential equations, 187--199, World Sci. Publ., River Edge, NJ, 1995.


\bibitem{HMMV} V. Hutson, S.  Martinez, K.  Mischaikow, G.T.  Vickers, \emph{The evolution of dispersal}, J. Math. Biol., 47 (2003), 6, 483--517.


\bibitem{knerr}
B.F.~Knerr, {\em Parabolic interior {S}chauder estimates by the maximum principle.} Arch. Rational Mech. Anal.  75 (1980/81), no. 1, 51--58.


\bibitem{KPP} A.N. Kolmogorov, I. G.  Petrovsky, N. S.  Piskunov, \emph{\'Etude de l'\'equation de la diffusion avec croissance de la quantit\'e de mati\`ere et son application \`a un probl\`eme biologique}, Bulletin Universit\'e d'\'Etat \`a Moscow (Bjul. Moskowskogo Gos. Univ), S\'erie Internationale, (1937), Section A, 1--26.

\bibitem{KM}
M. Kot,  J. Medlock,
\emph{Spreading disease: integro-differential equations old and new}, Math. Biosci., 184 (2003), 2, 201--222.

\bibitem{krein-rutman}
M.G.~Krein, M.A.~Rutman, {\em Linear operators leaving invariant a cone in a Banach space.}  Uspehi Matem. Nauk (N. S.)  3,  (1948), no. 1(23), 3--95.

%

\bibitem{MNL}
H. Matano, K.I. Nakamura, B. Lou, \emph{Periodic traveling waves in a two-dimensional cylinder with saw-toothed boundary and their homogenization limit}, Netw. Heterog. Media, 1, (2006), 4, 537--568.

\bibitem{Mellet} A. Mellet, J.-M. Roquejoffre, Y. Sire. \emph{Generalized fronts for onedimensionnal
reaction-diffusion equations}. Discrete Contin. Dyn. Syst. A, 26, (2010), 1, 303--312.

\bibitem{M} J. D. Murray,  \emph{Mathematical biology}, Biomathematics, 19, Second Ed., Springer-Verlag, Berlin, 1993.

\bibitem{nadin}
G. Nadin, {\em Traveling fronts in space-time periodic media.} J. Math. Pures Appl. (9) 92 (2009), no. 3, 232--262.

\bibitem{nadin-rossi} G. Nadin, L. Rossi. \emph{Propagation phenomena for time heterogeneous
KPP reaction-diffusion equations}, preprint.

\bibitem{nolen-roquejoffre} J. Nolen, J.-M. Roquejoffre, L. Ryzhik, A. Zlatos.\emph{Existence and Nonexistence
of Fisher-KPP Transition Fronts}, preprint.

\bibitem{Nolen-Ryzhik}
 J. Nolen, L. Ryzhik, {\em Traveling waves in a one-dimensional heterogeneous medium.} Ann. Inst. H. Poincar\'e Anal. Non Lin\'eaire 26 (2009), no. 3, 1021--1047.

\bibitem{nussbaum}
R.D. Nussbaum, {\em The radius of the essential spectrum}. Duke
Math. J.  37  (1970), 473--478.

\bibitem{shen-die}
W.-X. Shen, {\em Traveling waves in time dependent bistable equations.}  Differential Integral Equations  19  (2006),  no. 3, 241--278.

\bibitem{shen-zhang}
W.-X. Shen,  A. Zhang,
{\em Spreading speeds for monostable equations with nonlocal dispersal in space periodic habitats.}
\change{J. Differential Equations 249 (2010), no. 4, 747--795.}

\change{\bibitem{shen-zhang2}
W.-X. Shen,  A. Zhang,
{\em Traveling wave solutions of spatially periodic nonlocal monostable equations}.}

\bibitem{SKT1} N. Shigesada,  K. Kawasaki, E. Teramoto, \emph{Traveling periodic waves in heterogeneous environments}, Theoret. Population Biol., 30 (1986), 1, 143--160.

\bibitem{SKT2} N. Shigesada, K. Kawasaki, E. Teramoto,\emph{The speeds of traveling frontal waves in heterogeneous environments}, Lecture Notes in Biomath., 71, 88--97, Springer, Berlin, 1987.


\bibitem{SSN}
F.M. Schurr, O.  Steinitz, R. Nathan,
\emph{Plant fecundity and seed dispersal in spatially heterogeneous environments: models, mechanisms and estimation}, J. Ecol., 96 (2008), 4, 628-641.




\bibitem{W2}
H.F. Weinberger,{\em On spreading speeds and travelling waves for growth and migration models in a periodic habitat.} J. Math. Biol. 45 (2002), no. 6, 511--548.
%

\bibitem{Xin1}
J.  Xin, \emph{Existence of planar flame fronts in convective-diffusive
periodic media}, Arch. Rational Mech. Anal., 121 (1992), 3, 205--233.

\bibitem{Xin2}
J. Xin, \emph{Existence and stability of travelling waves in periodic media governed by a bistable nonlinearity}, J. Dynam. Differential Equations, 3 (1991), 4, 541--573.

\bibitem{Xin3}
J. Xin, \emph{Front propagation in heterogeneous media}, SIAM Rev., 42, (2000), 2, 161--230.

\bibitem{zeidler}
E. Zeidler, Nonlinear functional analysis and its applications. I.
Fixed-point theorems. Springer-Verlag, New York, 1986.

\bibitem{Zlatos} A. Zlatos. \emph{ Generalized travelling waves in disordered media: Existence,
uniqueness, and stability}, preprint.
\end{thebibliography}
\end{document}